%% file: jet_0.tex
\documentclass[12pt]{book}
\usepackage{amsmath}
\usepackage{amssymb}
\pdfoutput=1 
\usepackage{ifpdf}

\title{\vskip -45pt\bf {Elements for a metric tangential calculus}}
\author{\Large{\bf Elisabeth Burroni}\\ \Large{\bf Jacques Penon}}
\date{\ }

\usepackage[french,english]{babel}

\input xypic
\usepackage[all]{xy}
\usepackage{picinpar, graphicx}

\input{macro}

\def\tang{\succ\!\!\!\prec}
\def\jf{_{_{\textup{\emph{free}}}}}

\begin{document}
\maketitle

\vfill\eject\null\vfill\eject

\vspace{1cm}
\centerline{\small\textbf{Abstract}}
The metric jets, introduced in the first chapter, generalize the jets (at order one) of Charles Ehresmann. In short, 
for a ``good'' map $f$ (said to be ``tangentiable'' at $a$), we define its metric jet tangent at $a$ (composed of all the maps which are locally lipschitzian at $a$ and tangent to $f$ at $a$) called the ``tangential'' of $f$ at $a$, and denoted T$f_a$ (the domain and codomain of $f$ being metric spaces).

Furthermore, guided by the heuristic example
of the metric jet T$f_a$, tangent  to
a map $f$ differentiable at $a$, which can be canonically represented by the unique continuous affine  map it contains, we will extend,
in the second chapter,
into a specific metric context, this property of representation of a metric jet.This yields a lot of relevant examples of such representations.
\par
\vspace{10mm}
\centerline{\textbf{INTRODUCTION}}
\vspace{2mm}

The first chapter contains a reference to two talks (given, one at the SIC [2], the other at the Calais conference [3]). Our aim being a deep thought inside the fundamentals of  differential calculus. 
Focussing ourselves on what is at the heart  of the notion of differential, it is the concept of  ``tangency'' which imposed itself on us in its great simplicity.
Now, amazingly, this concept of tangency can be formulated without resorting to the whole traditionnal structure of normed vector space 
(here on $\R$, what will be denoted n.v.s.):
see section 1.1. It is thus the more general structure  of metric space in which we have decided to work
from now on. Therefore, it is legitimate to ask ourselves if it is possible to construct a meaningful ``metric differential calculus''.
We will see that this aim has been essentially reached, even if, on the way,  it required  the help of an additionnal structure (the ``transmetric'' structure).

The first challenge we came accross was in the very formulation of a ``differential'' if we want to remain strictly in a metric context: what can we replace the continuous affine maps with, though they are essential to the definition of the classical differentials?
We cope with it with the help of jets 
(that we call metric jets in order to emphasize that
the metric structure is enough to define them)
which will  play the part of these continuous affine maps, 
willingly forgetting their algebraic feature.
At this point, we have to precise that, in the n.v.s. context, the continuous affine maps are not metric jets,
our metric jets being much more numerous than these maps (some of these jets possess at most one of these maps).
Actually, in the general metric context, our metric jets allow us to introduce a ``new differential'' for a map $f$
which admits a tangent at $a$ which is locally lipschitzian at $a$
(such a map being said to be ``tangentiable'' at $a$): 
it is a metric jet, tangent to $f$ at $a$, called the ``tangential'' of $f$ at $a$, and denoted T$f_a$.
So, our ``new differentiability'' (called ``tangentiability'') will extend 
the classical differentiability to new maps, traditionally pathologic.

As is well known, jets were first introduced by C.Ehresmann in 1951  [6],
in order to adapt Taylor's expansions to differential geometry;
more precisely, his infinitesimal jets (at order one) can be seen as equivalence classes of maps of class $C^1$ between differentiable varieties, under an equivalence relation of tangency.
The metric jets we are proposing here are more general
(being  equivalence classes of locally lipschitzian maps between metric spaces,
under an analogous equivalence relation of tangency).

Referring to the fact that, within the differential framework, we can define, with the help of the operator norm, a distance between two continuous affine maps, we will, more generally, 
construct a metric to evaluate the distance between two metric jets:
fixing a pair of points ($a$ and $a'$, respectively in the metric spaces $M$ and $M'$), the set  $\J$et$((M,a),(M',a'))$
of the metric jets from $(M,a)$ to $(M',a')$ can, itself, be equipped with a metric structure! 
Now, since these metric jets are ``tied'' to base-points,
it would be judicious to call them, from now on, ``linked metric jets''  (as one speaks of linked vectors, as opposed to free vectors), in opposition to the free metric jets quoted just below.

Noticing that our distance between linked metric jets
does not fit to speak of the distance between 
the tangentials of one map at two different points (these tangentials being linked metric jets which are tangent to this map at different points), we have added
a geometrical structure to our metric spaces, inspired by the translations of the usual vector space framework.  
These particular metric spaces are called ``transmetric''.
Inside these transmetric spaces, we can now speak of
``free metric jets'' (invariant by given ``translation jets'').
Now, if $M$ and $M'$ are such transmetric spaces, it is natural to consider the set $\J$et$_{\!{_{f\! r\! e\! e}}}(M,M')$
of the free metric jets from $M$ to $M'$; 
and this set $\J$et$_{\!{_{f\! r\! e\! e}}}(M,M')$ can then be equipped with a metric structure. 

Among these free metric jets, we find the
free metric jets of the form t$f_a$
``associated'' to the tangential T$f_a$
for an $f$ supposed to be tangentiable at $a$. 
Finally, when $f : M\lra M'$ is a tangentiable map (i.e. tangentiable at every point of $M$), we are now able to construct its tangential
t$f:M\lra \J$et$_{\!{_{f\! r\! e\! e}}}(M,M')$, whose domain and codomain are here metric spaces (since $M$ and $M'$ are transmetric spaces).
We thus obtain a new map on which, in the same way, we can apply the different techniques of the new theory, as, for instance, 
to study the continuity or the tangentiability of this tangential t$f$.

\vspace{2mm}

The second chapter takes up again and develops two talks (given, the first one at the SIC at Paris $[4]$, 
the second one at the internationnal conference at Calais $[5]$.
The (linked) metric jets being equivalence classes (of locally
lipschitzian maps), they are hard to ``manipulate'' concretely, and,
worse still, have a weak intuitive support.
In order to lighten their manipulation, we intend, 
referring to the vectorial case, to select a canonical element in each metric jet which could represent it and be compatible with the composition.
We cope with it in an algebraico-metric context: the one of 
``$\Sigma$-contracting'' metric space ($\Sigma$ being a ``valued monoid'' of actions which are compatible with the distance 
of the metric space, itself equipped with a
``central point'' denoted $\omega$);
the morphisms between such spaces being called ``$\Sigma$-homogeneous''.
All these $\Sigma$-homogeneous maps 
have a fundamental property in common, called the 
``$\Sigma$-uniqueness property'': two $\Sigma$-homogeneous
maps which are tangent at $\omega$ are equal.

Carrying on the analogy with the classical differential calculus, 
we are interested (still in the $\Sigma$-contracting world)
in the maps which possess, at the point $\omega$, a 
tangent which is lipschitzian and $\Sigma$-homogeneous.
In other words, maps which are tangentiable at $\omega$, and whose tangential at this point 
possesses a lipschitzian
$\Sigma$-homogeneous element which can represent it;
such a map is said
``$\Sigma$-contactable'' at $\omega$, 
the unique lipschitzian $\Sigma$-homogeneous element tangent to it at 
$\omega$ being its\break
``$\Sigma$-contact'' at $\omega$.
In many respects, the properties of
``$\Sigma$-contactibility''
are very similar to those of differentiability,
as, for instance, the search of extrema of a map 
(see section 2.6).

We will mainly be interested in two great special cases in the n.v.s. world (once we have realized that a map which is differentiable at $a$
is fortunately $\R$-contactable at $a$, its $\R$-contact being 
self-evidently the continuous affine map which is tangent to it at $a$;
although the inverse is false: see 2.5.12).
First,  the one where $\Sigma=\R_+$ which brings back the ``old''  interesting notion of maps which are
differentiable in the sense of Gateaux [7], which are very close to our $\R_+$-contactable maps (weaker than differentiable maps, and that we will call G-differentiable). 
Here, the lipschitzian $\R_+$-homogeneous maps can be easily seen as generalized cones (see 2.4.1).
On the way, we show that, in finite dimension, 
the continuity of the G-differential of a G-differentiable map $f$, forces this $f$ to be of class $C^1$ (see section 2.4)!
Secondly, the even weaker one 
where $\Sigma=\N'_r$ 
(i.e. $\N\cup\{\infty \}$ equipped with a
``valuation'' depending on $r$).
Here, the $\Sigma$-homogeneous maps are precisely the\break
$r$-fractal maps.
We let the reader go to the implications diagram of 2.5.10 and to 
2.5.11 to see all the links 
between these different contactibilities;
it gives also the opportunity to finally notice
that the notion of contactibility does not entirely exhaust the one of tangentiability,
since there exist maps which are tangentiable (at a central point) and $\Sigma$-contactable (at this point) for 
no $\Sigma$ ... see example 4) in 2.5.11.

\vspace{1mm}
For general definitions in category theory (for instance cartesian or enriched categories), see $[1]$.

\vspace{2mm}
\noindent\textsf{Acknowledgements}: 
We wish to thank our dear colleagues Christian Leruste for his constant help
with linguistical matters; and
Maryvonne Teissier for having designed for us the figures of this paper, with the help of the software Mathematica; and also
Vincent Padovani and Ren\' e Cori for 
their friendly useful help with computer matters.
We don't forget our ``brother in category'' Ren\' e Guitart: 
his opening to listening, his interest
in our work and his encouragements
have been a precious stimulating help.

\tableofcontents








 





 



\input{jet_I1}
\input{jet_I2}
\input{jet_I3}

\input{jet_II1}
\input{jet_II2}
\input{jet_II3}

\end{document}

%% file: macro.tex
\def \lra{\longrightarrow}

\def\Lra{\Longrightarrow}
\def\Lla{\Longleftarrow}
\def\Llra{\Longleftrightarrow}

\def\cqfd{\hfill$\square$\par}

\def\proof{\vspace{0mm}$\underline{\hbox{\sl Proof}}$ :\ }

\newtheorem{thm}{Théorème}[section]

\newtheorem{theo}[thm]{Theorem}
\newtheorem{prop}[thm]{Proposition}

\newtheorem{lema}[thm]{Lemma}

\newtheorem{cory}[thm]{Corollary}
\newtheorem{defi}[thm]{Definition}
\newtheorem{remk}[thm]{Remark}
\newtheorem{remks}[thm]{Remarks}

\newtheorem{exam}[thm]{Example}
\newtheorem{exams}[thm]{Examples}
\newtheorem{examsc}[thm]{Examples and counter-examples}
\newtheorem{cexam}[thm]{Counter-example}
\newtheorem{cexams}[thm]{Counter-examples}

\def\R{\mathbb R}
\def\L{\mathbb L}
\def\F{\mathbb F}
\def\J{\mathbb J}
\def\M{\mathbb M}
\def\C{\mathbb C}
\def\T{\mathbb T}
\def\I{\mathbb I}
\def\N{\mathbb N}
\def\Z{\mathbb Z}

\def\K{\mathbb K}
\def\Q{\mathbb Q}


%% file: jet_I1.tex
\chapter{A metric tangential calculus}

\section{The relation of tangency}
As said in the introduction, the relation of tangency
(at a given point $a$), which is at the heart of the differential calculus,
is an essentially metric notion, since it can merely be written
(the distance of any metric space being always be denoted $d$):
$\ \lim_{a\not= x\rightarrow a}\frac{d(f(x),g(x))}{d(x,a)}=0$
for two maps  $f,g:M\lra M'$ (with $M,M'$ metric spaces 
and $a\in M$); we should thus begin by such a definition.
However, this definition uses the term ``$\lim_{a\not= x\rightarrow a}$''
which makes sense solely for a point $a$, non isolated in $M$.
In order to avoid this snag, we have opted for a more general definition (see 1.1.2 and 1.1.3 below).

Yet, a first difficulty appears if we want  this relation of tangency to be compatible with the composition.
This will compel us to restrict the type of maps on which 
we will work. Once the relation of tangency is made composable, we will be able to define the concepts of metric jets and of tangentials (in the following sections).

\vspace{3mm}
First, let us begin by giving the right definition of the notion of tangency; it will arise from the following equivalences (easy to verify), where
$M$ and $M'$ are two metric spaces, $a$ a fixed point in $M$ and 
$f,g:M\lra M'$ two maps (a priori without any hypothesis). 

\vfill\eject
\begin{prop}
the  following properties are equivalent :

(i) $\ \forall\varepsilon>0\ \,\exists\eta>0\ \,\forall x\in M$\par
\qquad\qquad\quad$(d(x,a)\leq\eta\Lra d(f(x),g(x))\leq \varepsilon d(x,a))$,

(ii) $f(a)=g(a)$\  and the map $C:M\longrightarrow\mathbb R_+$, defined by:

\qquad\qquad\quad$C(a)=0$ and $\ \,\forall x\not=a$ $\ C(x)=
\frac{d(f(x),g(x))}{d(x,a)}$ 

\noindent is continuous at $a$,

(iii) there exists a map $c:M\lra\mathbb R_+$ which is continuous at $a$ and which 
verifies: 
$c(a)=0\ $ and $\ \forall x\in M\  
(d(f(x),g(x))=c(x)d(x,a))$,  

(iv) there exist a neighborhood $V$ of $a$ in $M$ and 
a map $c:V\lra\mathbb R_+$ which
is continuous at $a$ and which verifies: \par

\qquad\qquad\quad\ $c(a)=0\ $ and $\ \forall x\in V 
\ (d(f(x),g(x))
\leq c(x)d(x,a))$.  
\end{prop}

\begin{defi}
We say that $f$ and $g$ are 
 
\textup{tangent}
at $a$ (which is denoted by
$\, f\tang_a g$) 
if they verify anyone of the equivalent conditions of the above 1.1.1.
\end {defi}

\begin{remks}
\par\hfill

1) When $a$ is not an isolated point in $M$ (i.e 
$a\in \overline{M -\{a\}}$), we have: $f\tang_a g\ $ iff
$\ (f(a)=g(a)\ $ and $\ \lim_{a\not= x\rightarrow a}\frac{d(f(x),g(x))}{d(x,a)}=0)$.

2) When $a$ is an isolated point in $M$, we have $f\tang_a g$
for any\break 
$f$ and $g$ verifying $f(a)=g(a)$.
\end{remks}

\begin{prop}
The relation $\tang_a$ is an equivalence relation on the set of maps from $M$ to $M'$ ;
this relation $\tang_a$ is called the \textup{relation of tangency } at $a$.
\end{prop}

\begin{prop}
If $f\tang_a g$, then $f$ is continuous at $a$ iff $g$ is continuous at $a$.
\end{prop}

\vspace{3mm}
Let us now study the behaviour of the relation of tangency towards composition. Thus, we consider the following situation $(S)$ :
$$\xymatrix{
M_0\ar[rr]^{f_0}&&M_1\ar@<1ex>[rr]^{f_1}\ar@<-1ex>[rr]_{g_1}&&M_2\ar[rr]^{f_2}&&M_3
}$$
where $M_0$, $M_1$, $M_2$, $M_3$ are metric spaces with $a_o\in M_0$, $a_1=f_0(a_0)$,
$a_2=f_1(a_1)=g_1(a_1)$. Let us
ask ourselves under what conditions  do we have one of  the implications :
$\quad  f_1\tang_{a_1}g_1\quad\Lra\quad f_1.f_0\tang_{a_0}g_1.f_0$\qquad  and $\quad  f_1\tang_{a_1}g_1\quad\Lra \quad 
f_2.f_1\tang_{a_1}f_2.g_1$.

\begin{remk}
The above implications are not true in general, even if  the maps are continuous. Consider $\ M_0=M_1=M_2=M_3=\R$, with\break 
$f_0=f_2 :x\mapsto x^{1/3}$, $f_1:x\mapsto x^3$ and $g_1:x\mapsto 0$;
even though\break
$f_1\tang_0 g_1$, however 
$f_1.f_0\not\!\tang_0 g_1.f_0$ (since $f_1.f_0=Id_\R$ and $g_1.f_0=0$)
and\break
$f_2.f_1\not\!\tang_0 f_2.g_1$ (since 
$ f_2.f_1=Id_\R$ and $f_2.g_1=0$).
\end{remk}

We are now going to give sufficient conditions to make the above implications true.

\begin{defi}
Let $M$ and $M'$ be metric spaces, $f:M\lra M'$ a map, and $a\in M$; let also $k$ be a strictly positive real number.
We say that :

1)$\ f$ is \textup{locally} k-\textup{lipschitzian} at $a$ (in short $k$-$LL_a$) if there exists a neighborhood $V$ of $a$ in $M$ for which the restriction $f:V\lra M'$ is\break
$k$-lipschitzian.
$f$  \textup{locally\ lipschitzian} at $a$ (in short $LL_a$) means that there exists $k>0$ such that $f$ is $k$-$LL_a$.

2)$\ f$ is $k$-\textup{semi-lipschitzian} at $a$ (in short $k$-$SL_a$) if we have:\break
$\ \forall x\in M\ (d(f(x),f(a))
\leq kd(x,a)))$.
$f$ \textup{semi-lipschitzian} at $a$ (in short $SL_a$) means that there exists $k>0$ such that $f$ is $k$-$SL_a$.

3)$\ f$ is \textup{locally} k-\textup{semi-lipschitzian} at $a$ (in short $k$-$LSL_a$) if, on a neighborhood $V$ of $a$ in $M$, the restriction $f:V\lra M'$ is 
$k$-$SL_a$.\break
$f$ \textup{locally\ semi-lipschitzian} at $a$ (in short $LSL_a$) means that there exists $k>0$ such that $f$ is $k$-$LSL_a$.

Naturally, $f$ is $LL$ or $LSL$ will mean that $f$ is $LL_a$
or $LSL_a$ at every point $a\in M$.

\end{defi}

\begin{prop}
Let $M$ and $M'$ be metric spaces, $f:M\lra M'$ a map, and $a\in M$. We have the implications:  $\ f\ \, LL_a\ \Lra\ f\ \, LSL_a\ $\break
\noindent and $\ f\ \, LSL_a\ \Lra\ f\ $ continuous at $a$.
\end{prop}

\begin{remk}
The inverses of the previous implications are not true (see examples 1) and 3) of 1.3.9, with 1.1.12 and 1.3.2).
\end{remk}

\begin{prop}
In the above situation $(S)$, let us assume that\break
$f_1\tang_{a_1} g_1$ ; 
we then have the two following implications : 

1) $\ f_0\ \, LSL_{a_0}\ \Lra\ f_1.f_0\tang_{a_0} g_1.f_0$,

2) $\ f_2\ \, LL_{a_2}$ and $f_1,g_1$ continuous at $a_1\ \Lra\ f_2.f_1\tang_{a_1} f_2.g_1$.
\end{prop}

\proof
First, since $f_1\tang_{a_1} g_1$, we know that there exists a map $c:M_1\lra\R_+$ which is continuous at $a_1$ verifying: $c(a_1)=0$ and 
$\ \forall x_1\in M_1\ (d(f_1(x_1),g(x_1))=c(x_1)d(x_1,a_1))$.
Let us now come to the proposed implications:

1) Here, since $f_0$ is $LSL_{a_0}$, there exists a neighborhood $V_0$ of $a_0$ in $M_0$ and $k>0$ such that:
$\ \forall x_0\in V_0\ (d(f_0(x_0),f_0(a_0))\leq kd(x_0,a_0))$.
Then, for $x_0\in V_0$, we have:
$d(f_1.f_0(x_0),g_1.f_0(x_0))=c(f_0(x_0))d(f_0(x_0),f_0(a_0))\break
\leq kc(f_0(x_0))d(x_0,a_0)$.
From the fact that the map $c_0:V_0\lra\R_+:x_0\mapsto kc.f_0(x_0)$ is continuous at $a_0$ (since $f_0$ is $LSL_{a_0}$), and verifies $c_0(a_0)=kc(a_1)=0$,
we deduce that $f_1.f_0\tang_{a_0}g_1.f_0$ (by $(iv)$
of 1.1.1).

2) Since $f_2$ is $LL_{a_2}$, there exist a neighborhood $V$ of $a_2$ in $M_2$ and $k>0$ such that: 
$\ \forall x_2,y_2\in V\ (d(f_2(x_2),f_2(y_2))\leq kd(x_2,y_2))$.
The maps $f_1$ and $g_1$ being continuous at $a_1$,  $W=\buildrel{-1}\over{f_1}(V)\cap
\buildrel{-1}\over {g_1}(V)$ is a neighborhood of $a_1$ in $M_1$ and, for 
$x_1\in W$, we have $f_1(x_1),g_1(x_1)\in V$ and so, 
$d(f_2.f_1(x_1),f_2.g_1(x_1))\leq
k d(f_1(x_1),g_1(x_1))= kc(x_1)d(x_1,a_1)$ which implies
$f_2.f_1\tang_{a_1}f_2.g_1$ (since the map $c_1(x_1)=kc(x_1)$
is continuous at $a_1$
and verifies $c_1(a_1)=0$).
\cqfd

\begin{prop}
Let $M,M'$ be metric spaces, $a\in M$ and $f,g$ two maps $M\lra M'$ such that
$f\tang_a g$ ; we then have the equivalence:
$\ f\ \, LSL_a\ \iff\ g\ \, LSL_a$.
\end{prop}

\proof
If $f$ is $LSL_a$, there exist a neighborhood $V$ of $a$ and $k>0$ verifying:
$\ \forall x\in V\ d(f(x),f(a))\leq kd(x,a)$. Moreover, since $f\tang_a g$, there exists an $\eta>0\,$ ($B(a,\eta)$ denotes an open ball)
such that: \par
\noindent$\forall x\in M\
(d(x,a)\leq\eta\Lra d(f(x),g(x))\leq d(x,a)$). Thus, for\break
$x\in V\cap B(a,\eta)$, we have (using the fact that $g(a)=f(a)$):\par
\noindent$d(g(x),g(a))\leq d(g(x),f(x))+d(f(x),f(a))\leq
d(x,a)+kd(x,a)$; we thus obtain $d(g(x),g(a))\leq
(1+k)d(x,a)$, so that $g$ is $LSL_a$.
\cqfd

\begin{prop}
Let $E,E'$ be two n.v.s., $U$ an open subset of $E$,\break
$a\in U$ and $f:U\lra E'$ a map;
let us denote $L(E,E')$ the set of continuous linear maps from $E$ to $E'$.
We have  the implications :

1) $f$ differentiable at $a\ \Lra\ f\ \, LSL_a$,

2) $f\!$ differentiable and $\textup{d}f\! :\! U\!\lra\! L(E,E')$ continuous at $a$
$\Lra\ f\ \, LL_a$ (in particular,
$f\!$ of class $C^1$ $\ \Lra\ f\ \, LL$).
\end{prop}

\proof
1) If $f$ is differentiable at $a$, there exists a continuous affine map
$\alpha :E\lra E'$ 
such that $f\tang_a\alpha|_U$. But, since $\alpha|_U$ is $LSL_a$ (even lipschitzian!), $f$ is thus $LSL_a$ (by 1.1.11).

2) d$f$ being continuous at $a$, there exists $\eta>0$ such that
the open ball $B(a,\eta)\subset U$
and: $\forall x\in  B(a,\eta)$, we have $\|$d$f_x-$d$f_a||<1$, and thus 
$\|$d$f_x\|-\|$d$f_a\|\leq\|$d$f_x-$d$f_a\|<1$, which implies 
$\|$d$f_x\|\leq\|$d$f_a\|+1$. Then, by the mean value theorem 
($B(a,\eta)$ being convex), the restriction of $f$ to
$B(a,\eta)$ is $(\|$d$f_a\|+1)$-lipschitzian, so that $f$ is $LL_a$.
\cqfd

\begin{remks}
\par\hfill

1) The inverses of the previous implications are not true (see example 4) of 1.3.9).

2) We will give a generalization of the above implication 2) in
1.5.10.
\end{remks}

\begin{prop}
$\ $ Let $M_0,M_1,M_2$ be metric spaces ; \break
$f_0:M_0\lra M_1$,
$f_1:M_1\lra M_2$ two maps,
and $a_0\in M_0$, $a_1=f_0(a_0)$.
We have the implications:

1) $f_0\ \, LSL_{a_0}$ and $f_1\ \, LSL_{a_1}\ \Lra\ f_1.f_0\ \, LSL_{a_0}$.

2) $f_0\ \, LL_{a_0}$ and $f_1\ \, LL_{a_1}\ \Lra\ f_1.f_0\ \, LL_{a_0}$.
\end{prop}

\begin{prop}
Let $M_0,M_1,M_2$ be metric spaces; consider also maps
$f_0,g_0:M_0\lra M_1$,
$f_1,g_1:M_1\lra M_2$, and $a_0\in M_0$,\break
$a_1=f_0(a_0)=g_0(a_0)$.
We assume that $f_0\tang_{a_0} g_0$ and $f_1\tang_{a_1}g_1$ where $g_0$ is $LL_{a_0}$ 
and $g_1$ is $LL_{a_1}$ ; then
$f_1.f_0\tang_{a_0}g_1.g_0$. 
\end{prop}

\proof
Since $g_1$ is $LL_{a_1}$ and $f_0$, $g_0$ are continuous at $a_0$ (see 1.1.8,  1.1.5 or 1.1.11), we have 
$g_1.f_0\tang_{a_0}g_1.g_0$ by 1.1.10. On the other hand 
(still by 1.1.10), since $f_0$ is $LSL_{a_0}$
(see 1.1.8 and 1.1.11) and
$f_1\tang_{a_1} g_1$, we also have $f_1.f_0\tang_{a_0}g_1.f_0$.  Then, we use the fact that the relation of tangency is an equivalence relation.
\cqfd

\vfill\eject
\begin{remks}
\par\hfill

1) In 1.1.15, we could have weakened the hypothesis : 
$g_0\ \, LSL_{a_0}$ would have been enough. 

2) The theorem of composition of differentiable maps may be seen as an immediate consequence
of 1.1.15, since continuous affine maps are lipschitzian and are preserved by composition.
\end{remks}

\begin{prop}
Let $M,M_0,M_1$ be metric spaces, $a\in M$; consider also maps 
$f_0,g_0:M\lra M_0$ and $f_1,g_1:M\lra M_1$. We have the implication:
$\quad f_0\tang_a g_0$ and $f_1\tang_a g_1\  \Lra\ 
(f_0,f_1)\tang_a (g_0,g_1)$.
\end{prop}

\proof
We use here 1.1.1 and 1.1.2.
For each
$i\in\{0,1\}$, let $c_i:M\lra\R_+$ be a map which is continuous at $a$ and which verifies
$c_i(a)=0$ and $\ \forall x\in M\ (d(f_i(x),g_i(x))=c_i(x)d(x,a))$.
Let us denote $c(x)=\sup(c_0(x),c_1(x))$ and equip $M_0\times M_1$ with its product distance $d((x_0,x_1),(y_0,y_1))=\sup_id(x_i,y_i)$. The map $c$ being continuous at $a$ and verifying $\ c(a)=0$ and $\ \forall x\in M$
($d((f_0,f_1)(x),(g_0,g_1)(x))=\sup_i d(f_i(x),g_i(x))=c(x)d(x,a)$), we
thus obtain
$(f_0,f_1)\tang_a (g_0,g_1)$.
\cqfd

\begin{prop}
Let $M,M_0,M_1$ be metric spaces and $a\in M$; consider also
maps $f_0:M\lra M_0$ and $f_1:M\lra M_1$. Then:

1) $\ f_0,f_1\ \, LSL_a\ \Lra\ (f_0,f_1)\ \, LSL_a$,

2)  $\ f_0,f_1\ \, LL_a\ \Lra\ (f_0,f_1)\ \, LL_a$.
\end{prop}

\begin{remk}
$\ $ 1.1.18 implies that the categories whose objects are metric spaces and whose morphisms are maps which are $LSL$ (resp. $LL$) at a point, are cartesian categories.
\end{remk}

We conclude this section by giving a new link between $LSL$ maps and lipschitzian maps
(see 1.1.8)).
It is a generalization of the mean value theorem (here we weaken the hypothesis of being differentiable by the one of being $LSL$ ... see  1.1.12).

\begin{prop}
Let $M$ be a metric space, $[a,b]$ a compact interval of $\R$,
$k>0$ a fixed real number and
$f:[a,b]\lra M$ a continuous map which is $k$-$LSL_x$ for all $x$ in the open interval $]a,b[$.
Then we have $d(f(b),f(a))\leq k(b-a)$.
\end{prop}

For the proof (see below), we use the following well-known lemma :

\begin{lema}
Let $g:[a,b]\lra\R$ be a continuous map and $k$ a real number such that the following property is true:\par
\noindent$\ \forall x\in\,]a,b[\ \exists x'\in\,]a,b]\ (x'>x$ and $g(x')-g(x)\leq k(x'-x))$.\par
\noindent Then we have $g(b)-g(a)\leq k(b-a)$.
\end{lema}

Let us now go back to the proof of 1.1.20:

\proof
Let us denote $g(x)=d(f(x),f(a))$ ; then $g$ is continuous by composition.
Let $x\in\,]a,b[$. Since $f$ is $k$-$LSL_x$, there exists $\eta>0$ verifying:
$\ ]x-\eta,x+\eta[\,\subset\, ]a,b[\ $ and $\ \forall x'\in\,]x-\eta,x+\eta[$:\par
\noindent$d(f(x'),f(x))\leq k|x'-x|$.
So, if we take $x<x'<x+\eta$, we have $x'\leq b$ and 
$g(x')-g(x)\leq d(f(x'),f(x))
\leq k(x'-x)$.
So, thanks to the above lemma, we obtain 
$d(f(b),f(a))=g(b)-g(a)\leq k(b-a)$.
\cqfd

\begin{cory}
Let $M$ be a metric space, $[a,b]$ a compact interval of $\R$,
$F$ a finite subset of $]a,b[$; let also $f:[a,b]\lra M$ be a continuous map. Let us assume that, for all $x\in\,]a,b[-F$, the map $f$ is 
$k$-$LSL_x$ ;
then $d(f(b),f(a))\leq k(b-a)$.
\end{cory}

\proof
Let us write $F=\{a_1<\cdots <a_n\}$ with $a_0=a$ and $a_{n+1}=b$.\break
Since, for every $i\in\{0,...,n\}$, the restriction $f|_{[a_i,a_{i+1}]}$
is continuous and $k$-$LSL_x$ at every $x\in]a_i,a_{i+1}[$, 
we have, for every such $i$:  $d(f(a_{i+1}),f(a_i))\break
\leq k(a_{i+1}-a_i)$,
so that $d(f(b)-f(a))\leq
\sum_{i=0}^{n}d(f(a_{i+1}),f(a_i))\leq\break
k\sum_{i=0}^n(a_{i+1}-a_i)=k(b-a)$.
\cqfd

\begin{cory}
Let $E$ be a n.v.s., $U$ an open subset of $E$, $a,b\in U$ such that $[a,b]\subset U$ and $F$ a finite subset of $]a,b[$;
let also $M$ be a metric space and $f:U\lra M$ a continuous map.
Let us assume that, for all $x\in\,]a,b[-F$, the map $f$ is $k$-$LSL_x$ ;
then, we have again:\par
\noindent $d(f(b),f(a))\leq k\|b-a\|$.
\end{cory}

\proof
We just have to apply 1.1.22 to the composite
$f.\alpha$ (where $\alpha:[0,1]\lra U :t\mapsto a+t(b-a)$) which is 
$k\|b-a\|$-$LSL_t$ for all $t\in\, ]0,1[-\buildrel{-1}\over\alpha(F)$
(since $\buildrel{-1}\over\alpha(F)$ is still a finite subset, 
when $a\not= b$).
\cqfd
\section{Linked metric jets}

The linked metric jets (in short, here, the jets), which are merely equivalence classes for the relation of tangency, will play the part of the continuous affine maps of the classical differential calculus
(but here, without any algebraic properties); 
it seems natural,
in such a solely metric context,  
to define a distance between these jets, i.e to
equip the set of jets with a metric structure. Thanks to this metric structure, 
we will be able to enrich the category of jets, between pointed metric spaces, in the category $\M$et
(a well chosen category of metric spaces).

$M$ and $M'$ being metric spaces, with $a\in M,\ a'\in M'$, let us denote
$\L$L$((M,a),(M',a'))$ the set of maps $f:M\lra M'$ which are $LL_a$ and which verify $f(a)=a'$. These sets $\L$L$((M,a),(M',a'))$ are the ``Hom''
of a category, denoted $\L$L, whose objects are pointed metric spaces; this category $\L$L is a cartesian category (i.e. it
has a final object and finite products: see 1.1.14 and 1.1.18).
Now, since $\tang_a$ is an equivalence relation on $\L$L$((M,a),(M',a'))$, we  set  
$\J$et$((M,a),(M',a'))=\L$L$((M,a),(M',a'))/\tang_a$.
 
\begin{defi}
An element of $\J$\textup{et}$((M,a),(M',a'))$ is called a \textup{(linked metric) jet} from $(M,a)$ to $(M',a')$.
\end{defi}

Let $q:\L$L$((M,a),(M',a'))\lra\J$et$((M,a),(M',a'))$ be the canonical surjection. 
Thanks to 1.1.15, we can compose the jets: we have 
$q(g.f)=q(g).q(f)$ when $g$ and $f$ are composable.

So, we are now in a position to construct a category, denoted $\J$et, called the category of linked metric jets, whose:

- objects are pointed metric spaces $(M,a)$,

- morphisms $\varphi:(M,a)\lra(M',a')$ are jets (i.e. elements of\break
$\J$et$((M,a),(M',a'))$). In particular, the map $Id_M:(M,a)\lra(M,a)$ being $LL_a$, it provides a jet $q(Id_M):(M,a)\lra(M,a)$ denoted 
$Id_{(M,a)}$.

The previous canonical surjections extend to a functor 
$q:\L$L$\lra\J$et (constant on the objects) which makes 
$\J$et a quotient category of $\L$L.

\begin{prop}
The functor $q:\L$\textup L$\lra\J$\textup{et} creates a cartesian structure
on the category $\J$\textup{et} ($q$ being constant on the objects, it means that $\J$\textup{et} is cartesian and $q$ a strict morphism of cartesian categories). 
\end{prop}

\proof
The pair $(\{0\},0)$ is clearly a final object in $\L$L and in 
$\J$et. If $(M_0,a_0),(M_1,a_1)\in |\J$\textup{et}$|$, 
let us show that
$(M_0\times M_1,(a_0,a_1))$ is the wished product in the category 
$\J$et (where 
$M_0\times M_1$ is equipped with the product distance).
First, since the canonical projections\break
$p_i:M_0\times M_1\lra M_i$ 
are $1$-lipschitzian and verify $p_i(a_0,a_1)=a_i$, we can thus consider the  jets 
$\pi_i=q(p_i):
(M_0\times M_1,(a_0,a_1))\lra (M_i,a_i)$. On the other hand, if $(M,a)$ is another object in 
$\J$\textup{et}, and if, for every $i\in\{0,1\}$,
$\varphi_i:(M,a)\lra(M_i,a_i)$ is a jet, we can denote
$(\varphi_0,\varphi_1)=q(f_0,f_1)$, where $f_i\in\varphi_i$
(this definition is non-ambiguous thanks to 1.1.17); 
furthermore, $(\varphi_0,\varphi_1)$ is a jet
$(M,a)\lra(M_0\times M_1,(a_0,a_1))$ verifying $\pi_i.(\varphi_0,\varphi_1)=\varphi_i$.
Now, if
$\psi:(M,a)\lra(M_0\times M_1,(a_0,a_1))$ is 
another jet verifying $\pi_i.\psi=\varphi_i$ for
$i\in\{0,1\}$, then for $g\in\psi$ and $g_i=p_i.g$, we have:
$q(g_i)=\pi_i.q(g)=\varphi_i$, so that $g_i\tang_a f_i$.
Thus,
$g=(g_0,g_1)\tang_a(f_0,f_1)$ and
$\psi=q(g)=q(f_0,f_1)=(\varphi_0,\varphi_1)$.
Consequently the jets $\pi_i$ are the canonical projections in $\J$et.
\cqfd

\vspace{3mm}
We now give some particular morphisms in the category $\J$et which will be usefull further on.

\begin{prop}
Let  $\varphi:(M,a)\lra(M',a')$ be a morphism in $\J$\textup{et} and $f\in\varphi$.
If $f$ is locally ``anti-lipschitzian'' at $a$ (i.e if there exist $k>0$ and 
a neighborhood $V$ of $a$ on which we have $d(f(x),f(y))\geq kd(x,y)$), then $\varphi$ is a monomorphism in $\J$\textup{et}.
\end{prop}

\proof
Let $\varphi_1,\varphi_2:(N,b)\lra(M,a)$ be two parallel morphisms in $\J$et such that 
$\varphi.\varphi_1=\varphi.\varphi_2$.
Let us take $f_1\in\varphi_1$ and $f_2\in\varphi_2$; we thus have $f.f_1\tang_b f.f_2$,
so that there exists $c:N\lra\R_+$, continuous at $b$ verifying $c(b)=0$ and:
$\ \forall x\in N\ (d(f.f_1(x),f.f_2(x))=c(x)d(x,b))$. 
Then, for $x\in\ \buildrel{-1}\over{f_1}(V)\,\cap\buildrel{-1}\over{f_2}(V)$,
we have $d(f_1(x),f_2(x))\leq c_1(x)d(x,b)$ where $c_1(x)=(1/k)c(x)$.
Thus $f_1\tang_b f_2$, which provides $\varphi_1=\varphi_2$.
\cqfd

\begin{remk}
The jet of an isometric embedding is thus a monomorphism
(in particular, in the case of a metric subspace). 
\end{remk}

\begin{prop}
Let $M$ be a metric space, $V$ a neighborhood of\break
$a\in M$.
Let us set $j_a=q(j):(V,a)\lra(M,a)$ where $j:V\hookrightarrow M$ is the canonical injection.
Then, the jet $j_a$ is an isomorphism in $\J$\textup{et}.
\end{prop}

\proof
Let us consider $g:M\lra V$, the map defined by $g(x)=x$ if $x\in V$ and $g(x)=a$
if $x\notin V$ (clearly $g$ is $LL_a$ since $V$ is a  neighborhood of $a$); then $j_a^{-1}=q(g)$
(since $g.j=Id_V$ and $j.g\tang_a Id_M$).
\cqfd
\vspace{3 mm}
Time has now come to equip the category $\J$et$((M,a),(M',a'))$
with a metric structure (where $(M,a),(M',a')\in|\J$et$|$).

First, we define $d(f,g)$ for
$f,g\in \L$L$((M,a),(M',a'))$; at first, this $d$ will not be a distance on $\L$L$((M,a),(M',a'))$: see 1.2.6 below (but it will provide a ``true'' distance for the quotient 
$\J$et$((M,a),(M',a'))$).\break
For such $f,g$, we consider the map $C:M\lra\R_+$ 
defined by
$C(x)=\frac{d(f(x),g(x))}{d(x,a)}$ if $x\not= a$ and $C(a)=0$.
We notice that $C$ is bounded on a neighborhood of $a$: indeed, since $f$ and $g$ are $LL_a$, there exist a neighborhood $V$ of $a$ and a real number $k>0$ such that the restrictions
$f|_V$ and $g|_V$ are $k$-lipschitzian. Then, for $x\in V$, we have:
$d(f(x),g(x))\leq d(f(x),a')+d(a',g(x))
\leq d(f(x),f(a))+d(g(a),g(x))
\leq 2kd(x,a)$, so
that
$C(x)\leq2k$
for all $x\in V$.

Now, for each $r>0$, we set
$d^r(f,g)=\sup\{C(x)\,|\,x\in B'(a,r)\cap V\}$
(where $B'(a,r)$ is a closed ball; this definition does not depend on $V$ for small $r$).
The map $r\mapsto d^r(f,g)$ is increasing and positive, we can put:
\qquad\qquad{$d(f,g)=\lim_{r\rightarrow 0}d^r(f,g)=\inf_{r>0}d^r(f,g)$}.

\begin{prop}
Let $d: (\L\textup{L}((M,a),(M',a')))^2\lra\R_+$ be the map defined just above.
For each $f,g,h\in \L\textup{L}((M,a),(M',a'))$, this map\break
$d$ verifies the following properties:

1) $d(f,g)=d(g,f)$,

2)  $d(f,h)\leq d(f,g)+d(g,h)$,

3) $d(f,g)=0\iff f\tang_a g$.
\end{prop}

\proof
1) Arises from the fact that, for all $r>0$, we have\break
$d^r(f,g)=d^r(g,f)$.

2) Let $V$ be a neighborhood of $a$ and $k>0$ such that the restrictions $f|_V$, $g|_V$
and $h|_V$ are $k$-lipschitzian. Then, for $r>0$ and\break
$x\in B'(a,r)\cap V$, $\, x\not= a$, we have:
\noindent$\frac{d(f(x),h(x))}{d(x,a)}\leq\frac{d(f(x),g(x))}{d(x,a)}+\frac{d(g(x),h(x))}{d(x,a)}\leq
d^r(f,g)+d^r(g,h)$, which gives $d^r(f,h)\leq d^r(f,g)+d^r(g,h)$.
Doing $r\rightarrow 0$, we finally obtain
$d(f,h)\leq d(f,g)+d(g,h)$.

3) $d(f,g)=0\iff \forall\varepsilon>0\,\ \exists r>0\,\ d^r(f,g)<\varepsilon$\par
$\,\qquad\qquad\qquad\iff\forall\varepsilon>0\,\ \exists r>0\,\ \forall x\in B'(a,r)\quad C(x)<\varepsilon$\par
$\,\qquad\qquad\qquad\iff C$ is continuous at $a$\par
$\,\qquad\qquad\qquad \iff f\tang_a g$.
\cqfd

\begin{prop}
The map $d:(\L\textup{L}((M,a),(M',a')))^2\lra\R_+$, studied in 1.2.6, factors through the quotient, giving a ``true'' distance on $\J$et$((M,a),(M',a'))$, defined by 
$d(q(f),q(g))=d(f,g)$ for all\break
$f,g\in LL((M,a),(M',a'))$.
\end{prop}

It turns out that this distance provides for $\J$et a structure of category enriched in $\M$et,
where $\M$et is the cartesian category whose objects are the metric spaces and whose morphisms are the locally semi-lipschitzian maps 
(see 1.1.19).
But before proving this (see 1.2.19 and its corollary), we need to establish some technical properties about what we call  
the lipschitzian ratio of a jet (that we also need in section 1.4).

\begin{defi}
For $\varphi\in\J$\textup{et}$((M,a),(M',a'))$, we set $\rho(\varphi)=\inf K(\varphi)$, where 
$K(\varphi)=\{k>0\, |\, \exists f\in\varphi,\, f$ is $k$-$LL_a\}$.
It is this $\rho(\varphi)$ that we call  the
\textup{lipschitzian ratio} of $\varphi$.
Furthermore, we will say that $\varphi$ is $k$-\textup{bounded} if $\rho(\varphi)\leq k$.
\end{defi}

\begin{prop}
Let $(M_0,a_0)$, $(M_1,a_1)$, $(M_2,a_2)\in|\J$\textup{et}$|$; and also jets
$\varphi_0:(M_0,a_0)\lra(M_1,a_1)$, $\,\varphi_1:(M_1,a_1)\lra(M_2,a_2)$.
Then, $\rho(\varphi_1.\varphi_0)\leq\rho(\varphi_1)\rho(\varphi_0)$.
\end{prop}

\proof
Let $f_0\in\varphi_0$, $f_1\in\varphi_1$, and also $k_0,k_1>0$ such that, for each 
$i\in\{0,1\}$,
$f_i$ is $k_i$-$LL_{a_i}$.
Clearly, $f_1.f_0$ is $k_1k_0$-$LL_{a_0}$,
and $f_1.f_0\in\varphi_1.\varphi_0$; so that $\rho(\varphi_1.\varphi_0)\leq k_1k_0$.
Now, if we fix $k_1$, we can write $\frac{\rho(\varphi_1.\varphi_0)}{k_1}\leq k_0$ for each $k_0\in K(\varphi_0)$; so that  $\frac{\rho(\varphi_1.\varphi_0)}{k_1}\leq\rho(\varphi_0)$.
When $\rho(\varphi_0)\not=0$, we have also
$\frac{\rho(\varphi_1.\varphi_0)}{\rho(\varphi_0)}\leq k_1$ for each $k_1\in K(\varphi_1)$.
We thus obtain $\frac{\rho(\varphi_1.\varphi_0)}{\rho(\varphi_0)}\leq\rho(\varphi_1)$\break 
which implies
$\rho(\varphi_1.\varphi_0)\leq\rho(\varphi_1)\rho(\varphi_0)$.
If $\rho(\varphi_0)=0$, we have $\frac{\rho(\varphi_1.\varphi_0)}{k_1}=0$,
and then 
$\rho(\varphi_1.\varphi_0)=0$; which gives again
$\rho(\varphi_1.\varphi_0)\leq\rho(\varphi_1)\rho(\varphi_0)$.
\cqfd

\begin{prop}
For each $\varphi\in\J$\textup{et}$((M,a),(M',a'))$, we have\break
$d(\varphi,O_{aa'})\leq\rho(\varphi)$
(where $O_{aa'}=q(\widehat {a'})$, and $\widehat{a'}:M\lra M'$ is the constant map on $a'$).
\end{prop}

\proof
Let $f\in\varphi$, $\, k>0$ and $V$ be a neighborhood of $a$ for which $f|_V$ is $k$-lipschitzian.
Then, for each $r>0$ verifying $B'(a,r)\subset V$, we have:
\noindent$\forall x\in B'(a,r)\ \,(d(f(x),a')=d(f(x),f(a))\leq kd(x,a))$ which implies
$d^r(f,\widehat{a'})\leq k$. Thus, when $r\rightarrow 0$, we obtain
$d(\varphi,O_{aa'})=d(f,\widehat{a'})\leq k$;  
this being true for all $k\in K(\varphi)$, we finally obtain $d(\varphi,O_{aa'})\leq\rho(\varphi)$.
\cqfd

\vspace{3mm}
\begin{exams}
{}
\end{exams}

In all that follows $(M,a),(M',a'),(M_i,a_i)$ are objects of $\J$et, i.e pointed metric spaces. Here we calculate some lipschitzian ratios for what we could call \textit{good jets} $\varphi$, since, for them, the inequality of 1.2.10, becomes an equality! (refer to  1.2.21).
We will see in 2.5.3 that it's not always the case.
 
0)  We have $\rho(O_{aa'})=0$ (since, for all $\varepsilon>0$, we have
$\rho(O_{aa'})\leq\varepsilon$,
the constant map $\widehat{a'}:M\lra M'$ verifying $\widehat{a'}\in O_{aa'}$
and being $\varepsilon$-lipschitzian for all these $\varepsilon$).

1) For every jet $\varphi:(M,a)\lra(M',a')$, where $a$ or $a'$ are isolated (respectively in $M$ or $M'$), then $\rho(\varphi)=0$.
Indeed, if $a$ is isolated in $M$, then $\varphi=O_{aa'}$
(see 1.1.3)!
Now, if $a'$ is isolated in $M'$, it means that $\{a'\}$ is open in $M'$, so that, for all 
$f\in\varphi$, $V=\buildrel{-1}\over f(\{a'\})$ is a neighborhood of $a$; as $f|_V$ is constant on $a'$, it is $\varepsilon$-lipschitzian for all $\varepsilon>0$, so that
$\rho(\varphi)\leq\varepsilon$ for all these $\varepsilon$.

2) As in 1.2.2, we denote $\pi_i:(M_1,a_1)\times (M_2,a_2)\lra(M_i,a_i)$
the canonical projections in $\J$et; 
then $d(\pi_i,O_{aa_i})=\rho(\pi_i)=1$ (where $a=(a_1,a_2)$, 
with $a_i$ non isolated in $M_i$).
Indeed, as $p_i\in\pi_i=q(p_i)$ is\break
1-lipschitzian, we have 
$d(\pi_i,O_{aa_i})\leq\rho(\pi_i)\leq 1$.
If $d(\pi_i,O_{aa_i})<1$,  we would have $d(p_i,\widehat{a_i})=d(\pi_i,O_{aa_i})<1$,
so that (see just before 1.2.6) there would exist $r>0$ such that
$d^r(p_i,\widehat{a_i})<1$ and thus, for all $a\not=x\in B'(a,r)$,  
$\frac{d(p_i(x),a_i)}{d(x,a)}\leq d^r(p_i,\widehat{a_i})<1$  (where $x=(x_1,x_2))$.
This strict inequality will be contradicted choosing $x_i\not=a_i$ (which is possible since $a_i$ is not isolated in $M_i$) and $x_j=a_j$: we then obtain 
$\frac{d(p_i(x),a_i)}{d(x,a)}=1$ (since now 
$d(x,a)=d(x_i,a_i)$).

3) Let $M,M'$ be metric spaces, $f:M\lra M'$ an isometric embedding,
$a$ a non isolated point in $M$ and $a'=f(a)$); then
$d(q(f),O_{aa'})=\rho(q(f))=1$.
Indeed, as $f\in q(f)$, we have $\rho(q(f))\leq 1$.
Referring to 1.2.10, we have thus $d(q(f),O_{aa'})\leq 1$; 
if $d(q(f),O_{aa'})<1$, we would have $d(f,\widehat{a'})=d(q(f),O_{aa'})<1$, and thus (reasoning just like in example 2), with $f$ instead of $p_i$), there would exist $r>0$ such that $\frac{d(f(x),f(a))}{d(x,a)}\leq d^r(f,\widehat{a'})<1$ for all $a\not=x\in B'(a,r)$; which
contradicts the fact that $f$ is isometric.

Thus, $a$ being non isolated in $M$, $d(Id_{(M,a)},O_{aa})=
\rho(Id_{(M,a)})=1$ 
(where $Id_{(M,a)}=q(Id_M)$);
and $d(j_a,O_{aa})=\rho(j_a)=1$ and\break
 $d(j^{-1}_a,O_{aa})=\rho(j^{-1}_a)=1$, where $j_a=q(j)$, 
with $j:V\hookrightarrow M$ the canonical injection
($V$ is a neighborhood of $a$ in $M$ (see 1.2.5)).

\vspace{3mm}
\begin{prop}
Let $(M,a),(M',a')\in\J\textup{et}$.

1) Let us assume that there exists $\varphi\in\J\textup{et}((M,a),(M'a'))$
which is an isomorphism in $\J\textup{et}$. Then, $a$ is isolated in $M$ iff $a'$ is isolated in $M'$.

2) $(M,a)$ is a final object in $\J\textup{et}$ iff $a$ is isolated in $M$.
\end{prop}

\proof
1) Let us assume that $a$ is isolated in $M$; using 1.2.9 and 
examples 1) and 3) in 1.2.11, we obtain $\rho(Id_{(M',a')})=
\rho(\varphi.\varphi^{-1})\leq\rho(\varphi)\rho(\varphi^{-1})=0$,
so that $\rho(Id_{(M',a')})=0$, which implies that $a'$ is isolated in $M'$.

2) We use the fact that, in any category,  there exists a unique isomorphism between two final objects. 
Thus, if $(M,a)$ is a final object in $\J$et, we have an isomorphism jet $\,!:(M,a)\lra(\{0\},0)$, which implies (thanks to the above 1)) 
that $a$ is isolated in $M$.
Conversely, if $a$ is isolated in $M$, it forces the above jet $\, !\,$ to be an isomorphism,
its inverse being the jet $O_{0a}:(\{0\},0)\lra(M,a)$;
indeed, the equality $O_{0a}\,.\, !=Id_{(M,a)}$ comes from the fact that $a$ is an isolated point (see 1.1.3),
and the equality
$!\, .\, O_{0a}=Id_{(\{0\},0)}$ from the fact that $(\{0\},0)$ is a final object
in $\J$et.
\cqfd

\vfill\eject

\begin{prop}
$\ $ Let $(M,a),\ (M_1,a_1),\ (M_2,a_2)\in|\J\textup{et}|$, and\break
$\varphi_1:(M,a)\lra(M_1,a_1)$, $\ \varphi_2:(M,a)\lra(M_2,a_2)$ be two jets.\break
Then $\ \rho(\varphi_1,\varphi_2)=\sup_i(\rho(\varphi_i))$.
\end{prop}

\proof
Let $\varepsilon>0$; then, for $i\in\{1,2\}$, there exists 
$f_i\in\varphi_i$ which is
$k_i$-$LL_a$, where
$k_i=\rho(\varphi_i)+\varepsilon$. 
As $(f_1,f_2)$ is $k$-$LL_a$ with $k=\sup(k_1,k_2)$ and $(f_1,f_2)\in(\varphi_1,\varphi_2)$, we have
$\rho(\varphi_1,\varphi_2)\leq k=\sup_i(\rho(\varphi_i))+\varepsilon$ ... for all $\varepsilon>0$,
so that $\rho(\varphi_1,\varphi_2)\leq\sup_i(\rho(\varphi_i))$.
Besides, as $\varphi_i=\pi_i.(\varphi_1,\varphi_2)$, where the $\pi_i$ are the
canonical projections in $\J$et, we have, using 1.2.9, $\rho(\varphi_i)\leq
\rho(\pi_i)\rho(\varphi_1,\varphi_2)$, and thus $\rho(\varphi_i)\leq\rho(\varphi_1,\varphi_2)$, since $\rho(\pi_i)\leq 1$ (see 1.2.11); so that
$\sup_i(\rho(\varphi_i))\leq\rho(\varphi_1,\varphi_2)$.
\cqfd

\begin{cory}
For each $i\in\{1,2\}$, let $(M_i,a_i),(M'_i,a'_i)\in|\J\textup{et}|$, 
and $\psi_i:(M_i,a_i)\lra(M'_i,a'_i)$ be two jets. Then
$\rho(\psi_1\times\psi_2)\leq\sup_i\rho(\psi_i)$.
\end{cory}

\proof
Since $\psi_1\times \psi_2=(\psi_1.\pi_1,\psi_2.\pi_2)$, we have $\rho(\psi_1\times\psi_2)=\sup_i\rho(\psi_i.\pi_i)\leq\sup_i(\rho(\psi_i)\rho(\pi_i))\leq\sup_i\rho(\psi_i)$
(because $\rho(\pi_i)\leq 1$). 
\cqfd

\begin{theo}
Let us consider the following diagram in $\J$\textup{et}:\par 
$$\xymatrix{
(M_0,a_0)\ar@<1ex>[rr]^{\varphi_0}\ar@<-1ex>[rr]_{\psi_0}
&&(M_1,a_1)\ar@<1ex>[rr]^{\varphi_1}\ar@<-1ex>[rr]_{\psi_1}&&
(M_2,a_2)
}$$
We then have the inequalities:

1) $d(\psi_1.\psi_0,\varphi_1.\varphi_0)\leq
d(\psi_1,\varphi_1)d(\psi_0,O)+\rho(\varphi_1)d(\psi_0,\varphi_0)$
(where\break
\noindent $O=O_{a_0a_1}$: see 1.2.10).

2) $\, d(\psi_1.\psi_0,\varphi_1.\varphi_0)\leq d(\psi_1,\varphi_1)+d(\psi_0,\varphi_0)$
if $\psi_0$ and $\varphi_1$ are 1-bounded (see 1.2.8).
\end{theo}

\proof
1) Let $i\in \{0,1\}$; and $f_i\in\varphi_i$ and $g_i\in\psi_i$. Then, there exist $k_i,k'_i>0$ and $V_i$ a neighborhood of $a_i$ in $M_i$ such that the restrictions $f_i|_{V_i}$ and
$g_i|_{V_i}$ are respectively $k_i$-lipschitzian and $k'_i$-lipschitzian.

Referring to the definition of the distance on the jets given
in 1.2.7, 
we have to prove the inequality (where $\widehat{a_1}$ is the constant map on $a_1$):
$d(g_1.g_0,f_1.f_0)\leq
d(g_1,f_1)d(g_0,\widehat{a_1})+k_1d(g_0,f_0)$ where this $d$ has been defined just above 1.2.6 for the locally lipschitzian maps.
Let $R_0,R_1>0$ such that $B(a_0,R_0)\subset V_0\,\cap\!
\buildrel{-1}\over{f_0}\!(V_1)\,\cap\buildrel{-1}\over{g_0}(V_1)$ and $B(a_1,R_1)\subset V_1$,
and let us put $R=\inf(R_0,\frac{R_1}{k'_0})$.
If $0<r<R$ and  if $x\in B'(a_0,r)$, we have $k'_0r<k'_0R\leq k'_0\frac{R_1}{k'_0}=R_1$
and then
$d(g_0(x),a_1)=d(g_0(x),g_0(a_0))\leq k'_0d(x,a_0)<k'_0 r<R_1$, so that$g_0(x)\in B'(a_1,k'_0r)\subset B(a_1,R_1)\subset V_1$; we obtain the inequalities:

\noindent$d(g_1.g_0(x),f_1.f_0(x))\!\leq d(g_1.g_0(x),f_1.g_0(x))+d(f_1.g_0(x),f_1.f_0(x))$\par
\qquad\qquad\qquad\quad $\leq d(g_1.g_0(x),f_1.g_0(x))+k_1d(g_0(x),f_0(x))$\par
\qquad\qquad\qquad\quad $\leq d^{k'_0r}(g_1,f_1)d(g_0(x),a_1)+k_1d^r(g_0,f_0)d(x,a_0)$\par
\qquad\qquad\qquad\quad $\leq d^{k'_0r}(g_1,f_1)d^r(g_0,\widehat{a_1})d(x,a_0)+
k_1d^r(g_0,f_0)d(x,a_0)$.

Thus, dividing by  $d(x,a_0)$ when possible (i.e. if $x\not=a_0$), we obtain:
$\frac{d(g_1.g_0(x),f_1.f_0(x))}{d(x,a_0)}\leq d^{k'_0r}(g_1,f_1)d^r(g_0,\widehat{a_1})+
k_1d^r(g_0,f_0)$, which finally implies
$d^r(g_1.g_0,f_1.f_0)\leq d^{k'_0r}(g_1,f_1)d^r(g_0,\widehat{a_1})+
k_1d^r(g_0,f_0)$. It remains to do $r\rightarrow 0$ to obtain the foretold inequality.

Now, in order to obtain the wanted inequality 1), it suffices, when\break
$\psi_0\not=\varphi_0$
(i.e when $d(g_0,f_0)\not=0$: see 1.2.7), to write:\par
\noindent$\frac{d(\psi_1.\psi_0,\varphi_1.\varphi_0)-
d(\psi_1.\varphi_1)d(\psi_0,O)}{d(\psi_0,\varphi_0)}=
\frac{d(g_1.g_0,f_1.f_0)-
d(g_1.f_1)d(g_0,\widehat{a_1})}{d(g_0,f_0)}\leq k_1$; this being true for all $k_1\in K(\varphi_1)$, we obtain
$\frac{d(\psi_1.\psi_0,\varphi_1.\varphi_0)-
d(\psi_1.\varphi_1)d(\psi_0,O)}{d(\psi_0,\varphi_0)}\leq\rho(\varphi_1)$.
We have thus proved the inequality 1) in the case where 
$\psi_0\not=\varphi_0$.
But this inequality is still true when $\psi_0=\varphi_0$, since, in this case, we have 
$d(g_1.g_0,f_1.f_0)\leq d(g_1,f_1)d(g_0,\widehat{a_1})$, which implies the inequality:\par
\noindent$d(\psi_1.\psi_0,\varphi_1.\psi_0)\leq d(\psi_1,\varphi_1)d(\psi_0,O)$. 

2) We apply 1.2.10 in 1).
\cqfd

\begin{cory}
The sets of jets being equipped with their distance, the maps
$\J$\textup{et}$((M_0,a_0),(M_1,a_1))\lra\J$\textup{et}$((M_0,a_0),(M_2,a_2)):\psi\mapsto\varphi_1.\psi$ 
and
$\J$\textup{et}$((M_1,a_1),(M_2,a_2))\lra\J$\textup{et}$((M_0,a_0),(M_2,a_2)):\psi\mapsto\psi.\varphi_0$
are respectively $\rho(\varphi_1)$-lipschitzian and $d(\varphi_0,O)$-lipschitzian (where $\varphi_0$ and $\varphi_1$ are jets as in 1.2.15).
\end{cory}

\proof
We just have to use the inequality 1) proved in 1.2.15 with 
$\psi_1=\varphi_1$
in the first case and $\psi_0=\varphi_0$ in the second one.\cqfd

\begin{cory}
$\psi_0,\psi_1,O\!$ being jets as 
in 1.2.15 (with $O=O_{a_0a_2},
O_{a_1a_2}$ or 
$O_{a_0a_1}$), we have the inequality:
$d(\psi_1.\psi_0,O)\leq 
d(\psi_1,O)d(\psi_0,O)$.
\end{cory}

\proof
We just put $\varphi_0=O_{a_0a_1}$ and $\varphi_1=O_{a_1a_2}$ in 1.2.15 (since $\rho(O)=0$, according to 1.2.11).
\cqfd

\begin{remk}
The inequalities obtained in 1.2.9 and 1.2.17 are both generalisations of the well-known inequality
$\|l_1.l_0\|\leq\|l_1\|\,\|l_0\|$ for composable continuous linear maps (see 1.2.28 below).
\end{remk}

\vspace{1mm}
\begin{prop}
The composition of jets:\par
\noindent$\J$\textup{et}$((M_0,a_0),\!(M_1,a_1))\!\times\!\J$\textup{et}
$((M_1, a_1),\!(M_2,a_2))
\!\lra^{^{\!\!\!\!\!\!\!\!\!\!\!\!\!\!\!{comp}}}\!
\J$\textup{et}$((M_0,a_0),\!(M_2,a_2))$ is $LSL$.
\end{prop}

\proof
Let us set $\J$et$_{ij}=\J$et$((M_i,a_i),(M_j,a_j))$ for $i,j\in\{0,1,2\}$; and let us fix
$(\varphi_0,\varphi_1)\in \J$et$_{01}\times \J$et$_{12}$, this set
being equipped with the product distance.
If $B=B((\varphi_0,\varphi_1),1)$ is a unit open ball in\break
$\J$et$_{01}\times \J$et$_{12}$, then, 
for $(\psi_0,\psi_1)\in
B$, we can write (where here\break
$O=O_{a_0a_1}$):
$d(\psi_0,O)\leq d(\psi_0,\varphi_0)+d(\varphi_0,O)\leq 1+d(\varphi_0,O)$
and thus also, using 1.2.15: 
$d(\psi_1.\psi_0,\varphi_1.\varphi_0)\leq
d(\psi_1,\varphi_1)d(\psi_0,O)+\rho(\varphi_1)d(\psi_0,\varphi_0)$\par
\noindent $\leq d((\psi_0,\psi_1),(\varphi_0,\varphi_1))(d(\psi_0,O)+\rho(\varphi_1))$\par
\noindent $\leq
d((\psi_0,\psi_1),(\varphi_0,\varphi_1))(d(\varphi_0,O)+\rho(\varphi_1)+1)$.\par
\noindent The map $comp$ is thus 
$(d(\varphi_0,O)+\rho(\varphi_1)+1)$-semi-lipschitzian on $B$.
\cqfd

\begin{cory}
The category $\J$\textup{et} can thus be enriched in the cartesian category $\M$\textup{et} previously quoted just before 1.2.8.
\end{cory}

\begin{prop}
First, refer to 1.2.11 about good jets;
and let us denote
$Jeg((M_1,a_1),(M_2,a_2))$ the following set:\par
\noindent $\{\varphi\in\J$\textup{et}
$((M_1,a_1),(M_2,a_2))\,|\,d(\varphi,O_{a_1a_2})=\rho(\varphi)\}$.
Then the restriction
$\J$\textup{et}$((M_0,a_0),\!(M_1,a_1))\times Jeg((M_1,a_1),\!(M_2,a_2))
\!\lra^{^{\!\!\!\!\!\!\!\!\!\!\!\!\!\!\!{comp}}}
\!\J$\textup{et}$((M_0,a_0),\!(M_2,a_2))$ is $LL$.
\end{prop}

\proof
We have to consider the restriction
$\J$et$_{01}\times Jeg_{12}\lra^{^{\!\!\!\!\!\!\!\!\!\!\!\!\!\!\!{comp}}} \J$et$_{02}$ where
the {$\J$et$_{ij}$} are as in 1.2.19 and $Jeg_{12}=
Jeg((M_1,a_1),\!(M_2,a_2))$.
Let $(\varphi_0,\varphi_1)\in \J$et$_{01}\times Jeg_{12}$ and
$B=B((\varphi_0,\varphi_1),1)$ for the product distance in $\J$et$_{01}\times Jeg_{12}$.
Now, if $(\psi_0,\psi_1),(\psi'_0,\psi'_1)\in B$, we write, using the inequality 1) of 1.2.15 (with
$O=O_{a_0a_1}$ or $O_{a_1a_2}$):\par
\noindent $d(\psi_1.\psi_0,\psi'_1.\psi'_0)\leq
d(\psi_1,\psi'_1)d(\psi_0,O)+\rho(\psi'_1)d(\psi_0,\psi'_0)=$\par
\noindent$d(\psi_1,\psi'_1)d(\psi_0,O)+d(\psi'_1,O)
d(\psi_0,\psi'_0)\leq$\par
\noindent$d(\psi_1,\psi'_1)(d(\psi_0,\varphi_0)+d(\varphi_0,O))+
d(\psi_0,\psi'_0)(d(\psi'_1,\varphi_1)+d(\varphi_1,O))\leq$\par
\noindent$d(\psi_1,\psi'_1)(1+d(\varphi_0,O))+d(\psi_0,\psi'_0)(1+d(\varphi_1,O))\leq$\par
\noindent$d((\psi_0,\psi_1),(\psi'_0,\psi'_1))(2+d(\varphi_0,O)+d(\varphi_1,O))$.\par
We thus have obtained that the restriction $comp|_B$ is $k$-lipschitzian
with $k=2+d(\varphi_0,O)+d(\varphi_1,O)$, which ends the proof.
\cqfd

\begin{prop}
$(M,a),(M_0,a_0),(M_1,a_1)$ being objects in the category $\J$\textup{et}, the following canonical map $can$ is an isometry:
$$\xymatrix{
\J\textup{et}((M,a),(M_0,a_0)\times (M_1,a_1))\ar[d]^{can}\\
\J\textup{et}((M,a),(M_0,a_0))\times\J\textup{et}
((M,a),(M_1,a_1))
}$$
\end{prop}

\proof
Consider a pair of jets
$\varphi,\psi:(M,a)\lra(M_0,a_0)\times(M_1,a_1)$;
we thus must prove that we have the equality:\par
\noindent$d(\varphi,\psi)=\sup(d(\pi_0.\varphi,\pi_0.\psi),
d(\pi_1.\varphi,\pi_1.\psi))$; or else, if $f\in\varphi$,
$g\in\psi$ and $f_i=p_i.f$, $g_i=p_i.g$ (see 1.2.2 for the definitions of $p_i$ and $\pi_i$), to prove the equality:
$d(f,g)=sup_id(f_i,g_i)$ (for the definition of the distance
between jets, see 1.2.7). It suffices to use the equality
$d^r(f,g)=\sup_id^r(f_i,g_i)$, which is true for small $r$, and to
do $r\rightarrow 0$.
\cqfd

\begin{remk}
As isometries are  isomorphisms in $\M$\textup{et}, it means that
$\J$\textup{et} is an enriched cartesian category.
\end{remk}

\begin{prop}
Let $M,M'$ be metric spaces, $V,V'$ be two neighbor-\break
hoods,
respectively of $a\in M$ and $a'\in M'$. Then, the map:\par
\noindent$\Gamma:\J$\textup{et}$((V,a),(V',a'))\lra
\J$\textup{et}$((M,a),(M',a')):\varphi\mapsto j'_{a'}.\varphi.j_a^{-1}$\break
is an isometry (where $j_a$ and $j'_{a'}$ have been defined in 1.2.5).
\end{prop}

\proof
We can deduce from 1.2.16 and 1.2.10  that, for
$\varphi,\varphi'\in 
\J$et$((V,a),(V',a'))$, we have:
\noindent$d(\Gamma(\varphi),\Gamma(\varphi'))=
d( j'_{a'}.\varphi.j_a^{-1},
j'_{a'}.\varphi'.j_a^{-1})$
$\leq \rho(j'_{a'})d(\varphi.j_a^{-1},
\varphi'.j_a^{-1})$
$\leq \rho(j'_{a'})d(\varphi,\varphi')
d(j_a^{-1},O)$
$\leq \rho(j'_{a'})d(\varphi,\varphi')
\rho(j_a^{-1})$
$\leq d(\varphi,\varphi')$,
\noindent the last inequality resulting from the fact that
$\rho(j_a)$ and $\rho(j_a^{-1})$ are 1-bounded
(see 1.2.11);
same things for $j'_{a'}$.
But $\Gamma$ is bijective with $\Gamma^{-1}(\psi)={j'}_{\!a'}^{-1}.\psi.j_a$,
so that we have also $d(\Gamma^{-1}(\psi),\Gamma^{-1}(\psi'))\leq
d(\psi,\psi')$. Finally, setting $\psi=\Gamma(\varphi)$ and 
$\psi'=\Gamma(\varphi')$, we obtain the
equality:
$d(\Gamma(\varphi),\Gamma(\varphi'))
=d(\varphi,\varphi')$.
\cqfd

\vspace{3mm}
We conclude this section with a come back to vectorial considerations.
\begin{prop}
Let $M$ be a metric space (with $a\in M$) and $E$ a n.v.s.. 
Then, we can canonically equip the sets
$\L\textup{L}((M,a),(E,0))$ and $\J\textup{et}((M,a),(E,0))$
with vectorial space structures, making linear the canonical surjection $q:\L\textup{L}((M,a),(E,0))\lra\J\textup{et}((M,a),(E,0))$.\break
Besides, the distance on $\J\textup{et}((M,a),(E,0))$, defined in 1.2.7, derives from a norm (providing a structure of n.v.s. on 
$\J\textup{et}((M,a),(E,0))$).
\end{prop}

\proof
The vectorial structure on
$\L$L$((M,a),(E,0))$ results, as usual, from the one
of $E$. In
$\widehat E=\J\textup{et}((M,a),(E,0))$, the zero is $O_{a0}$;
as for the vectorial operations, we can set, for $\varphi,\psi\in \widehat E$, $\varphi+\psi=q(\sigma).(\varphi,\psi)$ (where 
$\sigma$ is the addition of $E$) and, for $\lambda\in\R$, $\lambda\varphi=q(m_\lambda).\varphi$ (where
$m_\lambda:E\lra E:x\mapsto\lambda x$).
The verification of the axioms of vectorial space and the linearity of $q$ arises from the fact that $q$ is a morphism of cartesian categories (see 1.2.2).

Now, we set $\|\varphi\|=d(\varphi, O_{a0})$ for 
$\varphi\in\widehat E$ and we verify that 
$\|\varphi-\psi\|=d(\varphi,\psi)$ and $\|\lambda\varphi\|\leq 
|\lambda|\,\|\varphi\|$ (which is enough to have a norm);
actually, for that, it suffices to prove that, for $f\in\varphi$, $g\in\psi$ and $r$ small enough, we have $d^r(f-g,\widehat 0)=
d^r(f,g)$ and $d^r(\lambda f,\widehat 0)\leq|\lambda| d^r(f,\widehat 0)$.
\cqfd

\begin{prop}
Let $M,M'$ be metric spaces, with $a\in M$ and\break
$a'\in M'$; let also
$\varphi\in \J\textup{et}((M',a'),(M,a))$ and $E$ a n.v.s..
Then, the map $\widetilde\varphi: \J\textup{et}((M,a),(E,0))\lra \J\textup{et}((M',a'),(E,0)):
\psi\mapsto \psi.\varphi$ is linear and continuous.
\end{prop}

\proof
For the linearity of $\widetilde\varphi$, we must verify that
$(\psi_1+\psi_2).\varphi=\psi_1.\varphi+\psi_2.\varphi$ and $(\lambda\psi).\varphi=\lambda(\psi.\varphi)$;
but, this arises from the fact that $q$ is functorial and also linear on the ``Hom'' (see 1.2.25).

Now, the continuity of $\widetilde\varphi$ arises from the fact that 
$\J$et is enriched in $\M$et;
we can also use 1.2.16).
\cqfd

\begin{prop}
$E$ and $E'$ being n.v.s., the canonical map :\par
\noindent$j:L(E,E')\lra \J$\textup{et}$((E,0),(E',0)):l\mapsto q(l)$ is a linear isometric
embedding ($L(E,E')$ being equipped with the norm
$\|l\|=\sup_{\|x\|\leq 1}\|l(x)\|$).
\end{prop}

\proof
The linearity of $j$ comes from the fact that 
$L(E,E')$ is a vecto-\break
rial subspace of $\L$L$((E,0),(E',0))$ on which
$q$ is linear (see 1.2.25).
 
Now, to prove that $j$ is isometric, we have just to verify that, for\break 
$l\in \L(E,E')$ and $r>0$, we have $\|l\|=d^r(l,0)$ 
(doing $r\rightarrow 0$, we obtain then $\|l\|=d(l,0)=
d(q(l),q(0))=\|q(l)\|=\|j(l)\|$). In order to obtain this equality, 
we fix $r>0$; first, for $x\in B'(0,r)$, 
we have $\|l(x)\|\leq\|l\|\, \|x\|$ which implies
$d^r(l,0)\leq\|l\|$.
Secondly, if $0<\|x\|\leq 1$ and
$0<r'<r$, we have $\|r'x\|\leq r'<r$
and then $\frac{\|l(x)\|}{\|x\|}=
\frac{\|l(r'x)\|}{\|r'x\|}
\leq d^r(l,0)$, which implies
$\|l(x)\|\leq d^r(l,0)\|x\|\leq d^r(l,0)$; so that
$\|l\|\leq d^r(l,0)$.
\cqfd

\begin{cory}
If $\ l:E\lra E'$ is a continuous linear map, we have
$q(l)\in Jeg((E,0),(E',0))$; more precisely, we have the equalities:\par
\noindent$d(q(l),O)=\rho(q(l))=\|l\|$.
\end{cory}

\proof
Using 1.2.25, 1.2.27 and 1.2.10, we obtain $\|l\|=
\|q(l)\|=d(q(l),O)
\leq\rho(q(l))
\leq\|l\|$, the last inequality being due to the fact that $l$ is
$\|l\|$-lipschitzian.
\cqfd

\begin{remk}
In 1.2.28, we can replace ``linear'' by ``affine''.
The reader can go back to 1.2.18 (one can refer to 1.5.12).
\end{remk}


%% file: jet_I2.tex
\section{Tangentiability}
In this new context, the notion of tangentiability plays the part of the one of differentiability. However, this new notion keeps lots of the properties 
of the classical differentiability. 
This section is a good opportunity to give some examples and 
counter-examples which will be indispensable to understand and visualize these new notions.
Still here, ``jet'' will mean ``linked metric jet''.

\begin{defi}
Let $f:M\lra M'$ be a map between metric spaces and $a\in M$.
We say that $f$ is \textup{tangentiable} at $a$ (in short $Tang_a$)
if there exists a map $g:M\lra M'$ which is $LL_a$ such that
$g\tang_a f$.

When $f$ is $Tang_a$, we set
\textup{T}$f_a=\{g:M\lra M'\,|\, g\tang_a f\,;\ g \ LL_a \}$; 
\textup{T}$f_a$ is a jet
$\, (M,a)\lra(M',f(a))$, said tangent to $f$ at $a$, and
that we can call the 
\textup{tangential} of $f$ at  $a$ (not to be mixed up with 
\textup{t}$f_a$, defined in a special context (see before 1.5.1)).
\end{defi}

\begin{prop}
We have the implications:
$f\ LL_a\ \Lra\ f\ Tang_a$ and
$f\ Tang_a\ \Lra f\ LSL_a\ \Lra$ continuous at $a$.
The inverse implications are not true (see 1.3.9).
\end{prop}

\proof
The first implication is obvious; the second one comes from
1.1.8 and 1.1.11, and the last one has been seen in 1.1.8.
\cqfd

\begin{remk}
$f$ is $LL_a$ iff $f$ is $Tang_a$ with $f\in \textup{T}f_a$ 
(then $\textup{T}f_a=q(f)$); 
in particular, $\textup{T}(Id_M)_a=q(Id_M)=Id_{(M,a)}$
(see the definition of the canonical surjection
$q:\L\textup{L}((M,a),(M',a'))\lra\J\textup{et}((M,a),(M',a'))$ just after  1.2.1);
moreover, for every jet $\varphi:(M,a)\lra(M',a')$, we have $\varphi=
\textup{T}g_a$ for every $g\in \varphi$.
\end{remk}

\begin{prop}
Let $E,E'$ be n.v.s., $U$ an open subset of $E$, $a\in U$ and 
$f:U\lra E'$ a map. We have the implication
(the inverse being not true (see 1.3.9):
$f\, $ differentiable at $\,a\ \Lra\ f\, $ tangentiable at $\,a$. 
\end{prop}

\proof
A continuous affine map being lipschitzian.
\cqfd

\begin{remk}
Actually, $f$ is differentiable at $a$ iff $f$ is $Tang_a$ where 
its tangential
$\textup{T}f_a$ at $a$ possesses an affine map.
For such a differential map, it is the unique affine map $Af_a$
defined by $Af_a(x)=f(a)+\textup{d}f_a(x-a)$(the translate at $a$ of 
its differential $\textup{d}f_a$); thus $\textup{T}f_a=q(Af_a)$.
In particular, if $f$ is a continuous
affine map, then $f$ is $Tang_a$ at every $a$
with $f\in \textup{T}f_a$ for all $a$.
\end{remk}

\begin{prop}
Let $M,M',M''$ be metric spaces, $f:M\lra M'$ and $g:M'\lra M''$ two maps, and $a\in M$, $a'=f(a)$.
If $f$ is $Tang_a$ and $g$ is $Tang_{a'}$, then $g.f$ is $Tang_a$
and we have $\textup{T}(g.f)_a=\textup{T}g_{a'}.\textup{T}f_a$.
\end{prop}

\proof
Comes from 1.1.14,  1.1.15, and from the functoriality of $q$.
\cqfd

\begin{prop}
Let $M,M_0,M_1$ be metric spaces, $f_0:M\lra M_0$ and 
$f_1:M\lra M_1$ be two maps, and $a\in M$.
If $f_0$ and $f_1$ are tangentiable at $a$, then
$(f_0,f_1):M\lra M_0\times M_1$ is tangentiable at $a$
and we have $\textup{T}(f_0,f_1)_a=
(\textup{T}f_{0a},\textup{T}f_{1a})$.
\end{prop}

\proof
Comes from 1.1.17, 1.1.18 and 1.2.2.
\cqfd

\begin{prop}
Let $M_0,M_1,M'_0,M'_1$ be four metric spaces,\break
$f_0:M_0\lra M'_0$ and 
$f_1:M_1\lra M'_1$ two maps, and $a_0\in M_0$, $a_1\in M_1$.
If $f_0$ is $Tang_{a_0}$ and $f_1$ $Tang_{a_1}$, then
$f_0\times f_1:M_0\times M_1\lra M'_0\times M'_1$ is 
$Tang_{(a_0,a_1)}$
and we have $\textup{T}(f_0\times f_1)_{(a_0,a_1)}=
\textup{T}f_{0a_0}\times 
\textup{T}f_{1a_1}$.
\end{prop}

\proof
Since $f_0\times f_1=(f_0.p_0,f_1.p_1)$ where the canonical projections $p_i$ are lipschitzian. Beyond 1.3.7, we use the fact that, for all 
$(a_0,a_1)\in M_0\times M_1$, we have 
T$p_{i(a_0,a_1)}=q(p_i)=\pi_i$, the canonical projections in $\J$et
(see 1.3.3 and 1.2.2).
\cqfd

\vspace{3mm}

\begin{examsc}
{}
\end{examsc}
All the maps considered below are functions $\R\lra\R$, where $\R$ is equipped with its usual structure of normed vector space
(these different functions give counter-exemples to the inverses of the implications given in 1.3.2 and 1.3.4).

1) Consider $f_0(x)=x^{1/3}$; this function is continuous but not 
$LSL_0$ (see 1.1.6 and 1.1.10).

2) Consider $f_1(x)=x\sin\frac{1}{x}$ if $x\not=0$ and $f_1(0)=0$; this function is obviously $LSL_0$, however not $Tang_0$: indeed, if $f_1$ was $Tang_0$, there would exist a fonction $g:\R\lra\R$ and a neighborhood $V$ of 0 such that $g\tang_0f_1$ and $g|_V$ is 
$k$-lipschitzian for a $k>0$. Let us consider the two sequences of reals defined by $x_n=1/2n\pi$ and $y_n=1/(4n+1)\frac{\pi}{2}$; since
$g\tang_0 f_1$ and $\lim_n x_n=\lim_n y_n=0$, we have
$\lim_n|\frac{f_1(x_n)-g(x_n)}{x_n}|=0=
\lim_n|\frac{f_1(y_n)-g(y_n)}{y_n}|$, so  that 
$\lim_n2\pi n g(x_n)=0$ and $\lim_n(4n+1)\frac{\pi}{2} g(y_n)=1$.
Now, since $g|_V$ is $k$-lipschitzian, we have
$|\frac{g(x_n)-g(y_n)}{x_n-y_n}|\leq k$ for $n$ big enough,
which is equivalent to $|\frac{4n+1}{n}(2n\pi g(x_n))
-4((4n+1)\frac{\pi}{2}g(y_n))|\leq\frac{k}{n}$.
It remains to do $n\rightarrow +\infty$ which leads to a contradiction.

3) Consider $f_2(x)=x^2\sin\frac{1}{x^2}$ if $x\not=0$ and $f_2(0)=0$; this function is $Tang_0$ (since it is differentiable at 0); however
not $LL_0$ (because $lim_{k\rightarrow +\infty}f'_2(\frac{1}
{\sqrt{2k\pi}})=-\infty$).

4) Consider $\vartheta(x)=|x|$; this function is $Tang_0$ (it is 1-lipschitzian!),
however not differentiable at 0. 

\begin{remk}
We will see in the second chapter
that this function $\vartheta(x)=|x|$ can be viewed as a model of 
G-differentiable maps;
in this second chapter, we will also see lots of other 
tangentiable maps which are not differentiable, some of them being not even G-differentiable.
\end{remk}

\vspace{3mm}

Consider now two metric spaces $M,M'$, and two maps\break
$f,g:M\lra M'$ which are $Tang_a$ where $a\in M$.
Like for general jets (see 1.2.7), we can speak of the distance 
between the two jets T$f_a$, T$g_a$.
The question is: do we still have $d$(T$f_a$,T$g_a)=d(f,g)$?\break
That is a natural question, except that, until now, this $d(f,g)$ does not mean anything here, since
it has been solely defined for maps which are $LL_a$! (see 
just before 1.2.6). Here, our maps $f$ and $g$ are
only supposed to be $Tang_a$. Therefore,
it remains to extend this definition $d(f,g)$ to such maps.

Let $M,M'$  be metric spaces, $a\in M$, $a'\in M'$ and 
$f,g:M\lra M'$\break
two maps which are $Tang_a$ and which verify $f(a)=g(a)=a'$.
Let us consider $f_1\in\, $T$f_a$ and $g_1\in\, $T$g_a$; then $f_1\tang_a f$ and 
$g_1\tang_a g$; so that, for every $\varepsilon>0$, there exists a neighborhood $V$ of $a$ on which
we have $d(f_1(x),f(x))\leq \varepsilon d(x,a)$ and 
$d(g_1(x),g(x))\leq \varepsilon d(x,a)$.
Furthermore, since $f_1$ and $g_1$ are $LL_a$, we know that there exists also a neighborhood $W$ of $a$ on which the map 
$x\mapsto \frac{d(f_1(x),g_1(x))}{d(x,a)}$, if $x\not= a$, is bounded
(let us say by $R$).
Now, if we take $x\in V_1=V\cap W$, $x\not=a$, we obtain:\par
\noindent$\frac{d(f(x),g(x))}{d(x,a)}\leq\frac{d(f(x),f_1(x))}{d(x,a)}+
\frac{d(f_1(x),g_1(x))}{d(x,a)}+\frac{d(g_1(x),g(x))}{d(x,a)}
\leq 2\varepsilon+R$.
So, the map $C(x)=\frac{d(f(x),g(x))}{d(x,a)}$, if $x\not=a$ and $C(a)=0$,
is still bounded on $V_1$.
Thus, for $r>0$, we can again set $d^r(f,g)=\sup\{C(x)|x\in V_1
\cap B'(a,r)\}$
and finally again $d(f,g)=\lim_{r\rightarrow 0}d^r(f,g)=\inf_{r>0}d^r(f,g)$.

\begin{prop}
If $f,g:M\lra M'$ are $Tang_a$ where $a\in M$, we have
$d(\textup{T}f_a,\textup{T}g_a)=d(f,g)$, this $d$ being defined just above.
\end{prop}

\proof
T$f_a=$T$f_{1a}$ and T$g_a=$T$g_{1a}$ where $f_1\in$T$f_a$ and
$g_1\in$T$g_a$;\break
using the fact that $f_1$ and $g_1$ are $LL_a$, 
we obtain $d($T$f_{1a},$T$g_{1a})=d(f_1,g_1)$.
The result comes from the fact that, since  
$f\tang_a f_1$ and $g\tang_a g_1$,
for every $\varepsilon>0$ and $r>0$, we have,
according to 1.3.10,\break 
$d^r(f,g)\leq d^r(f_1,g_1)+2\varepsilon$ and also
$d^r(f_1,g_1)=d^r(f,g)+2\varepsilon$, which implies
$d(f,g)=d(f_1,g_1)$.
\cqfd

\begin{cory}
Let $f:M\lra M'$ be a map, $a\in M$ and $b=f(a)$.\break
We assume that $f$ is $Tang_a$ and that $d(\textup{T}f_a,O_{ab})<k$.
Then, $f$ is\break $k$-$LSL_a$.
\end{cory}

\proof
Since $O_{ab}=q(\widehat b)=$T$\widehat{b}_a$ 
(where $\widehat b$ is the constant map on $b$), we have, using  1.3.11, $d(f,\widehat b)=d($T$f_a,$T$\widehat b_a)=
d($T$f_a,O_{ab})<k$.
So that,
referring to before 1.3.11,  there exists $r>0$ such that  
$d^r(f,\widehat b)<k$.\break
This means that, when $r$ is small enough, 
we have $d(f(x),f(a))=d(f(x),b)\leq kd(x,a)$  for every $x\in B'(a,r)$. We have thus proved that
$f$ is $k$-$SL_a$ on $B'(a,r)$. 
\cqfd

\section{Transmetric spaces}

In order to define a tangential map t$f:M\lra``\J$et$(M,M')"$ for a tangentiable map 
$f:M\lra M'$ (see the section 1.5 below), we need, of course, to 
define first such a set $``\J$et$(M,M')"$ equipped with an adequate distance. We have succeeded in it,
assuming that $M$ and $M'$ are transmetric spaces
(see 1.4.1 below). Then, such a metric space $``\J$et$(M,M')"$ can really be defined, its elements being new metric jets,  
called free metric jets (in opposition to the linked metric jets
we have used up to now and that we continue to call briefly jets, still here);
it's why the set $``\J$et$(M,M')"$ will be in fact denoted 
$\J$et$\jf(M,M')$. 

\begin{defi}
A \textup{transmetric space} is a metric space $M$, supposed to be non empty, equipped with a functor
$\gamma:Gr(M)\lra\J$et (where the category $\J$\textup{et} has been defined at the beginning of section 1.2, and $Gr(M)$ refers to the groupoid of the pairs, defined by $|Gr(M)|=M$ and $Hom(a,b)=\{(a,b)\}$ for 
all $a,b\in M$) verifying:

- for every $a\in M$, $\,\gamma(a)=(M,a)$,

- for every morphism $(a,b):a\lra b$ in $Gr(M)$, the jet\break
$\gamma(a,b):(M,a)\lra(M,b)$ is 1-bounded (i.e. verifies
$\rho(\gamma(a,b))\leq 1$: see 1.2.8); every $\gamma(a,b)$ is invertible
in $\J$\textup{et}. 

An attentive reader can be ``puzzled'' by the restrictive hypothesis
``non empty'' for transmetric spaces; this reader will find some justifications in 1.4.14 and 1.4.22.
\end{defi}

Before giving some examples, let us give the following special case:

\begin{defi}
A \textup{left isometric group} is a metric space $G$ equipped with a
group structure verifying the following condition:\par
\noindent$\forall g,g_0,g_1\in G\quad(d(g.g_0,g.g_1)=d(g_0,g_1))$.
\end{defi}

\begin{remks}
\par\hfill

1) In an equivalent manner, we could have only assume
in 1.4.2  that:
$\ \forall g,g_0,g_1\in G\quad(d(g.g_0,g.g_1)\leq d(g_0,g_1))$.

2) Let $G$ be a left isometric group; then,  for all $g\in G$, the map $G\lra G:g'\mapsto g.g'$  is isometric.
\end{remks}

\begin{prop}
Every left isometric group $G$ can be equipped with a canonical structure of transmetric space.
\end{prop}

\proof
Here, $\gamma$ is the composite
$Gr(G)\!\buildrel{\theta}\over\lra\!\L$L$
\ \buildrel{q}\over\lra\! \J$et,
where $\theta$ is the functor defined by $\theta(g)=(G,g)$ and
$\theta(g_0,g_1)$ is the morphism $(G,g_0)\lra(G,g_1)$ 
in $\L$L which assignes $g_1.g_0^{-1}.g$
to $g$. The $\theta(g_0,g_1)$ being isometries, the 
$\gamma(g_0,g_1)$ are 1-bounded (see 1.2.11).
\cqfd

\vfill\eject

\begin{exams}
{}
\end{exams}
1) Here are examples of transmetric spaces which are even
left isometric groups:

\quad$\ $ a) Every n.v.s. is a left isometric (additive) group: here, $\gamma(a,b)=q(\theta(a,b))$, where
$\theta(a,b)$ is the translation
$x\mapsto b-a+x$.

\quad$\ $ b) The multiplicative group $S^1=\{z\in\C|\,|z|=1\}$ is also
a left isometric group.

\quad$\ $ c) $E$ being an euclidian space, the orthogonal group $O(E)$, 
equipped with its operator norm, is a left isometric group.

\quad$\ $ d) The additive subgroups of $\R$ (as, for example $\Q$)
are left isometric groups whuch are not e.v.n..

2) Every non empty discreet space has a unique structure of transmetric space; conversely, if a transmetric space possesses an isolated point, it is a discreet space. Indeed, since $M$ discreet
means that all its points are isolated, we know, by 1.2.12, that,
for all $a\in M$, $(M,a)$ is a final object in $\J$et, so that 
there exists a unique functor $\gamma:Gr(M)\lra\J$et\break
verifying 
$\gamma(a)=(M,a)$;
this $\gamma$ defines a structure of transmetric space on $M$, since, referring to 1.2.11, we know that, for all $a,b\in M$
$\rho(\gamma(a,b))=0\leq 1$, since $a$ and $b$ are isolated.
Conversely, if $a$ is an isolated point in a transmetric space $M$,
we use the fact that, by definition,  for every $b\in M$, the jet
$\gamma(a,b):(M,a)\lra(M,b)$ is an isomorphism, to deduce (still by  1.2.12) that $b$ is also isolated.

3) The below 1.4.7 will provide a lot of transmetric spaces which are not 
left isometric groups (an open subset of a left isometric group
being not necessarily a subgroup!).

\vspace{3mm}
\begin{prop}
Let $M_0$ and $M_1$ be transmetric spaces; then\break 
$M_0\times M_1$ has a canonical stucture of transmetric space.
\end{prop}

\proof
Let $\gamma_0:Gr(M_0)\lra\J$et and $\gamma_1:Gr(M_1)\lra\J$et be
the functors defining the transmetric structrures on 
$M_0$ and $M_1$, and consider the functor 
$\gamma:Gr(M_0\times M_1)\lra\J$et
defined by $\gamma((a_0,a_1),(b_0,b_1))=
\gamma_0(a_0,b_0)\times\gamma_1(a_1,b_1)$
(which makes sense, since $\J$et is a cartesian category
(see 1.2.2).
The fact that 
$\rho(\gamma((a_0,a_1),(b_0,b_1)))\leq 1$ arises
immediately from 1.2.14, since 
the $\gamma_i(a_i,b_i)$ are 1-bounded.
\cqfd

\begin{prop}
Let $M$ be a transmetric space and $U\not=\emptyset$, an open subset of $M$;
then, $U$ has a canonical structure of transmetric space.
\end{prop}

\proof
Let $\gamma:Gr(M)\lra\J$et be the functor defining the transmetric structure on $M$.
Now, for each $a,b\in U$, we set $\breve{\!\gamma}(a)=(U,a)$ and 
$\breve{\!\gamma}(a,b)=j_b^{-1}.\gamma(a,b).j_a$
(refer to 1.2.5 for the definition of the jets $j_a$ and $j_b$; and 
peep at the diagram just below); this $\breve{\!\gamma}$
is obviously a functor $\ Gr(U)\lra\J$et.
Besides, using 1.2.9, we have
$\rho(\breve{\!\gamma}(a,b))\leq
\rho(j_b^{-1})\rho(\gamma(a,b))\rho(j_a)\leq\rho(\gamma(a,b))\leq 1$,
since the jets $j_a$ and $j_b^{-1}$ are 1-bounded (see 1.2.11).
$$\xymatrix{
(U,a)\ar[d]_{j_a}\ar[r]^{\breve{\!\gamma}(a,b)}&(U,b)\ar[d]^{j_b}\\
(M,a)\ar[r]_{\gamma(a,b)}&(M,b)
}$$
\cqfd

\begin{defi}
Let $M,M'$ be transmetric spaces; a map\break
$f:M\lra M'$ is called a \textup{morphism
of transmetric spaces} if it is a tangentiable map (at every point of $M$)
such that, for every $a,b\in M$, the following diagram commutes in the category
$\J$\textup{et}:
$$\xymatrix{
(M,a)\ar[d]_{\gamma(a,b)}\ar[r]^<<<<<{\textup{T}f_a}&(M',f(a))
\ar[d]^{\gamma(f(a),f(b))}\\
(M,b)\ar[r]_<<<<<{\textup{T}f_b}&(M',f(b))
}$$
We will denote $\T$\textup{rans} the category whose objects are the transmetric spaces
and whose morphisms are the morphisms of transmetric spaces.
\end {defi}

\vfill\eject

\begin{exams}
{}
\end{exams}
1) Every continuous affine map $f:E\lra E'$ (between n.v.s.) is a morphism of transmetric spaces: it is tangentiable since it is lipschitzian;
and the commutativity of the diagram of 1.4.8 results from the fact that, for every 
$a,b,x\in E$, we have the
following commutative diagram in the category $\L$L
(where the $\theta(a,b)$ are the translations of a n.v.s. defined in  1.4.5):
$$\xymatrix{
(E,a)\ar[d]_{\theta(a,b)}\ar[r]^<<<<<{f}&(E',f(a))\ar[d]^{\theta(f(a),f(b))}\\
(E,b)\ar[r]_<<<<<{f}&(E',f(b))
}$$
We have just to transport this diagram in the category $\J$et with the help of  the functor
$q:\L$L$\lra\J$et (see just before 1.2.2, and 1.4.5).

For instance, the translations $\theta(a,b)$ of a n.v.s. are morphisms of transmetric spaces.

2) If $M$ is a transmetric space, $Id_M:M\lra M$ is a morphism of transmetric spaces since, according to 1.3.3,
T${(Id_M)}_a=Id_{(M,a)}$ is an identity for each $a\in M$.

3) Every constant map (on $c\in M'\ $)
$\widehat c :M\lra M'$ between transmetric spaces is a morphism of transmetric spaces, since, it is lipschitzian (thus tangentiable, with
T$\widehat c_a=O_{ac}$ for all $a\in M$).
For instance, for every transmetric space $M$, 
the unique map $!_M:M\lra\I$ is a morphism of transmetric spaces
(where $\I=\{0\}$).

4) Let $M$ be a transmetric space and $U\not=\emptyset$, an open subset of $M$; the diagram included into the proof of 1.4.7 shows that the canonical injection 
$j:U\hookrightarrow M$ is a morphism of transmetric spaces. 

5) $M_0$ and $M_1$ being transmetric spaces, the canonical projections
$p_i:M_0\times M_1\lra M_i$ are morphisms of transmetric spaces.
It comes from the fact that, in a cartesian category, the projections
are natural.

\begin{prop}
The category $\T$\textup{rans} is a cartesian category.
\end{prop}

\proof
$\I=\{0\}$ is a final object in $\T$rans (see the example 3) of 1.4.9).
Now, referring to the above example 5) of 1.4.9,
let $M,M_0,M_1$ be three transmetric spaces and, for each
$i\in\{0,1\}$,
$f_i:M\lra M_i$ a morphism of transmetric spaces.
By 1.3.7, we know that $(f_0,f_1):\break
M\lra M_0\times M_1$ is tangentiable and that
for every $a\in M$, we have T$(f_0,f_1)_a =($T$f_{0a},$T$f_{1a})$;  
using then the definition of $\gamma:\break
Gr(M_0\times M_1)\lra\J$et
given in 1.4.6, 
we deduce the commutativity of the following diagram in the cartesian category $\J$et
(where $a_i=f_i(a)$ and $b_i=f_i(b)$):
$$\xymatrix{
(M,a)\ar[d]_{\gamma(a,b)}
\ \ar[rrrr]^<<<<<<<<<<<<<<<<<<<<{(\textup{T}f_{0a},
\textup{T}f_{1a})}&&&&\ 
(M_0,a_0)\times (M_1,a_1)
\ar[d]^{\gamma(a_0,b_0)\times \gamma(a_1,b_1)}\\
(M,b)\ \ar[rrrr]_>>>>>>>>>>>>>>>>>>>>>{(\textup{T}f_{0b},
\textup{T}f_{1b})}&&&&\ 
(M_0,b_0)\times (M_1,b_1)
}$$
Therefore, the map $(f_0,f_1):M\lra M_0\times M_1$ is a morphism
of transmetric spaces.
\cqfd

\vspace{3mm}
We are now ready to 
construct the foretold set $\J$et$\jf(M,M')$ of the free metric jets,
when $M$ and $M'$ are transmetric spaces.
First, we consider the set
$J(M,M')=\coprod_{(a,a')\in M\times M'}\J$et$((M,a),(M',a'))$
on which we define the following equivalence relation:
$(\varphi,a,a')\sim(\psi,b,b')$ if the following diagram commutes 
in the category $\J$et:
$$\xymatrix{
(M,a)\ar[d]_{\gamma(a,b)}\ar[r]^<<<<<{\varphi}&(M',a')\ar[d]^{\gamma(a',b')}\\
(M,b)\ar[r]_<<<<<{\psi}&(M',b')
}$$
The transmetric structures of $M$ and $M'$ are essential to this
definition.
Then, we set $\J$et$\jf(M,M')=J(M,M')/\sim\,$.

\begin{defi}
The elements of $\J$\textup{et}$\jf(M,M')$ are called 
\textup{free metric jets} (here from $M$ to $M'$),
in opposition to the elements of\break
$\J$\textup{et}$((M,a),(M',a'))$
called linked metric jets (here from $(M,a)$ to\break
$(M',a')$), or still shortly jets.
\end{defi}

If $q:J(M,M')\lra\J$et$\jf(M,M')$ is the canonical surjection, we set
$[\varphi,a,a']=q(\varphi,a,a')$ when $(\varphi,a,a')\in J(M,M')$.

We are now going to canonically equip $\J$et$\jf(M,M')$ with a structure of metric space.
First, for $(\varphi,a,a'),(\psi,b,b')\in J(M,M')$, we denote\break
$d((\varphi,a,a'),(\psi,b,b'))$ the
distance between the linked metric jets\break 
$\gamma(a',b').\varphi$ and $\psi.\gamma(a,b)$ (given in 1.2.7).
However, the first member of this equality is not a distance on 
$J(M,M')$ (see 1.4.12 below); but it will become a ``true'' distance for the quotient 
$\J$et$\jf(M,M')$.
 
\begin{prop}\par\hfill
For each $(\varphi,a,a'),(\psi,b,b'),(\xi,c,c')\in J(M,M')$, we have the following properties:

1) $d((\varphi,a,a'),(\psi,b,b'))=d((\psi,b,b'),(\varphi,a,a'))$,

2) $d((\varphi,a,a'),(\xi,c,c'))\leq d((\varphi,a,a'),(\psi,b,b'))+d((\psi,b,b'),(\xi,c,c'))$,

3) $d((\varphi,a,a'),(\psi,b,b'))=0\ \Longleftrightarrow\ (\varphi,a,a')\sim(\psi,b,b')$. 
\end{prop}

\proof
1.2.10, 1.2.15 and 1.2.16 are very usefull here (reminding that the jets $\gamma(a,b)$ are 
1-bounded and invertible: with $\gamma(a,b)^{-1}=\gamma(b,a)$).

1) $d(\psi,b,b'),(\varphi,a,a'))=d(\gamma(b',a').\psi,\varphi.\gamma(b,a))$\par
\qquad\qquad\qquad\qquad\qquad = $d(\gamma(b',a').\psi,\gamma(b',a').\gamma(a',b').\varphi.\gamma(b,a))$\par
\qquad\qquad\qquad\qquad\qquad $\leq d(\psi,\gamma(a',b').\varphi.\gamma(b,a))$\par
\qquad\qquad\qquad\qquad\qquad $=d(\psi.\gamma(a,b).\gamma(b,a),\gamma(a',b').\varphi.\gamma(b,a))$\par
\qquad\qquad\qquad\qquad\qquad $\leq d(\psi.\gamma(a,b),\gamma(a',b').\varphi)$\par
\qquad\qquad\qquad\qquad\qquad $=d(\gamma(a',b').\varphi,\psi.\gamma(a,b))$\par
\qquad\qquad\qquad\qquad\qquad $=d((\varphi,a,a'),(\psi,b,b'))$.\par
\noindent Same for $d((\varphi,a,a'),(\psi,b,b'))\leq d((\psi,b,b'),(\varphi,a,a'))$.

2) $d((\varphi,a,a'),(\xi,c,c'))=d(\gamma(a',c').\varphi,\xi.\gamma(a,c))\leq$\par
\noindent $d(\gamma(a',c').\varphi,\gamma(b',c').\psi.\gamma(a,b))+d(\gamma(b',c').\psi.\gamma(a,b),\xi.\gamma(a,c))=$\par
\noindent $d(\gamma(b',c').\gamma(a',b').\varphi,\gamma(b',c').\psi.\gamma(a,b))+$\par
\qquad\qquad\qquad\qquad\quad\ \, $d(\gamma(b',c').\psi.\gamma(a,b),\xi.\gamma(b,c).\gamma(a,b))\leq$\par
\noindent $d(\gamma(a',b').\varphi,\psi.\gamma(a,b))+d(\gamma(b',c').\psi,\xi.\gamma(b,c))=$\par
\noindent $d((\varphi,a,a'),(\psi,b,b'))+d((\psi,b,b'),(\xi,c,c'))$.

3) Obvious.
\cqfd

\begin{prop}
The map $d:(J(M,M'))^2\lra\R_+$, studied in 1.4.12,
factors through the quotient, 
giving a ``true'' distance on\break
$\J$\textup{et}$\jf(M,M')$, defined by 
$d([\varphi,a,a'],[\psi,b,b'])=d((\varphi,a,a'),(\psi,b,b'))$ for all
$(\varphi,a,a'),(\psi,b,b')\in J(M,M')$. 
\end{prop}

\begin{remk}
\par\hfill

Let $M,M'$ be transmetric spaces;
then $\J$\textup{et}$\jf(M,M')$
possesses a particular element denoted $O$: indeed, we notice that, for 
every $a,b\in M$, $a',b'\in M'$, we have $(O_{aa'},a,a')\sim(O_{bb'},b,b')$, so that
$[O_{aa'},a,a']$ does not depend on the choice of $(a,a')\in M\times M'$. It's this free metric jet that we denote $O$ and that we call the \textup{free zero} of $\J$\textup{et}$\jf(M,M')$.
In the same way,  $\J$\textup{et}$\jf(M,M)$ possesses also a \textup{free identity}, denoted $I_M$, which is equal to
$[Id_{(M,x)},x,x]$: this free metric jet not depending on $x\in M$,  since $Id_M$ is a morphism of transmetric spaces (see 1.4.9).

We notice that the existence of such $O$ and $I_M$ arises from the fact that transmetric spaces are supposed to be non empty.
\end{remk}

\begin{prop}
For each $a\in M$, $a'\in M'$, the following composite (denoted $can$):
$\J$\textup{et}$((M,a),(M',a'))\!\buildrel c\over\lra\! J(M,M')\!\buildrel q\over\lra\!
\J$\textup{et}$\jf(M,M')$
(where $c(\varphi)=(\varphi,a,a')$) is an isometry.
\end{prop}

\proof
$can$ is bijective: if $[\psi,b,b']\in\J$et$\jf(M,M')$, then 
$[\psi,b,b']=can(\varphi)$, where $\varphi=\gamma(b',a').\psi.\gamma(a,b)$ is the
unique element of\break
$\J$et$((M,a),(M',a'))$ to verify $(\varphi,a,a')\sim(\psi,b,b')$.

$can$ is isometric since, if $\varphi,\varphi'\in \J$et$((M,a),(M',a'))$,
we have\break
$d([\varphi,a,a'],[\varphi',a,a'])=d((\varphi,a,a'),(\varphi',a,a'))=d(\varphi,\varphi')$, 
since\break
$\gamma(a,a)=Id_{(M,a)}$ and $\gamma(a',a')=Id_{(M',a')}$.
\cqfd

\begin{prop}
Let $M,M'$ be transmetric spaces. For each\break
$a,b\in M$,
$a',b'\in M'$, we set $\Omega(\varphi)=\gamma(a',b').\varphi.\gamma(b,a)$; it defines a map
$\Omega:\J$\textup{et}$((M,a),(M',a'))\lra
\J$\textup{et}$((M,b),(M',b'))$ which is an isometry, and the following diagram commutes
(where the maps $can$ are defined in 1.4.15):
$$\xymatrix{
\J\textup{et}((M,a),(M',a'))\ar[rr]^{\Omega}\ar[dr]_{can}&&
\J\textup{et}((M,b),(M',b'))\ar[dl]^{can}\\
{}&\J\textup{et}\jf(M,M')&{}
}$$
\end{prop}

\proof
The fact that $(\varphi,a,a')\sim(\Omega(\varphi),b,b')$
provides the commutativity of the diagram.
$\Omega$ is an isometry just as the maps $can$.
\cqfd

\vspace{3mm}
Let now $M_0,M_1,M_2$ be three transmetric spaces, 
and also:\par
\noindent $\varphi_0:(M_0,a_0)\lra(M_1,a_1)$ and 
$\varphi_1:(M_1,b_1)\lra(M_2,b_2)$ be two jets; 
then, we consider the natural composition (as in the diagram below)\par
\noindent $(\varphi_1,b_1,b_2).(\varphi_0,a_0,a_1)=
(\varphi_1.\gamma(a_1,b_1).\varphi_0,a_0,b_2)$:
$$\xymatrix{
(M_0,a_0)\ar[r]^{\varphi_0}& (M_1,a_1)
\ar[d]^{\gamma(a_1,b_1)}&{}\\
{}&(M_1,b_1)
\ar[r]^{\varphi_1}&(M_2,b_2)
}$$
\noindent This composition defines a map:\par
\qquad\qquad$comp:J(M_0,M_1)\times J(M_1,M_2)
\lra\ J(M_0,M_2)$

\begin{prop}
\noindent This map $comp$ factors through the quotient:
$$\xymatrix{
J(M_0,M_1)\times J(M_1\times M_2)\ar[d]_{q\times q}
\ \ar[rr]^<<<<<<<<<<<<{comp}&&\ J(M_0,M_2)
\ar[d]^{q}\\
\J\textup{et}\jf(M_0,M_1)\times\J\textup{et}\jf(M_1,M_2)\ \ar[rr]_>>>>>>>>>>
{comp}&&\ \J\textup{et}\jf(M_0,M_2)
}$$
\noindent (and we will write $[\varphi_1,b_1,b_2].
[\varphi_0,a_0,a_1]=
[\varphi_1.\gamma(a_1,b_1).\varphi_0,a_0,b_2]$, which gives simply
$[\varphi_1,b_1,b_2].[\varphi_0,a_0,a_1]
\buildrel\star\over=[\varphi_1.\varphi_0,a_0,b_2]$ when $a_1=b_1$).
\end{prop}

\proof
This factorisation is merely due to the fact that, if\break
$(\varphi'_0,a'_0,a'_1)\sim(\varphi_0,a_0,a_1)$ and
$(\varphi'_1,b'_1,b'_2)\sim(\varphi_1,b_1,b_2)$, then we have
$\varphi'_1.\gamma(a'_1,b'_1).\varphi'_0\sim \varphi_1.\gamma(a_1,b_1).\varphi_0$,
since the following diagrams commute:
$$\xymatrix{
(M_0,a_0)\ar[d]_{\gamma(a_0,a'_0)}\ar[r]^<<<<<{\varphi_0}&(M_1,a_1)
\ar[d]_{\gamma(a_1,a'_1)}\ar[rr]^{\gamma(a_1,b_1)}&&(M_1,b_1)\ar[d]^{\gamma(b_1,b'_1)}
\ar[r]^{\varphi_1}&(M_2,b_2)\ar[d]^{\gamma(b_2,b'_2)}\\
(M_0,a'_0)\ar[r]_<<<<<{\varphi'_0}&(M_1,a'_1)\ar[rr]_{\gamma(a'_1,b'_1)}&&(M_1,b'_1)
\ar[r]_{\varphi'_1}&(M_2,b'_2)
}$$
\cqfd

\begin{prop}
The composition defined just above:\par
\noindent$comp:\J\textup{et}\jf(M_0,M_1)\times\J\textup{et}\jf(M_1,M_2)
\lra\J\textup{et}\jf(M_0,M_2)$ is $LSL$.
\end{prop}

\proof
We just have to use 1.2.19 and to notice that the following diagram commutes (where $a_0\in M_0$, $a_1\in M_1$, $a_2\in M_2$; and $can$ has been defined in 1.4.15;
to let the diagram enter the page we have denoted 
$\J\textup{et}((M_i,a_i),(M_j,a_j))$ by $\J$et$_{ij}$):
$$\xymatrix{
\J\textup{et}_{01}\times \J\textup{et}_{12}
\ar[d]_{can\times can}
\ \ar[rr]^<<<<<<<<<<<<<<<<<<<<{comp}&&\ \J\textup{et}_{02}
\ar[d]^{can}\\
\J\textup{et}\jf(M_0,M_1)\times\J\textup{et}\jf(M_1,M_2)\ \ar[rr]_>>>>>>>>>>
{comp}&&\ \J\textup{et}\jf(M_0,M_2)
}$$
\cqfd

So, we are now in a position to construct a new category, enriched 
in $\M$et, denoted $\J$et$\jf$,
called the category of free metric jets,
whose:

- objects are the transmetric spaces $M$, 
 
- ``Hom'' are the metric spaces $\J$et$\jf(M,M')$,
 
- identity $\I\lra\J$et$\jf(M,M)$, is the map giving 
the free identity\break
$I_M=[Id_{(M,a)},a,a]$ defined in 1.4.14,

- composition $\J\textup{et}\jf(M,M')\times\J\textup{et}\jf(M',M'')\lra
 \J\textup{et}\jf(M,M'')$ is the previous $comp$.

We notice that the hypothesis for transmetric spaces to be non empty is essential here (see 1.4.14 and 1.4.22).

\begin{defi}
We will denote $\J\textup{et}'$ the full subcategory of
the category $\J\textup{et}$
whose objects are the pointed transmetric spaces. Just as $\J\textup{et}$
(see 1.2.20),
it is a category enriched in $\M\textup{et}$. 
\end{defi}

\begin{prop}
Actually, we have a forgetfull enriched functor\break
$U:\J\textup{et}'\lra\J\textup{et}$.
\end{prop}

\proof
On the object, $U$ only keeps the metric structure; 
$U$ is an identity on each Hom (which provides the enriched 
property).
\cqfd

\begin{prop}
We have a functor $can':\J\textup{et}'\lra\J\textup{et}\jf$, defined by $can'(M,a)=M$ and, for $\varphi:
(M,a)\lra(M',a')$, $can'(\varphi)=can(\varphi)=[\varphi,a,a']$
(the $can$ of 1.4.15); 
then, this functor is enriched in $\M\textup{et}$.
\end{prop}

\proof
The fact that $can'$ is a functor comes from the composition
in $\J$et$\jf$ given by the equality $\buildrel\star\over =$ of 
1.4.17; and from the fact that $can'(Id_{(M,a)})=[Id_{(M,a)},a,a]=
I_M$. 
This functor $can'$ is enriched in $\M$et.
\cqfd

\vspace{3mm}
We now come back to cartesian considerations.

\begin{remk}
By definition, $f:M\lra M'$ is a morphism of transmetric spaces
iff we have $(\textup{T}f_x,x,f(x))\sim(\textup{T}f_y,y,f(y))$
for all $x,y\in M$;
so that 
the free metric jet $[\textup{T}f_x,x,f(x)]$ is independent on the choice of $x\in M$; we denote $\kappa(f)$ this element of 
$\J$\textup{et}$\jf(M,M')$. In particular, if $c\in M'$, we have seen in 1.4.9
that the constant map $\widehat c$  is a morphism of transmetric spaces; as $\kappa(\widehat c)=O$ (see 1.4.14), it does not depend on $c$.\break
As noticed in 1.4.14, these $\kappa(f)$ would not exist if the transmetric spaces could be empty.
\end{remk}

\begin{prop}
The map $\T\textup{rans}(M,M')\lra\J\textup{et}\jf(M,M'):\break
f\mapsto\kappa(f)$
extends to a functor $\ \kappa:\T\textup{rans}\lra\J\textup{et}\jf$ which is constant on the objects.
\end{prop}

\proof
Comes from 1.3.6. 
\cqfd

\begin{prop}
The functor $\kappa:\T\textup{rans}\lra\J\textup{et}\jf$ 
creates a cartesian structure on the category
$\J\textup{et}\jf$ ($\kappa$ being constant on the objects, it means that $\J\textup{et}\jf$ is cartesian and 
$\kappa$ a strict morphism of cartesian categories).
\end{prop}

\proof
Clearly, $\kappa(!_M)\ $ is the unique arrow $M\lra \I$ in $\J$et$\jf$, so that
$\kappa(\I)=\I$ is the final object of $\J$et$\jf$.

Now, let $M_0,M_1$ be transmetric spaces. We have seen in 1.4.9 that the canonical projection $p_i:M_0\times M_1\lra M_i$ are morphisms of transmetric spaces verifying T$p_{i(a_0,a_1)}=\pi_i$ 
for all $(a_0,a_1)\in M_0\times M_1$.\break
We are going to show that the 
$\kappa(p_i):M_0\times M_1\lra M_i$ are the canonical projections in 
$\J$et$\jf$. 
Thus, let $M$ be another transmetric space and 
$[\varphi_0,a,a_0]:M\lra M_0$, $[\varphi_1,b,b_1]:M\lra M_1$ be morphisms
in $\J$et$\jf$; we put $\bar a=(a_0,b_1)$ and
$\varphi=
(\varphi_0,\varphi_1.\gamma(a,b))$.
Then, $[\varphi,a,\bar a]:\break
M\lra M_0\times M_1$ is a morphism in $\J$et$\jf$
verifying (we use the equality $\buildrel\star\over =$ of 1.4.17):
$\,\kappa(p_0).[\varphi,a,\bar a]=
[\pi_0,\bar a,a_0].[\varphi,a,\bar a]\buildrel\star\over =[\varphi_0,a,a_0]$, since
$\pi_0$ is a projection in $\J$et (see 1.2.2);  
same for $\kappa(p_1).[\varphi,a,\bar a]=
[\pi_1,\bar a,b_1].[\varphi,a,\bar a]\buildrel\star\over =
[\varphi_1.\gamma(a,b),a,b_1]=[\varphi_1,b,b_1]$.
Now, let $\bar{a'}=(a'_0,a'_1)$ with $(a'_0,a'_1)
\in M_0\times M_1$
and $[\varphi',a',\bar{a'}]:M\lra M_0\times M_1$ 
also
verifying $\kappa(p_0).[\varphi',a',\bar{a'}]=[\varphi_0,a,a_0]$ and       $\kappa(p_1).[\varphi',a',\bar{a'}]=[\varphi_1,b,b_1]$,
i.e. $[\pi_0,\bar{a'},a'_0].[\varphi',a',\bar{a'}]=[\varphi_0,a,a_0]$
and $[\pi_1,\bar{a'},a'_1].[\varphi',a',\bar{a'}]=
[\varphi_1,b,b_1]=[\varphi_1.\gamma(a,b),a,b_1]$;
then we have
$[\pi_0.\varphi',a',a'_0]=[\varphi_0,a,a_0]$, i.e.
$(\pi_0.\varphi',a',a'_0)\sim(\varphi_0,a,a_0)$, i.e.
$\pi_0.\varphi'.\gamma(a,a')
=\gamma(a_0,a'_0).\varphi_0$, and in the same way
$\pi_1.\varphi'.\gamma(a,a')=\break
\gamma(b_1,a'_1).\varphi_1.\gamma(a,b)$.
On the other hand, since the $p_i$ are morphisms of transmetric spaces, we have also $\pi_0.\gamma(\bar a,\bar{a'}).\varphi=
\gamma(a_0,a'_0).\pi_0.\varphi=\break
\gamma(a_0,a'_0).\varphi_0$ 
and
$\pi_1.\gamma(\bar a,\bar{a'}).\varphi=
\gamma(b_1,a'_1).\pi_1.\varphi=
\gamma(b_1,a'_1).\varphi_1.\gamma(a,b)$;
so that
$\varphi'.\gamma(a,a')=\gamma(\bar a,\bar{a'}).\varphi$
which finally means that $(\varphi,a,\bar a)\sim(\varphi',a',\bar{a'})$ i.e.
$[\varphi,a,\bar a]=[\varphi',a',\bar{a'}]$.
\cqfd

\begin{prop}
The functor $can':\J\textup{et}'\lra\J\textup{et}\jf$ defined in 1.4.21 is a cartesian functor.
\end{prop}

\proof
Just as $\J$et, the category $\J\textup{et}'$ is cartesian
(see 1.4.6 for the product of transmetric spaces; the canonical projections in $\J\textup{et}'$ being the $\pi_i$, just as in $\J$et: 
see 1.2.2).
The fact that we have $can'(\pi_i)=can(\pi_i)=\kappa(p_i)$, the canonical projections in $\J$et$\jf$, provides the cartesian property of the functor $can'$.
\cqfd

\begin{prop}
$M,M_0,M_1$ being transmetric spaces, the following canonical 
map $can:\J\textup{et}\jf(M,M_0\times M_1)\lra
\J\textup{et}\jf(M,M_0)\times
\J\textup{et}\jf(M,M_1)$ is an isometry.
\end{prop}

\proof
Because, for $a\in M$, $a_0\in M_0$, $a_1\in M_1$, we have the
following commutative diagram (where $\J$et$_{(\ ,01)}$ and 
$\J$et$_{(\ ,i)}$ stand respectively for 
$\J$et$((M,a),(M_0\times M_1,(a_0,a_1)))$ and
$\J$et$((M,a),(M_i,a_i))$:
$$\xymatrix{
\J\textup{et}_{(\ ,01)}
\ar[d]_{can}
\ \ar[rr]^<<<<<<<<<<<<<<<<<<<<{can}&&\ \J\textup{et}_{(\ ,0)}
\times\J\textup{et}_{(\ ,1)}
\ar[d]^{can\times can}\\
\J\textup{et}\jf(M,M_0\times M_1)\ \ar[rr]_>>>>>>>>>>
{can}&&\ \J\textup{et}\jf(M,M_0)\times\J\textup{et}\jf(M,M_1)
}$$
The first horizontal $can$ refers to 1.2.22 and the vertical ones
to 1.4.15; we use the fact that the functor $can'$ of
1.4.25 is cartesian.
\cqfd

\begin{prop}
Let $M$ be a transmetric space and $U$ a non empty open subset of $M$.
Then, $\kappa(j):U\lra M$ is an isomorphism in $\J\textup{et}\jf$
(where $j:U\hookrightarrow E$ is the canonical injection).
\end{prop}

\proof
We have seen in 1.4.9 that $j$ is a morphism of
transmetric spaces; so that $\kappa(j)$ is a well defined
free metric jet.
Furthermore, referring to 1.2.5 and 1.3.3, we know that, for every $a\in U$,
T$j_a=q(j)=j_a:(U,a)\lra(M,a)$ is an isomorphism
in $\J$et and in $\J\text{et}'$, so that $\kappa(j)=can'(j_a)$ is itself
an isomorphism in $\J\textup{et}\jf$.
\cqfd

\begin{prop}
Let $M$, $M'$ be two transmetric spaces, and $U$, $U'$ two non empty open subsets of $M$ and $M'$ respectively.
Then, the canonical map
$\Gamma:\J\textup{et}\jf(U,U')\lra\J\textup{et}\jf(M,M'):
\Phi\mapsto \kappa(j').\Phi.\kappa(j)^{-1}$ is an isometry
(where $j:U\hookrightarrow M$ and $j':U'\hookrightarrow M'$
are the canonical injections).
\end{prop}

\proof
Because, for $a\in U$ and $a'\in U'$, the following diagram commutes (since $can'$ is a functor):
$$\xymatrix{
\J\textup{et}((U,a),(U',a'))\ar[d]_{can}\ar[rr]^{\Gamma}&&
\J\textup{et}((M,a),(M',a'))\ar[d]^{can}\\
\J\textup{et}\jf(U,U')\ar[rr]_{\Gamma}&&
\J\textup{et}\jf(M,M')
}$$
where the first horizontal $\Gamma$ refers to 1.2.24
and the vertical $can$ to 1.4.15. 
\cqfd

\begin{prop}
Let $M,M'$ be two transmetric spaces.
Then the map $J(M,M')\lra\R_+:(\varphi,a,a')\mapsto\rho(\varphi)$
factors through the quotient; we still denote
$\rho:\J\textup{et}\jf(M,M')\lra\R_+$ this factorization.
(see 1.2.8 for the definition of the lipschitzian ratio of a jet).
\end{prop}

\proof
Consider $(\varphi,a,a'),(\psi,b,b')\in J(M,M')$ such that
$(\varphi,a,a')\sim(\psi,b,b'$); we thus have
$\psi=\gamma(a',b').\varphi.\gamma(b,a)$, so that (using 1.2.9)
$\rho(\psi)\leq \rho(\gamma(a',b')).\rho(\varphi).\rho(\gamma(b,a))
\leq \rho(\varphi)$, since $\gamma(a',b')$ and $\gamma(b,a)$
are 1-bounded jets. As symetrically, we have also
$\rho(\varphi)\leq\rho(\psi)$, we obtain $\rho(\varphi)=\rho(\psi)$.
\cqfd

\begin{prop}
Let $M_0,M_1,M_2$ be three transmetric spaces; then,
if $\Phi_0:M_0\lra M_1$ and $\Phi_1:M_1\lra M_2$ are 
free metric jets, we have $\rho(\Phi_1.\Phi_0)\leq
\rho(\Phi_1)\rho(\Phi_0)$.
\end{prop}

\proof
Comes from 1.2.9.
\cqfd

\begin{prop}
Let $M_0,M_1,M_2$ be transmetric spaces;  we set\break
$Jeg\jf(M_1,M_2)=\{\Phi\in\J\textup{et}\jf(M_1,M_2)\,|\,
d(\Phi,O)=\rho(\Phi)\}$. Then the restriction 
$\J\textup{et}\jf(M_0,M_1)\times Jeg\jf(M_1,M_2)\
\buildrel{comp}\over\lra\ \J\textup{et}\jf(M_0,M_2)$
is $LL$.
\end{prop}

\proof
One can refer to 1.4.14 for the definition of the free metric jet $O$, and to 1.4.13 for the definition of the distance
between free metric jets; and also to 1.2.21 for the definition of
$Jeg((M_1,a_1),(M_2,a_2))$, when $a_1\in  M_1,a_2\in M_2$.

Now, the canonical map defined in 1.4.15 gives a restriction\break
$can:Jeg((M_1,a_1),(M_2,a_2))\lra Jeg\jf(M_1,M_2):\varphi
\mapsto[\varphi,a_1,a_2]$; indeed,
$\rho([\varphi,a_1,a_2])=\rho(\varphi)=d(\varphi,O_{a_1a_2})=
d([\varphi,a_1,a_2],[O_{a_1,a_2},a_1,a_2])\break
=d([\varphi,a_1,a_2],O)$, so that $can(\varphi)\in  
Jeg\jf(M_1,M_2)$.
This map is bijective:
if $\Phi=[\varphi, b_1,b_2]\in Jeg\jf(M_1,M_2)$, we 
know (referring to 1.4.15) that
$\Phi=can (\varphi')$ where $\varphi'=\gamma(b_2,a_2).\varphi.\gamma(a_1,b_1)$.
By similar arguments as before, we show that
$\varphi'\in Jeg((M_1,a_1),(M_2,a_2))$.
Then, we use the result of 1.2.21 and
the commutativity of the following diagram to obtain our result
(where, $\J$et$_{ij}$ and $Jeg_{ij}$ stand respectively for
$\J\textup{et}((M_i,a_i),(M_j,a_j))$ and
$Jeg((M_i,a_i),(M_j,a_j))$):
$$\xymatrix{
\J\textup{et}_{01}\times Jeg_{12}
\ar[d]_{can\times can}
\ \ar[rr]^<<<<<<<<<<<<<<<<<<<<{comp}&&\ \J\textup{et}_{02}
\ar[d]^{can}\\
\J\textup{et}\jf(M_0,M_1)\times Jeg\jf(M_1,M_2)\ 
\ar[rr]_>>>>>>>>>>
{comp}&&\ \J\textup{et}\jf(M_0,M_2)
}$$
\cqfd

\vspace{3mm}
Again, we conclude this paragraph with vectorial considerations.

\begin{prop}
Let $M$ be a transmetric space and $E$ a n.v.s..
Then, the
metric space $\J\textup{et}\jf(M,E)$ has also 
a canonical structure of n.v.s..
(its distance defined in 1.4.13
derives from a norm). Besides, for every $a\in M$, the map
$can:\J\textup{et}((M,a),(E,0))\lra\J\textup{et}\jf(M,E)$ 
(see 1.4.15) is a
linear isometry.
\end{prop}

\proof
We have just to transfer the structure of n.v.s. from\break
$\J$et$((M,a),(E,0))$ to $\J$et$\jf(M,E)$, via the
isometry $can$;
there is no problem since this definition does not depend on $a$
(thanks to the commutativity of the diagram of 1.4.16, in which we set $M'=E$ and $a'=b'=0$; we use the linearity of $\Omega$ which
comes from 1.2.26). Let us notice that the ``zero'' of 
$\J$et$\jf(M,E)$ is the free zero $O$ defined in 1.4.14.
\cqfd

\begin{prop}
Let $M$ be a transmetric space and $E$ a n.v.s..
Let also $\Phi\in \J\textup{et}\jf(M',M)$.
Then, the following  map $\widetilde\Phi$
is linear and continuous, where 
$\widetilde\Phi:\J\textup{et}\jf(M,E)\lra
\J\textup{et}\jf(M',E):\Psi\mapsto\Psi.\Phi$.
\end{prop}

\proof
We just have to use the commutativity of the following diagram (with
the help of 1.2.26 and 1.4.32; where the jet $\varphi$ arises from the definition of $\Phi=[\varphi,a,a']$):
$$\xymatrix{
\J\textup{et}((M,a),(E,0))\ar[rr]^{\widetilde\varphi}
\ar[d]_{can}&&\J\textup{et}((M',a'),(E,0))
\ar[d]^{can}\\
\J\textup{et}\jf(M,E)\ar[rr]_{\widetilde\Phi}&&
\J\textup{et}\jf(M',E)
}$$
\cqfd

%% file: jet_I3.tex
\section{Tangential}
Untill now, we have only spoken of the tangential T$f_a$, at $a$, of a 
map $f:M\lra M'$, tangentiable at $a$ (see section 1.3). 
From now on, as foretold in the introduction of section 1.4, we will be able to speak of the tangential 
t$f:M\lra\J\textup{et}\jf(M,M')$ when $f:M\lra M'$ is a tangentiable map between transmetric spaces
(as we speak of the differential d$f:U\lra L(E,E')$ for  a differentiable
map $f:U\lra E'$ when $U$ is an open subset of $E$, an e.v.n., 
just as $E'$).

\vspace{3mm}

First, consider $M,M'$ two transmetric spaces and $f:M\lra M'$ a tangentiable 
map at the point $x\in M$; then,  we set t$f_x=
[$T$f_x,x,f(x)]$.

Now, if $f:M\lra M'$ is 
tangentiable (at every point in $M$),
we can define the map
t$f:M\lra\J\textup{et}\jf(M,M'):x\mapsto $t$f_x$
(in fact, this map t$f$ is the following composite
$M\buildrel{\textup{T}f}\over\lra J(M,M')\buildrel q\over\lra
\J\textup{et}\jf(M,M')$, 
where $J(M,M')$, $\J\textup{et}\jf(M,M')$ and $q$ have been defined round
1.4.11, and where T$f(x)=($T$f_x,x,f(x))$.

\begin{defi}
The map $\textup{t}f$ defined just above (for a tangentiable map $f$ between
transmetric spaces) will be called the 
\textup{tangential} of $f$.
\end{defi}

\begin{prop}
Let $f:M\lra M'$ be a tangentiable map between transmetric spaces.
Then, the map $\textup{t}f:M\lra\J\textup{et}\jf(M,M')$ is constant iff
$f$ is a morphism of transmetric spaces; in this case
we have\break
$\textup{t}f_x=\kappa(f)$ for every $x\in M$.  
\end{prop}

\proof
See 1.4.8 and 1.4.22.
\cqfd

\vspace{3mm} 
Let us now begin with the particular vectorial context:
we assume here that $E,E'$ are n.v.s. and $U$ a non empty open subset of $E$; consider then the following composite
(denoted $J$):
$$\xymatrix{
L(E,E')\ar[r]^>>>>>{j} &\J\textup{et}((E,0),(E',0))
\ar[r]^<<<<{can} & \J\textup{et}\jf(E,E')
\ar[r]^{\Gamma^{-1}} & \J\textup{et}\jf(U,E')
}$$
See 1.2.27, 1.4.15 and 1.4.28 for the definitions of 
$j$, $can$ and 
$\Gamma$;
for every $l\in L(E,E')$, 
we have $J(l)=[q(l),0,0].\kappa(j)$
(where 
$j:U\hookrightarrow E$ is the canonical injection : see 1.4.27), and since, a priori, $0$ is not necessarilly in $U$, 
we take any $a\in U$, to write $\kappa(j)=[j_a,a,a]$ so that 
$J(l)=[q(l).\gamma(a,0).j_a,a,0]$.
We notice that $J(l)$ does not depend on $a\in U$, just like $\kappa(j)$. 
We set $Im(J)=J(L(E,E'))$, the image of $J$ in $\J$et$\jf(U,E')$.

\begin{prop}
The map $J:L(E,E')\lra \J\textup{et}\jf(U,E')$ defined just above is 
a linear isometric embedding.
\end{prop}

\proof
It comes from
1.2.27, 1.4.32, 1.4.28 and 1.4.33. 
\cqfd

\begin{prop}
Let $f:U\lra E'$ be a tangentiable map and $a\in U$; then, 
$f$ is differentiable at $a\,$ iff $\ \textup{t}f_a\in Im(J)$;
in this case, we have $\textup{t}f_a=J(\textup{d}f_a)$.
\end{prop}

The proof of this proposition uses the result of the following lemma:  

\begin{lema}
Let $[\varphi,a,b]\in\J\textup{et}\jf(U,E')$.
Then, $[\varphi,a,b]\in Im(J)$ iff there exists $l\in L(E,E')$
such that $Al|_U\in\varphi$, where $Al(x)=b+l(x-a)$;
and we have $J(l)=[\varphi,a,b]$.
\end{lema}

\proof
$[\varphi,a,b]\in Im(J)$ iff there exists $l\in L(E,E')$
such that $[q(l).\gamma(a,0).j_a,a,0]=J(l)=[\varphi,a,b]$
iff  $q(l).\gamma(a,0).j_a=\gamma(b,0).\varphi$ iff
$q(\theta(0,b).l.\theta(a,0).j)=\varphi$
(since, in the vectorial case, each $\gamma(a,b)$ 
is the jet of the translation $\theta(a,b)$ defined
in 1.4.5).
It remains just to notice that the map $\theta(0,b).l.\theta(a,0)$
is nothing but the map $Al$. 
\cqfd
\vspace{1mm}
We now come back to the proof of 1.5.4:

\proof
Since t$f_a=[$T$f_a,a,f(a)]$, we have: t$f_a\in Im(J)$ iff
there exists $l\in L(E,E')$ such that $Al|_U\in$T$f_a$,
i.e. such that $Al|_U\tang_a f$, which means that $f$ is differentiable
at $a$ with $l=$d$f_a$. 
\cqfd

\begin{cory}
Let $f:U\lra E'$ a differentiable map; then $f$ is tangentiable
and the following diagram commutes:
$$\xymatrix{
{}&U\ar[dl]_{\textup{d}f}\ar[dr]^{\textup{t}f}&{}\\
L(E,E')\ar[rr]_J && \J\textup{et}\jf(U,E')
}$$
\end{cory}

\vspace{3mm}
We now come back to the general transmetric context.

\begin{prop}
Let $M,M',M''$ be three transmetric spaces;\break
$f:M\lra M'$, $g:M'\lra M''$ two maps.
We assume that $f$ and $g$ are tangentiable and that 
their tangential $\textup{t}f:M\lra\J\textup{et}\jf(M,M')$ and
$\textup{t}g:M'\lra\J\textup{et}\jf(M',M'')$ are continuous, respectively at
$a\in M$ and $a'=f(a)\in M'$. Then the map
$\textup{t}(g.f):M\lra\J\textup{et}\jf(M,M'')$ is well defined and continuous at $a$, and we have $\textup{t}(g.f)_a=
\textup{t}g_{a'}.\textup{t}f_a$.
\end{prop}

\proof
We know (see 1.3.6) that $g.f$ is tangentiable and that, for every
$x\in M$, we have T$(g.f)_x=$T$g_{f(x)}.$T$f_x$; which gives
t$(g.f)_x=\break
$t$g_{f(x)}.$t$f_x$,
using the composition of the free metric jets (see 1.4.17).
In fact, the map t$(g.f)$ is the following composite:
$$\xymatrix{
M\ar[rr]^>>>>>>>>>>{(Id,f)}&& M\times M'
\ar[rr]^>>>>>>>>>>{\textup{t}f\times \textup{t}g}
&&\J\textup{et}\jf(M,M')\times\J\textup{et}\jf(M',M'')
\ar[d]^{comp}\\
{}&&&&\J\textup{et}\jf(M,M'')
}$$
which gives the result, thanks to 1.1.8, 1.3.2
and 1.4.18.
\cqfd

\begin{prop}
Let $M,M_0,M_1$ be transmetric spaces, with $a\in M$ and 
$f_0:M\lra M_0$, $f_1:M\lra M_1$ two maps.
We assume that, for each $i\in\{0,1\}$, $f_i:M\lra M_i$ is tangentiable and that  its tangential
\textup{t}$f_i:M\lra\J\textup{et}\jf(M,M_i)$ is continuous at $a$.
Then, the map\break
$\textup{t}(f_0,f_1):M\lra\J\textup{et}\jf(M,M_0\times M_1)$
is well defined and continuous at $a$, and we have
$\textup{t}(f_0,f_1)_a=(\textup{t}f_{0a},\textup{t}f_{1a})$.
\end{prop}

\proof
We know (see 1.3.7) that $(f_0,f_1)$ is tangentiable and that, for every $a\in M$, we have T$(f_0,f_1)_a=($T$f_{0a},$T$f_{1a})$; 
which gives t$(f_0,f_1)_a=($t$f_{0a},$t$f_{1a})$ (thanks to 1.4.24 and  1.4.25).
In fact, the map t$(f_0,f_1)$ is the following composite:
$$\xymatrix{
M\ar[rr]^>>>>>>>>>>{(\textup{t}f_0,\textup{t}f_1)}&& 
\J\textup{et}\jf(M,M_0)\times\J\textup{et}\jf(M,M_1)
\ar[r]^<<<<<{can^{-1}}
&\J\textup{et}\jf(M,M_0\times M_1)
}$$
which gives the wished result (thanks to 1.4.26).
\cqfd

\begin{prop}
Let $E$ be a n.v.s., $U$ an open subset of $E$, $a,b\in U$ such that $[a,b]\subset U$ and $F$ a finite subset of $]a,b[$; let also
$M$ be a transmetric space and $f:U\lra M$ a continuous map.
Now, if we assume that, for all $x\in\,]a,b[-F$, the map 
$f$ is $Tang_x$ and that 
$d(tf_x,O)\leq k$ (where $O$ is the free zero of 1.4.14 and $k$ a fixed positive real number), then we have
$d(f(b),f(a))\leq k||b-a||$.
\end{prop}

\proof
Let $\varepsilon>0$ and $x\in\,]a,b[-F$.
Since we have the equalities\break
$d($T$f_x,O_{xf(x)})=d([$T$f_x,x,f(x)],[O_{xf(x)},x,f(x)])=
d($t$f_x,O)<k+\varepsilon$, we know, by 1.3.12, that 
$f$ is $(k+\varepsilon)$-$LSL_x$. Using now 1.1.23, 
we obtain $d(f(b),f(a))\leq(k+\varepsilon)||b-a||$, thus the wished result (doing $\varepsilon\rightarrow 0$).
\cqfd

\vspace{3mm}
We now give a generalization of 2) in 1.1.12:

\begin{theo}
Let $E$ and $M$ be respectively a n.v.s. and a transmetric space, $U$ a non empty open subset of $E$ and $f:U\lra M$ a tangentiable map such that $\textup{t}f:U\lra\J\textup{et}\jf(U,M)$ is continuous at $a\in U$.
Then,
$f$ is $LL_a$ and we have $d(\textup{t}f_a,O)=
\rho(\textup{t}f_a)$; i.e $\textup{t}f_a$ is a good free metric jet (see 1.2.11
and 1.4.31).
\end{theo}

\proof
First, thanks to 1.2.10 and
by definition of the lipschitzian ratio $\rho$ for linked or free metric jets (see 1.2.8 and 1.4.29), we have
$d($t$f_a,O)=d($T$f_a,O_{af(a)})\leq \rho($T$f_a)=\rho($t$f_a)$, the first equality resulting from the fact that the map $can$ of 1.4.15 is an isometry.

Now, let $\bar f:U\lra \R$ be the following composite:
$$\xymatrix{
U\ar[r]^>>>>>{\textup{t}f}&\J\textup{et}\jf(U,M)
\ar[rr]^<<<<<<<<<<{d(-,O)}&&\R
}$$
The map $\bar f$ is continuous at $a$, by composition.
Let us fix $\varepsilon>0$; then, there exists $r>0$ such that $B(a,r)\subset U$ and $\bar f(x)<\bar f(a)+\varepsilon$ for all $x\in B(a,r)$.
Let us set $R=\bar f(a)+\varepsilon$;
then, for all $x\in B(a,r)$, we have $d($t$f_x,O)=\bar f(x)<R$,
so that (thanks to 1.5.9), the restriction $f|_{B(a,r)}$ is
$R$-lipschitzian (since $B(a,r)$ is convex).
Thus $f$ is $LL_a$, so that T$f_a=q(f)$ (see 1.3.3). Furthermore, 
we have $\rho($t$f_a)=\rho($T$f_a)=
\rho(q(f))\leq R=\bar f(a)+\varepsilon$;
we have thus obtained $\rho($t$f_a)\leq d($t$f_a,O)+\varepsilon$ for all $\varepsilon>0$, so that $\rho($t$f_a)\leq d($t$f_a,O)$.
\cqfd

\begin{prop}
Let $E$ be a n.v.s. and $U$ a non empty convex open subset of $E$, and $M$ a transmetric space.
Then, every morphism of transmetric spaces $f:U\lra M$ is $R$-lipschitzian
(where $R=d(\kappa(f),O)$; see 1.4.22 for $\kappa(f)$). 
Hence, $d(\kappa(f),O)=\rho(\kappa(f))$; i.e $\kappa(f)$ is a good free
metric jet (see 1.2.11 and 1.4.31).
\end{prop}

\proof
We use 1.5.9, since, for all $x\in U$, we have 
$d($t$f_x,O)=d(\kappa(f),O)=R$; thus $f$ is $R$-lipschitzian on $U$
(which is convex), and thus $\rho(\kappa(f))\leq R$. For the 
inverse inequality, we use  1.2.10.!
\cqfd

\begin{cory}
Under the same notations and hypothesis of 1.5.11, except that we don't assume here that $U$ is convex, then $f$ is $R$-$LL$ on 
$U$,
and $\kappa(f)$ is still a good free metric jet.
\end{cory}

\proof
Because every open subset of a n.v.s. is locally convex.
\cqfd

\vspace{3mm}
\begin{remk}
Referring to 1.4.9, 1.5.11 is true for every continuous affine map.
\end{remk}
\vspace{3mm} 

\begin{defi}
Let $M,M'$ be transmetric spaces and $f:M\lra M'$ a map.
We say that $f$ is \textup{continuously tangentiable} if $f$ is tangentiable and if 
its tangential $\textup{t}f:M\lra\J\textup{et}\jf(M,M')$ is continuous. 
\end{defi}

\begin{prop}
The continuously tangentiable maps are stable\break
under composition and under pairs.
\end{prop}

\proof
See 1.5.7 and 1.5.8.
\cqfd

\vspace{3mm}
\begin{exams}
{}
\end{exams}

1) Every morphism of transmetric spaces is continuously tangentiable (see 1.5.2).

2) Every map which is of class $C^1$ is continuously tangentiable (by 1.5.3 and 1.5.6).

\vspace{3mm}
\begin{theo}
Let $E,E'$ be n.v.s. where $E'$ is complete, $U$ a non empty open subset of $E$
and $f:U\lra E'$ a continuously tangentiable map which is also differentiable
at every point of a dense subset $D$ of $U$.
Then, $f$ is of class $C^1$.
\end{theo}

\proof
Referring to 1.5.4, we know that, a priori, the set of the points of $U$ at which $f$ is differentiable is nothing but 
$\buildrel{-1}\over{\textup{t}f}(Im(J))$; thus, we must first
show that this set is equal to $U$. Since, by hypothesis, we have the inclusions
$D\subset\ \buildrel{-1}\over{\textup{t}f}(Im(J))\subset U$ with $D$ dense in $U$, we have just to prove that 
$\buildrel{-1}\over{\textup{t}f}(Im(J))$ is closed in $U$.
Actually, since $Im(J)$ is complete (just as the n.v.s. $L(E,E')$, since
$J:L(E,E')\lra \J$et$\jf(U,E')$ is an isometric embedding:
see 1.5.3), it is closed in  $ \J$et$\jf(U,E')$, so that 
$\buildrel{-1}\over{\textup{t}f}(Im(J))$ is also closed in $U$ (since t$f:U\lra\J$et$\jf(U,E')$ is continuous);
$f$ is thus differentiable on $U$. 
Finally, the fact that the map
d$f:U\lra L(E,E')$ is continuous comes from
the fact it factors through t$f$ (see the commutative triangle in 1.5.6).
\cqfd

\vspace{3mm}
\begin{exam}
{}
\end{exam}
The function $\vartheta:\R\lra\R:x\mapsto |x|$ given in 1.3.9 is tangentiable (for it is lipschitzian);
however, it cannot be continuously tangentiable on $\R$, since it is differentiable on 
$\R^*$ and not of class $C^1$ on $\R$.

\vspace{3mm}
We now give a generalization of the fact that the composite of two maps of class $C^2$ is also of class $C^2$.

\begin{theo}
Let $M,M'$ be transmetric spaces, $E$ a n.v.s. and $U$ a non empty open subset of $E$.
Let also $f:M\lra U$ and $g:U\lra M'$\break
be two continuously tangentiable maps.
We assume that $\textup{t}f:\break
M\lra\J\textup{et}\jf(M,U)$ is $Tang_a$ where $a\in M$ and
$\textup{t}g:U\lra\J\textup{et}\jf(U,M')$ is $LL_{a'}$ where $a'=f(a)\in U$.
Then, $\textup{t}(g.f):M\lra\J\textup{et}\jf(M,M')$ is $Tang_a$.
\end{theo}

\proof
Since $g:U\lra M'$ is continuously tangentiable, we have (by 1.5.10) 
$d($t$g_y,O)=\rho($t$g_y)$ for all $y\in U$, which means that\break
t$g_y\in Jeg\jf(U,M')$ (see 1.4.31).We can thus consider the restriction
$\breve{\textup{t}g}:U\lra Jeg\jf(U,M')$ of t$g$.
Besides, by 1.5.7, we know that $g.f$ is continuously tangentiable
with t$(g.f)_x=$t$g_{f(x)}.$t$f_x$ for all $x\in M$; 
thus, we can write its tangential
t$(g.f)$ as the composite:
$$\xymatrix{ 
M\ar[rr]^>>>>>>>>>>>{(Id,f)}&& M\times U\ar[rr]^>>>>>>>>>>>>
{\textup{t}f\times\breve{\textup{t}g}}&&
\J\textup{et}\jf(M,U)\times Jeg\jf(U,M')
\ar[d]^{comp}\\
&&&& \J\textup{et}\jf(M,M')
}$$
Now, by 1.3.6, we only have to justify the fact that each forementionned arrow is tangentiable at the right point. Now, $(Id,f)$ is $Tang_a$ (see 1.3.7);
t$f\times\breve{\textup{t}g}$ is $Tang_{(a,a')}$ (see 1.3.2 and 1.3.8) and
$comp$ is $LL_{(\textup{t}f_a,\textup{t}g_{a'})}$ (see 1.4.31), and thus $Tang_{(\textup{t}f_a,\textup{t}g_{a'})}$.
Their composite t$(g.f)$ is thus $Tang_a$ as foretold in the theorem.
\cqfd
\vspace{12mm}
\noindent {\texttt{
Here are some open problems that we submit to the sagacity of the
reader (after having gone through this first chapter)}}
\vspace{3mm}

1) We have seen, in the examples given in 1.4.9, that every continuous affine map between n.v.s. is a morphism of transmetric spaces.
The question is: are there other morphisms of transmetric spaces in such a 
vectorial context?

2) Still in the case of the n.v.s., does it exist continuous tangentiable maps which are not of class $C^1$ (see 1.5.16 and 1.5.17)?

%% file: jet_II1.tex
\chapter{Representation of metric jets}
 
\section{Valued monoids,
contracting spaces}
We have already mentioned (in the initial abstract) that the unique 
continuous affine map of 
the tangential T$f_a$ of a map $f$, differentiable at $a$,
could be considered as canonically
representing this linked metric jet. Our aim, in this section, is
precisely to enlarge the domain of application of this uniqueness property.
We begin with giving an algebraico-metric context 
which is rich enough to obtain this uniqueness property.
\vspace{3mm}
\begin{defi}
a \textup{valued monoid} $\Sigma$ is a monoid (its law of composition 
is denoted here multiplicatively) possessing an absorbing element 0 
(it verifies $t.0=0$ for each $t\in \Sigma$) which is, in addition, equipped with a monoid homomorphism $v:\Sigma\lra\R_+$
(where $\R_+=[0,+\infty[$), called
the \textup{valuation} of $\Sigma$, verifying the two following conditions:

(1) $\forall t\in\Sigma\quad(v(t)=0\Llra t=0)$,

(2) $\exists t\in\Sigma\quad(0<v(t)<1)$.

\noindent Thanks to the above (1), we have, for all $t,t'\in\Sigma$, the implication:
$t.t'=0\Lra t=0$ or $t'=0$.
\end{defi}

\vfill\eject

\begin{exams}
{}
\end{exams}

1) $\R$ and its multiplicative submonoids $\R_+$ and 
$[0,1]$ are valued monoids, where $0$ is their absorbing element (their valuation being the absolute
value).

2) If $r$ is a real number verifying $0<r<1$, we denote $\N'_r$
the additive monoid $\N'_r=\N\cup\{\infty\}$, where $\infty$ is its
absorbing element, and whose valuation $v_r:\N'_r\lra\R_+$ is given
by $v_r(n)=r^n$ if $n\in\N$ and $v_r(\infty)=0$.

3) The previous examples are commutative valued monoids; here is a non commutative one: let us consider $\mathcal S$im$(E)$, the set of the ``similitudes'' on a given n.v.s $E$ (its elements are the $u=\lambda\varphi$ where $\lambda\in \R$ and $\varphi:E\lra E$ is an isometry).
$\mathcal S$im$(E)$ is a non commutative monoid (for its natural law of composition); it's valuation being $v(u)=\|u\|$ (the operator norm; in fact, $v(u)=|\lambda|$, if $u=\lambda\varphi$).

\vspace{3mm}
\begin{defi}
A \textup{morphism of valued monoids} $\sigma:\Sigma\lra\Sigma'$ is a\break
monoid homomorphism verifying the two following conditions:

(1) $\forall t\in\Sigma\quad(v(\sigma(t))=v(t))$,

(2) $\sigma(0)=0$.
\vspace{1mm}

A \textup{valued submonoid} of $\Sigma'$ is a
valued monoid $\Sigma$ verifying $\Sigma\subset\Sigma'$ such that the canonical injection $j:\Sigma\hookrightarrow\Sigma'$ is a 
morphism of valued monoids; it is the case of 
the canonical injection $j:\R_+\hookrightarrow\R$.
\end{defi}

\begin{prop}
If $\sigma:\Sigma\lra\Sigma'$ is a morphism of valued monoids, we have the equivalence: $\sigma(t)=0\Llra t=0$.
\end{prop}
\proof
Because $0=v(\sigma(t))=v(t)\Lra t=0$.
\cqfd

\vspace{3mm}
\begin{exam}
{}
\end{exam}

A valuation $v:\Sigma\lra\R_+$ is itself a morphism of valued monoids.
\vspace{3mm}

\begin{defi}
$\Sigma$ being a valued monoid, a $\Sigma$-\textup{contracting
metric space} (in short a $\Sigma$-contracting space) is a metric space $M$ equipped with a
\textup{central point} denoted $\omega$ and with an external operation $\Sigma\times M\lra M:(t,x)\mapsto t\star x$ which,
in addition of the usual properties:

(1) $\forall x\in M\,$ $(1\star x=x)$,

(2) $\forall t,t'\in\Sigma\ \forall x\in M\,$ $(t'\star(t\star x)=(t'.t)\star x)$,

\noindent verifies also the following conditions:

(3) $0\star \omega=\omega$,

(4) $\forall t\in\Sigma$ $\,\forall x,y\in M\ 
(d(t\star x,t\star y)=v(t)d(x,y))$

\noindent The central point of any generic $\Sigma$-contracting space
will always be denoted $\omega$.
\end{defi}

\begin{prop}
A $\Sigma$-contracting space $M$ verifies the following properties :

(1) $\forall t\in\Sigma\,$ $(t\star\omega=\omega$),

(2) $\forall x\in M\,$ $(0\star x=\omega$).
\end{prop}

\proof
1) $t\star\omega=t\star(0\star\omega)=(t.0)\star\omega=0\star\omega=\omega$.

2) $d(0\star x,\omega)=d(0\star x,0\star\omega)=v(0)d(x,\omega)=0$.
\cqfd

\vspace{3mm}
\begin{exams}
{}
\end{exams}

1) Let $E$ be a n.v.s.; fixing a point $a\in E$, it defines a pointed n.v.s. $(E,a)$ denoted $E_a$. The point $a$ being considered as the central point,
we easely make this $E_a$ a 
$\R$-contracting space, setting,
for $t\in\R$ and $x\in M$, $t\star x=a+t(x-a)$;
which gives $t\star x=tx$, when $a=0$.
If we restrict the scalars to be in $\R_+$ (resp. in $[0,1]$), 
$E_a$ is also vewed as a $\R_+$-contracting space
(resp. a $[0,1]$-contracting space).
These external operations on $E_a$ are said to be ``standard''.

2) Consider the non commutative valued monoid $\mathcal S$im$(E)$, 
whose elements are the similitudes on a given n.v.s $E$ (see 2.1.2). Then $E$ is a $\mathcal S$im$(E)$-contracting space,
with $u\star x=u(x)$ and $\omega=0$.

3) When $\Sigma$ is a valued monoid whose valuation 
$v:\Sigma\lra\R_+$ is injective, then $\Sigma$ becomes
itself a $\Sigma$-contracting space setting $\omega=0$ and, for $s,t\in\Sigma$,
$t\star s=t.s$ and $d(s,t)=|v(t)-v(s)|$. It is the case for 
$\Sigma=[0,1]$ and $\Sigma=\N'_r$ (see 2.1.2).

4) If $M,M'$ are $\Sigma$-contacting spaces, then so is $M\times M'$. 

\vspace{3mm}

\begin{prop} \textup{(Restriction of the scalars)}

Let $\sigma:\Sigma\lra\Sigma'$ be a morphism of valued monoids;
then every $\Sigma'$-contracting space can be canonically equipped with a structure of\break
$\Sigma$-contracting space.
\end{prop}

\proof
The new $\Sigma$-operation is $t\star x=\sigma(t)\star' x$.
\cqfd

\vspace{3mm}
\begin{remks}\par\hfill

\noindent From now on, we will denote every $\Sigma$-operation by the same symbol $\star$.
 
1) Referring to 2.1.5 and 2.1.9, we know that, if $\Sigma$ is a valued monoid,  every $\R_+$-contracting space $M$ has also 
a ``canonical'' structure of
$\Sigma$-contracting space (its external operation being $t\star x=
v(t)\star x$ for every $t\in\Sigma$ and $x\in M$). 
In particular, if $E$ is a n.v.s. (with a fixed point $a$) : $E_a$ (see 2.1.8)
has also a canonical structure of
$\Sigma$-contracting space (with $\omega=a$,
its external operation being, for $t\in\Sigma$ and $x\in E$, 
$t\star x=a+v(t)(x-a)$); for instance $E_a$ 
has a canonical structure of
$\N'_r$-contracting space 
(its canonical external operation being, for $n\in\N$ and $x\in E$, 
$n\star x=a+r^n(x-a)$ and $\infty\star x=a=\omega$).

2) We notice that we have two structures of $\R$-contracting space 
on $E_a$: the standard one (see 2.1.8) and the above canonical one,
whose external operation is $t\star x=a+|t|(x-a)$. 
Actually, enlightened by the point 7) in 2.1.17, it is the standard 
$\R$-contracting structure on $E_a$ which impose itself to be the most natural.

3) In a $\Sigma$-contracting space $M\not=\{\omega\}$, the point $\omega$
is never isolated in $M$ : indeed, if we choose $\,t\in\Sigma$ verifying $0<v(t)<1$, then we have, for all $x\in M$,
$d(t^n\star x,\omega)=d(t^n\star x,t^n\star\omega)=v(t^n)d(x,\omega)=
(v(t))^nd(x,\omega)$, which implies, on the one hand, $t^n\star x\not=\omega$ for all $n\in\N$ when $x\not=\omega$ 
and, on the other hand,
$\omega=\lim_n t^n\star x$ (since $\lim_n(v(t))^n=0$).

4) Why the term ``contracting''? Let us briefly say that it reminds the contracting maps (i.e $k$-lipschitzian maps with $0<k<1$; it is the case of the maps $x\mapsto t\star x$
when
$v(t) <1$) and the contractible spaces (indeed, when $\Sigma=[0,1]$, the external operation $[0,1]\times M\lra M$ verifies
$0\star x=\omega$ and $1\star x=x$ ... which is not very far from an homotopy: $\widehat\omega\simeq Id_M$).
\end{remks}

\vspace{3mm}
\begin{defi}\par\hfill 

1) We say that a valued monoid is a \textup{quasi-group}
if every $t\in\Sigma$ verifying $t\not=0$ is invertible.

2) We say that a $\Sigma$-contracting space is \textup{revertible} if,
for every $t\in\Sigma$ verifying $t\not=0$, the map
$t\star(-):M\lra M$ is bijective. In this case, we set
$\ t\buildrel{-1}\over\star x=(t\star(-))^{-1}(x)\ $ for $x\in M$.
\end{defi}

\begin{prop}
When $\Sigma$ is a quasi-group, every $\Sigma$-contracting space is revertible, with $t\buildrel{-1}\over\star x=t^{-1}\star x$ for every
$t\in\Sigma$ ($t\not=0$) and $x\in M$.
\end{prop}

\proof
Since $t^{-1}\star (t\star x)=(t^{-1}.t)\star x=1\star x=x$.
\cqfd

\begin{prop}
If $\sigma:\Sigma\lra\Sigma'$ is a morphism of valued monoids and 
$M$ a revertible $\Sigma'$-contracting space, then $M$ is also
a revertible\break
$\Sigma$-contracting space.
\end{prop}

\proof
For all $t\not=0$, we have $\ t\buildrel{-1}\over{\star} x=
\sigma(t)\buildrel{-1}\over\star x$.
\cqfd

\vspace{3mm}
\begin{examsc}
{}
\end{examsc}

1) $\R$ and $\R_+$ being (multiplicative) quasi-goups, every 
$\R$-contracting space or $\R_+$-contracting space is revertible.
(for their standard structure ... see 2.1.10).

2) If $E$ is a n.v.s., and $a\in E$, then $E_a$ is a 
revertible
$\R_+$-contracting space (resp. a revertible $\R$-contracting space)
for the standard structures,
with
$t\buildrel{-1}\over\star x=a+t^{-1}(x-a)$ for $t\in\R_+$ (or $t\in\R$), 
$t\not=0$.
Referring to 2.1.5 and 2.1.13, we know that, for every valued monoid $\Sigma$,  $E_a$ is also a
revertible $\Sigma$-contracting space for the canonical structure (with
$t\buildrel{-1}\over\star x=a+v(t)^{-1}(x-a)$ for $t\in\Sigma$, $t\not=0$).

3) Referring to 2.1.2 and 2.1.8, we know that $[0,1]$ and $\N'_r$ are\break
$\Sigma$- contracting spaces with, respectively, $\Sigma=[0,1]$, $\omega=0$,
and $\Sigma=\N'_r$, $\omega=\infty$ (these valued monoids being
not quasi-groups!). But none of them is revertible since, for the first one 
$\frac{1}{2}\buildrel{-1}\over\star x=2x\notin[0,1]$ for 
$\frac{1}{2}<x\leq 1$, and, for the second one
$n\buildrel{-1}\over\star m=m-n\notin\N'_r\,$ for $m<n$. 
\vspace{3mm}

\begin{prop}
Let $M$ be a revertible $\Sigma$-contracting space, then, if we set
$\Sigma^*=\Sigma-\{0\}$, we have:

(1) $\, \forall t\in\Sigma^*\ (t\buildrel{-1}\over\star\omega=\omega)$,

(2) $\,\forall t,t'\in\Sigma^*\ \forall x\in M\
(t\buildrel{-1}\over\star(t'\buildrel{-1}\over\star x)=
(t'.t)\buildrel{-1}\over\star x)$,

(3) $\,\forall t\in\Sigma^*\ \forall x,y\in M\
(d(t\buildrel{-1}\over\star x,t\buildrel{-1}\over\star y)=
\frac{1}{v(t)}d(x,y))$,

(4) $\,\forall t\in\Sigma^*\ \forall x,y\in M\ 
(t\star x=y\Llra x=t\buildrel{-1}\over\star y)$.
\end{prop}

\proof
(1): Since, here, $t\star(-)$ is injective, we have just to notice that
$t\star(t\buildrel{-1}\over\star\omega)=\omega=t\star\omega$,

(2) Since $t'.t\not=0$ (see 2.1.1), we write 
$(t'.t)\buildrel{-1}\over\star (-)=
((t'.t)\star(-))^{-1}=
((t'\star(-)).(t\star(-))) ^{-1}=
(t\star(-))^{-1}.(t'\star(-)) ^{-1}=
(t\buildrel{-1}\over\star(-)).(t'\buildrel{-1}\over\star(-))$.

(3) $d(x,y)=d(t\star(t\buildrel{-1}\over\star x),t\star(t\buildrel{-1}\over\star y))=v(t)d(t\buildrel{-1}\over\star x,t\buildrel{-1}\over\star y)$;
it remains to use the fact that $v(t)\not=0$.

(4) By definition of $\buildrel{-1}\over\star$.
\cqfd

\vspace{3mm}
\begin{defi}
Let $\Sigma$ be a fixed valued monoid and $M,M'$ two\break
$\Sigma$-contracting spaces.

1) We say that a map $h:M\lra M'$ is 
$\Sigma$-\textup{homogeneous} if\break
$h$ verifies:
$\forall t\in\Sigma\ \,\forall x\in M\quad (h(t\star x)=t\star h(x))$.
From now on, when we will say that a map $h$ is 
$\Sigma$-homogeneous, 
it will implicitly mean that its domain and codomain are 
$\Sigma$-contracting spaces.

2) We say that $M$ is a $\Sigma$-\textup{contracting metric subspace} (in short a $\Sigma$-\textup{contracting subspace})
of $M'$ if $M$ is a
metric subspace of $M'$ and if the canonical injection
$M\hookrightarrow M'$ is $\Sigma$-homogeneous.
\end{defi}

\begin{prop}\par\hfill

1) A $\Sigma$-homogeneous map $h:M\lra M'$ verifies 
$h(\omega)=\omega$.

2) If $M,M'$ are revertible $\Sigma$-contracting spaces, then every\break
$\Sigma$-homogeneous map
$h:M\lra M'$ verifies $h(t\buildrel{-1}\over\star x)=t\buildrel{-1}\over\star
h(x)$ for all $t\in\Sigma^*$ and $x\in M$.

3) The $\Sigma$-homogeneous maps are stable under composition and pairs.

4) $M'$ being a $\Sigma$-contracting space and $M$ a metric subspace of $M'$,
then $M$ is a $\Sigma$-contracting subspace of $M'$ iff $\omega\in M$ and $t\star x\in M$ for all
$t\in \Sigma$ and $x\in  M$.

5) Let $\sigma:\Sigma\lra\Sigma'$ be a morphism of valued monoids, 
$M$ and $M'$ two $\Sigma'$-contracting spaces.
Then, every $\Sigma'$-homogeneous map $M\lra M'$ is 
$\Sigma$-homogeneous;
the inverse is true when $\sigma$ is surjective.

6) If $\Sigma'=\R_+$ and $\sigma=v$, then every 
$\R_+$-homogeneous map\break
 $M\lra M'$ 
is $\Sigma$-homogeneous 
(for the canonical $\Sigma$-contracting structures on $M$ and $M'$ defined in 2.1.10);
the inverse is true when $v$ is surjective.

7) In particular, if $E$ and $E'$ are n.v.s., and $a\in E$, a map\break 
$h:E_a\lra E'_{h(a)}$ is $\R_+$-homogeneous iff it is $\R$-homogeneous for the canonical $\R$-contracting structures
on $E_a$ and $E'_{f(a)}$. 
So nothing new as for the homogeneousness with the canonical 
$\R$-contracting structures!
On the other hand, if $h:E_a\lra E'_{h(a)}$ is
$\R$-homogeneous 
for the standard $\R$-contracting structures on $E_a$ and $E'_{h(a)}$
defined in 2.1.8, it is also $\R_+$-homogeneous (see 2.1.3).

8) If $f:E\lra E'$ is an affine map (between n.v.s.), then, for all $a\in E$,
the map $f:E_a\lra E'_{f(a)}$ is $\R$-homogeneous
(for the standard $\R$-contracting structures on $E_a$ and $E'_{f(a)}$
defined in 2.1.8).
\end{prop}

\proof
1) $h(\omega)=h(0\star\omega)=0\star\omega=\omega$.

2) Using the injectivity of $t\star(-)$, we must verify that
$t\star(t\buildrel{-1}\over\star h(x))=t\star h(t\buildrel{-1}\over\star x)$
which is true thanks to the $\Sigma$-homogeneity of $h$.

3) Refer to 2.1.8.

4) Obvious.

5) According to 2.1.9, we know that
$M$ and $M'$ are also
$\Sigma$-contracting spaces (to make clearer this proof, we denote
here $\star'$ the $\Sigma'$-operations). The map $h$ being
$\Sigma'$-homogeneous, we have
$h(t'\star' x)=t'\star' h(x)$ for all $x\in M$ and $t'\in \Sigma'$
(thus for all  $t'\in\sigma(\Sigma)$); the fact that $h$ is\break
$\Sigma$-homogeneous results immediately from
the formula $\sigma(t)\star' x=t\star x$ of 2.1.9.
Now, if $h$ is $\Sigma$-homogeneous and $\Sigma'=\sigma(\Sigma)$,
we write, for $x\in M$ and $t'=\sigma(t)\in\Sigma'$,
$h(t'\star' x)=h(\sigma(t)\star' x)=h(t\star x)=t\star h(x)=
\sigma(t)\star' h(x)=t'\star' h(x)$ to obtain the 
$\Sigma'$-homogeneity
of $h$.

6) and 7) are immediate consequences of 5).

8) For every $x\in E$, we can write $f(x)=f(a)+l(x-a)$, where 
$l:E\lra E'$ is linear. Thus, for every $t\in\R$, we have 
$f(t\star x)=\break
f(a+t(x-a))=f(a)+tl(x-a))=
f(a)+t(f(x)-f(a))=t\star f(x)$.
\cqfd

\begin{prop}
For every $\Sigma$-homogeneous map $h:M\lra M'$ (where $M$ and $M'$ are $\Sigma$-contracting spaces), we have the equivalence:
\centerline{$\ h\ $ $k$-lipschitzian $\ \Llra\ $ $h\ $ $k$-$LL_\omega$}
\end{prop}

\proof
Referring to 1.1.8, we have just to prove the implication
$\Lla$. So,
let us assume that $h$ is $k$-$LL_\omega$; thus,  there exists $r>0$ such that
$h|_{B(\omega,r)}$ is $k$-lipschitzian.
Let us consider $x,y\in M$ and $t\in\Sigma$ verifying $0<v(t)<1$ 
(thus, this $t$ verifies $\lim_n v(t^n)=\lim_n(v(t))^n=0$).

- If $x\not=\omega$ and $y\not=\omega$, there exists $n\in\N$ such that\break
$v(t^n)<\inf(\frac{r}{d(x,\omega)},\frac{r}{d(y,\omega)})$.
Let us put $x'=t^n\star x$ and $y'=t^n\star y$; they verify
$x',y'\in B(\omega,r)$ (since $d(x',\omega)=d(t^n\star x,t^n\star\omega)
=v(t^n)d(x,\omega)<r$; same for $y'$). So that $v(t^n)d(h(x),h(y))=
d(h(x'),h(y'))\leq kd(x',y')=kv(t^n)d(x,y)$, which implies
$d(h(x),h(y))\leq kd(x,y)$ (for $v(t^n)\not=0$).

- If $x=\omega$ or $y=\omega$ (with $x\not=y$), we argue in the same way
(choosing $n\in\N$ such that $v(t^n)<\frac{r}{d(x,\omega)}$ or
$v(t^n)<\frac{r}{d(y,\omega)}$). The case $x=y=\omega$
is obvious.
\cqfd


\begin{theo}\textup{(of $\Sigma$-uniqueness)} 

\noindent Let $h_1,h_2:M\lra M'$ two
$\Sigma$-homogeneous maps (where $M,M'$ are\break
$\Sigma$-contracting spaces) verifying $h_1\tang_\omega h_2$;
then $h_1=h_2$.
\end{theo} 

\proof
Let us take $t\in\Sigma$ verifying $0<v(t)<1$ and fix $x\in M$.
We can assume that $x\not=\omega$ (since, as seen in 2.1.17, $h_1(\omega)=
h_2(\omega)=\omega$).
Let us set $x_n=t^n\star x$ for each $n\in\N$.
Then, as seen in 2.1.10, we have $x_n\not=\omega$ and
$\lim_n x_n=\omega$, so that we can write:
$\frac{d(h_1(x_n),h_2(x_n))}{d(x_n,\omega)}=
\frac{d(h_1(t^n\star x),h_2(t^n\star x))}{d(t^n\star x,t^n\star\omega)}=
\frac{d(t^n\star h_1(x),t^n\star h_2(x))}{d(t^n\star x,t^n\star\omega)}=
\frac{v(t^n)d(h_1(x),h_2(x))}{v(t^n)d(x,\omega)}=
\frac{d(h_1(x),h_2(x))}{d(x,\omega)}$.
But, as $\lim_n x_n=\omega$ and $h_1\tang_\omega h_2$, 
this provides
$0=\lim_n\frac{d(h_1(x_n),h_2(x_n))}{d(x_n,\omega)}\break
=\frac{d(h_1(x),h_2(x))}{d(x,\omega)}$ which implies 
$d(h_1(x),h_2(x))=0$, i.e $h_1(x)=h_2(x)$, which ends the proof.
\cqfd

\begin{theo}
Let $M,M'$ be $\Sigma$-contracting spaces with $M'$ reverti-\break
ble, 
$V$ a neighborhood of $\omega$ in $M$, and maps $f:V\lra M'$,
$h:M\lra M'$ such that $h$ is $\Sigma$-homogeneous 
and verifies $f\tang_\omega h|_V$. Then, for all $x\in M$, we have
$h(x)=\lim_{0\not=v(t)\rightarrow 0} t\buildrel{-1}\over\star f(t\star x)$.
\end{theo}

\proof
The above equality is clearly true for $x=\omega$
(using 2.1.15 and the fact that $f\tang_\omega h|_V$
implies $f(\omega)=h(\omega)=\omega$ (see 2.1.17)).
Now, if $x\not=\omega$, we have 
$\omega=\lim_{v(t)\rightarrow 0}t\star x$ (since
$d(t\star x,\omega)=
d(t\star x,t\star\omega)\break
=v(t)d(x,\omega)$),
which insures that 
for all $t\in\Sigma$ such that $0<v(t)<\varepsilon$,\break
 we have 
$t\star x\in V-\{\varepsilon\}$. Thus, for all these $t$,
we can write:
$\frac{d(f(t\star x),h(t\star x))}{d(t\star x,\omega)}=
\frac{d(f(t\star x),t\star h(x))}{d(t\star x,t\star\omega)}=
\frac{d(t\star(t\buildrel{-1}\over\star f(t\star x)) ,t\star h(x))}
{v(t)d(x,\omega)}= 
\frac{v(t)d(t\buildrel{-1}\over\star f(t\star x),h(x))}
{v(t)d(x,\omega)}=
\frac{d(t\buildrel{-1}\over\star f(t\star x),h(x))}{d(x,\omega)}$.
Since 
$f\tang_\omega h|_V$, we have 
$\lim_{0\not= v(t)\rightarrow 0}\frac{d(f(t\star x),h(t\star x))}{d(t\star x,\omega)}=0$, so that\break
$\lim_{0\not= v(t)\rightarrow 0}\frac{d(t\buildrel{-1}\over\star
f(t\star x),h(x))}{d(x,\omega)}=0$, i.e.
$\lim_{0\not= v(t)\rightarrow 0}d(t\buildrel{-1}\over\star 
f(t\star x),h(x))=0$,
which finally means that     $\lim_{0\not= v(t)\rightarrow 0}
d(t\buildrel{-1}\over\star f(t\star x))=h(x)$.
\cqfd

\begin{theo}
Let $M,M'$ be $\Sigma$-contracting spaces with $M'$ revertible,
$V$ a neighborhood of $\omega$ in $M$, 
$g:V\lra M'$ a $k$-$LL_\omega$ map, and\break
$h:M\lra M'$ 
a $\Sigma$-homogeneous map verifying
$g\tang_\omega h|_V$. Then,\break
$h$ is $k$-lipschitzian.
\end{theo}

\proof
Let $W$ be a neighborhood of $\omega$ in $V$ (thus in $M$) such that
$g|_W$ is $k$-lipschitzian. Let $x,y\in M$ and $t\in\Sigma$ such that\break
$0<v(t)<1$; then, there exists $N\in\N$ such that 
$t^n\star x,t^n\star y\in W$ for all\break
$n\geq N$ (since $\lim_n t^n\star x=
\lim_n t^n\star y=\omega$); so that, for all these $n$, we have
$d(g(t^n\star x),g(t^n\star y))\leq kd(t^n\star x,t^n\star y)=kv(t^n)d(x,y)$,
which provides, thanks to 2.1.15, 
$d(t^n\buildrel{-1}\over\star g(t^n\star x),
t^n\buildrel{-1}\over\star g(t^n\star y))=
(v(t^n))^{-1}d(g(t^n\star x),g(t^n\star y))
\leq kd(x,y)$. Now, $d$ being continuous, we obtain (doing 
$n\rightarrow +\infty$, so that $v(t^n)=(v(t))^n\rightarrow 0$):  
$d(h(x),h(y))\leq kd(x,y)$, thanks to 2.1.20.
\cqfd

\begin{cory}
$M,M'$ being as in 2.1.21, and $h:M\lra M'$ being a
$\Sigma$-homogeneous map, we have the equivalences:\par
\qquad\qquad $h\, $ lipschitzian $\ \Llra\ $ $h\, $ $LL_\omega$ 
$\ \Llra\ $ 
$h\, $ $Tang_\omega$
\end{cory}

\proof
Referring to 1.3.2 and 2.1.18, we
have only to prove that: $h\,$ $Tang_\omega$ $\ \Lra\ $
$h\,$ lipschitzian; since $h\,$ $Tang_\omega$ means that there exists
$g:M\lra M'$ which is $LL_\omega$ and verifies
$g\tang_\omega h$, we use 2.1.21.
\cqfd

\vspace{3mm}
\begin{cexam}
{}
\end{cexam}
\includegraphics {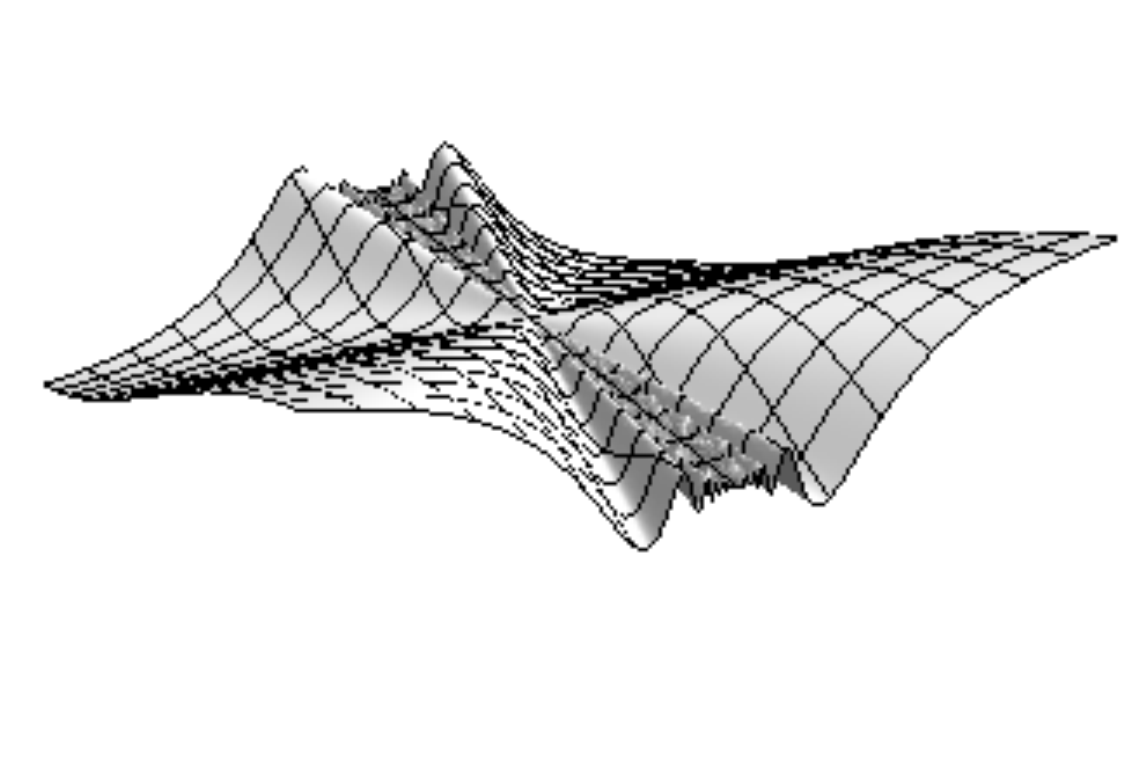}
\begin{picture}(0,0)
\put(-68,10){Figure 1}
\end{picture}

We give here an example of function $\R^2\lra\R$ which is 
$\R$-homogene-\break
ous but not lipschitzian (thus not $Tang_0$).
Consider the function\break
$f(x,y)=x\sin\frac{y}{x}$ if $x\not=0$ and $f(0,y)=0$ (Fig 1).
This $f$ is\break
 $\R$-homogeneous
$\R_0^2\lra\R_0$ since, for $t\in\R$ and $(x,y)\in\R^2$, we have:

- when $t\not=0$ and $x\not=0$, $f(tx,ty)=tx\sin\frac{y}{x}=tf(x,y)$,

- when $t=0$ or $x=0$, $f(tx,ty)=f(0,ty)=0=tf(0,y)=tf(x,y)$.

We also notice that $f$ is $LSL$ (at every point) since it is differentiable (then tangentiable and thus LSL: see 1.3.2 and 1.3.4)
on $\R^*\times\R$; and, for all $a\in\R$ and  $(x,y)\in \R^2$, we have
$|f(x,y)-f(0,a)|=|f(x,y)|\leq|x|\leq\|(x,y)-(0,a)\|_2$ (the euclidian norm).

Now, if $x\not=0$, we have $\frac{\partial f}{\partial x}(x,y)=
\sin\frac{y}{x}-\frac{y}{x}\cos\frac{y}{x}$ and thus, putting
$x_n=\frac{1}{n^2}$ and $y_n=\frac{2\pi}{n}$ (so that 
$\frac{y_n}{x_n}=2\pi n$), we obtain 
$\frac{\partial f}{\partial x}(x_n,y_n)=-2\pi n$, and thus
$\lim_n\frac{\partial f}{\partial x}(x_n,y_n)=-\infty$, where
$\lim_n(x_n,y_n)=(0,0)$.\break
Thus, this function $f$ cannot be lipschitzian!

\vspace{3mm}
\begin{remk}
By the way, we have proved that :
$LSL_a\not\!\!\Lra Tang_a$.\break
Thus, we cannot complete our equivalences of 
2.1.22, adding
the properties of being $LSL_\omega$ or continuous at $\omega$!
Even though, for linear maps, all these properties are equivalent
(using 1.3.2 and the fact that, for a linear map, 
continuous $\,=\,$ lipschitzian).
\end{remk}

\section{Representability and 
Contactibility}
Once we have restricted ourselves to the $\Sigma$-homogeneous maps which are lipschitzian,
the $\Sigma$-uniqueness property (see 2.1.19) allows us to choose
at most one canonical representative element in each 
linked metric jet (in short: jet)  between
$\Sigma$-contracting spaces (pointed by their central point $\omega$); hence the term of 
``$\Sigma$-representable'' jets. It is then natural to be interested in
the maps $f$ which are tangentiable at $\omega$, and whose tangent jet at $\omega$ (i.e its tangential T$f_\omega$) is 
$\Sigma$-representable; such maps are called
``$\Sigma$-contactable'' at $\omega$.
That way, we obtain a natural generalization of the notion of differentiability (at a point). 
We then refind for $\Sigma$-contactable maps some well-known properties of the differentiability.

\vspace{3mm}
\begin{remk}
Let $\Sigma$ be a valued monoid, and $M,M'$ two $\Sigma$-contrac-ting spaces. To lighten the text,  a map $h:M\lra M'$ which is\break 
$\Sigma$-homogeneous and lipschitzian will be called
$\Sigma$-Lhomogeneous 
(such a map is, in particular, $Tang_\omega$).
\end{remk}

\begin{remk}\par\hfill

Referring to 2.1.18, a $\Sigma$-Lhomogeneous map is a
$\Sigma$-homogeneous which is $LL_\omega$.
\end{remk}

\begin{prop}
Let $h:M\lra M'$ a $\Sigma$-homogeneous map
(between $\Sigma$-contracting spaces) which is $Tang_\omega$; 
let us also assume that we have at our disposal 
an isometric embedding $j:M'\lra \widehat{M'}$ which is\break 
$\Sigma$-homogeneous, where $\widehat{M'}$ is a revertible
$\Sigma$-contracting space.
Then $h$ is $\Sigma$-Lhomogeneous.
\end{prop}

\proof
By hypothesis, there exists a $LL_\omega$ map $g:M\lra M'$ such that
$g\tang_\omega h$. Since the maps $j.g$ and $j.h$ are respectively
$LL_\omega$ and
$\Sigma$-homogeneous, verify 
$j.g\tang_\omega j.h$  (for $j$ is lipschitzian)
and have a revertible $\Sigma$-contracting codomain, we can use 
2.1.21 to claim that $j.h$ is lipschitzian, which
implies that $h$ is also lipschitzian (since $j$ is isometric).
\cqfd

\vspace {3mm}
Let us now denote $\Sigma$-$\C$ontr, the category whose objects
are the\break 
$\Sigma$-contracting spaces and whose morphisms are the
$\Sigma$-Lhomogeneous maps (it is a suitable world for guarantying the whished $\Sigma$-uniqueness property).
When $\sigma:\Sigma\lra\Sigma'$
is a morphism of valued monoids,
there exists a
canonical functor
$\widehat\sigma:
\Sigma'$-$\C$ontr$\ \lra\Sigma$-$\C$ontr (see 2.1.9 and 
2.1.17). For every valued monoid $\Sigma$, we also have another canonical functor
$J:\Sigma$-$\C$ontr$\ \lra\J$et defined by $JM=(M,\omega)$ and $J(h)=$T$h_\omega$ (since $h$ is $LL_\omega$, actually lipschitzian, we have $h\in$ T$h_\omega$: see 1.3.3).

\begin{prop}
$\ \Sigma$-$\C$ontr is a cartesian category and the functor\break
$J:\Sigma$-$\C$ontr $\lra\J$et is a strict morphism of cartesian categories.
\end{prop}

\proof
Clearly, $\{\omega\}$ is a $\Sigma$-contracting space which is a final
object in $\Sigma$-$\C$ontr. Moreover, if 
$M_1,M_2\in|\Sigma$-$\C$ontr$|$, $M_1\times M_2$ has a canonical
structure of $\Sigma$-contracting space; the canonical
projections\break
$p_i:M_1\times M_2\lra M_i$ (are $\Sigma$-Lhomogeneous and,
if $M$ is another\break
$\Sigma$-contracting space, with
$h_i:M\lra M_i$ $\Sigma$-Lhomogeneous (with\break
$i\in\{1,2\}$), so is the pair
$(h_1,h_2):M\lra M_1\times M_2$ 
(see the beginning of section 1.2 and 2.1.17).

Besides, $J$ is the composite
$\Sigma$-$\C$ontr $\buildrel U\over\lra\L$L
$\buildrel q\over\lra\J$et, where (see 1.2.2 for $q$), the canonical surjection $q$ is, itself  a strict morphism of cartesian categories, just as the forgetful functor $U$ (which only keeps the 
lipschitzian property).
\cqfd
\vspace{3mm}
Now, since for each $M,M'\in|\Sigma$-$\C$ontr$|$,
the canonical map:\break
$\Sigma$-$\C$ontr$(M,M')\buildrel can\over\lra\J$et
$\!\!(JM,JM')$, defined by $can(h)=J(h)=$T$h_\omega$, 
is injective (thanks to the $\Sigma$-uniqueness property: see 2.1.19),
we can equip $\Sigma$-$\C$ontr$(M,M')$ with
the distance $d(h,h')=d(J(h),J(h'))$ (for the distance between jets
(see 1.2.6 and 1.2.7);
This fact provides\break
$\Sigma$-$\C$ontr with a structure of category enriched in $\M$et
(see 1.2.20).
More precisely, this will arise from 2.2.6.

\begin{remk} 
Apparently, we have two definitions for $d(h,h')$: the previous one and the 
one of 1.2.6 (where we can set $V=M$ since $h$ and $h'$ are lipschitzian here).
Actually, these two definitions coincide here (since 
$h\in$ \textup{T}$h_\omega$ and  
$h'\in $\textup{T}$h'_\omega$).
As, by definition, the new one is a distance, so is the old one 
(which, a priori, was not for maps solely $LL_\omega$); this is not surprising
if we think of the $\Sigma$-uniqueness theorem 2.1.19 for 
$\Sigma$-homogeneous maps.
\end{remk}

\begin{prop}
For each $M,M',M''\in|\Sigma$-$\C\textup{ontr}|$, the composition:
$\Sigma$-$\C\textup{ontr}(M,M')\times\Sigma$-$\C\textup{ontr}(M',M'')
\buildrel{comp}\over\lra\Sigma$-$\C\textup{ontr}(M,M'')$
is $LSL$.
\end{prop}

\proof
For (since $J$ is functorial) the following diagram commutes:
$$\xymatrix{
\Sigma\!-\!\C\textup{ontr}(M,M')\times\Sigma\!-\!\C\textup{ontr}(M',M'')
\ar[d]_{can\times can}
\ \ar[rr]^>>>>>>>>>>{comp}&&\ 
\Sigma\!-\!\C\textup{ontr}(M,M'')
\ar[d]^{can}\\
\J\textup{et}(JM,JM')\times \J\textup{et}(JM',JM'')\ 
\ar[rr]_>>>>>>>>>>>>>>
{comp}&&\ \J\textup{et}(JM,JM'')
}$$
Where, by definition of the distance on $\Sigma$-$\C$ontr$(M,M')$, the vertical maps are isometric embeddings, and the down horizontal
$comp$ is $LSL$ (since $\J$et is enriched in $\M$et: 
see 1.2.20).
\cqfd

Actually, the cartesian category $\Sigma$-$\C$ontr (see 2.2.4) is enriched in $\M$et, as shown in the following proposition:

\begin{prop}
For each $M,M_0,M_1\in|\Sigma\!-\!\C\textup{ontr}|$, the following\break canonical map is an isometry:\par
\noindent $\ \Sigma\!-\!\C\textup{ontr}(M,M_0\times M_1)
\buildrel{can}\over\lra
\Sigma\!-\!\C\textup{ontr}(M,M_0)\times
\Sigma\!-\!\C\textup{ontr}(M,M_1)$.
\end{prop}

\proof
It comes from the commutativity of the following diagram:
$$\xymatrix{
\Sigma\!\!-\!\!\C\textup{ontr}(M,M_0\!\times\! M_1)
\ar[d]_{can}
\ar[rr]^>>>>>>>>>>{can}&&
\Sigma\!\!-\!\!\C\textup{ontr}(M,M_0)\!\times\!
\Sigma\!\!-\!\!\C\textup{ontr}(M,M_1)
\ar[d]^{can\times can}\\
\J\textup{et}(JM,JM_0\times JM_1)
\ar[rr]_>>>>>>>>>>>>>>
{can}&&\J\textup{et}(JM,JM_0)\!\times\!\J\textup{et}(JM,JM_1)
}$$
The fact that the down horizontal $can$ is an isometry has been proved in 1.2.22 (we have used the fact that the functor 
$J$ commutes to the products).
\cqfd

\begin{prop}
Let $h,h'\in\Sigma\!-\!\C\textup{ontr}(M,M')$;  
then, $d(h,h')$ being defined as before 2.2.5,  we have
$d(h,h')=\sup\{C(x)\,|\,x\in B'(\omega,1)\}$, where 
$C(x)=\frac{d(h(x),h'(x))}{d(x,\omega)}$ if $x\not=\omega$ and $C(\omega)=0$.
\end{prop}

\proof
Referring to before 1.2.6, and 2.2.5, we must show that $d(h,h')=d^1(h,h')$;
which will arise from the fact that, for all $r\leq 1$, we have $d^r(h,h')=d^1(h,h')$.
First, $d^r(h,h')\leq d^1(h,h')$ comes from the inclusion
$B'(\omega,r)\subset B'(\omega,1)$.
Besides, if $x\in B'(\omega,1)$, we choose $t\in\Sigma$
and $n\in\N$ such that $0<v(t)<1$ and $t^n\star x\in B'(\omega,r)$
(which is possible, since $\lim_n t^n\star x=\omega$: see 2.1.10). Now, 
referring to the proof of  
2.1.19, we have 
$C(x)=C(t^n\star x)\leq d^r(h,h')$, which gives $d^1(h,h')\leq d^r(h,h')$.
It remains to do $r\rightarrow 0$.
\cqfd

\begin{prop}
Let $M,M'$ be $\Sigma$-contracting spaces where $M'$ is\break
revertible, and $h:M\lra M'$ a $\Sigma$-Lhomogeneous map.
Let us set\break
$k=\sup\{\frac{d(h(x),h(y))}{d(x,y)}\,|\, x,y\in B'(\omega,1);
x\not= y\}$. Then:
 
1) $h$ is $k$-lipschitzian,

2) $\rho(\textup{T}h_\omega)=k$ (see 1.2.8 for the
lipschitzian ratio $\rho$ of a jet).
\end{prop}

\proof
1) Since, by definition of $k$, $h|_{B'(\omega,1)}$ is 
$k$-lipschitzian, we know (see 2.1.18 and 2.2.2) that $h$ is actually $k$-lipschitzian.

2) Since $h\in $T$h_\omega$, we obtain immediately
$\rho($T$h_\omega)\leq k$. In order to obtain the inverse inequality, 
let us take $f\in $T$h_\omega$ supposed to be $R$-$LL_\omega$
(with $R>0$); then, by 2.1.21, we know that 
$h$ is $R$-lipschitzian; thus, for all $x,y\in B'(\omega,1)$ 
($x\not=y$), we have $\frac{d(h(x),h(y))}{d(x,y)}\leq R$, so that
$k\leq R$, still by definition of $k$; 
thus $k\leq\rho($T$h_\omega)$.
\cqfd

\begin{remk}
The previous proposition  gives us a calculus of 
the lipschitzian ratio $\rho(\varphi)$ for a $\Sigma$-representable jet
$\varphi:(M,\omega)\lra(M',\omega)$ (see the definition of the
$\Sigma$-representability of a jet just below).
\end{remk}

\vspace{3mm}
\begin{defi}
Consider two $\Sigma$-contracting spaces $M,M'$ and a jet\break
$\varphi:(M,\omega)\lra(M',\omega)$.
We say that:

1) $\varphi$ is $\Sigma$-\textup{representable} if there exists a  
$\Sigma$-Lhomogeneous element\break
$h:M\lra M'$ verifying $J(h)=\varphi$ 
(i.e. \textup{T}$h_\omega=\varphi$).
Thanks to 2.1.19, we know  that such a $h$ is then unique.

2) The $h$ mentioned just above (for a $\Sigma$-representable
jet $\varphi$) is the $\Sigma$-\textup{representative element}
of the jet $\varphi$.
\end{defi}

\begin{remk}
We will see, in the final summary of section 2.5 (more precisely in  2.5.11), that the jets between
$\Sigma$-contracting spaces are not always $\Sigma$-representable, even if 
$M=M'=\R_0$ (for $\Sigma=\R_+$ or $\N'_r$ ...).
\end{remk}

\begin{prop}
The $\Sigma$-representable jets are stable under composition, pairs (and thus products).
\end{prop}

\proof
Since $J$ is a functor which commutes with products (see 2.2.4).
\cqfd

\vspace{3mm}
For convenience, we now call $\Sigma$-\textit{contracting domain}
a pair $(M,V)$ where $M$ is a $\Sigma$-contracting space and $V$ a neighborhood of $\omega$ in $M$.
Now, if $(M,V)$ et $(M',V')$ are
$\Sigma$-contracting domains,
a \textit{centred map}\break 
$f:(M,V)\lra(M',V')$ is merely a map
$f:V\lra V'$ verifying\break
$f(\omega)=\omega$.

\begin{defi}
Let $(M,V)$ et $(M',V')$ be
$\Sigma$-contracting domains and
$f:(M,V)\lra(M',V')$ a centred map.
We say that $f$ is\break 
$\Sigma$-\textup{contactable} 
(at the central point $\omega$)
if $f:V\lra V'$ is tangentiable at $\omega$ and if the following
composite jet is $\Sigma$-representable:
$$\xymatrix{
(M,\omega)\ar[rr]^{\sim}&&(V,\omega)
\ar[rr]^{\textup{T}f_\omega}&&(V',\omega)
\ar[rr]^{\sim}&&(M',\omega)
}$$
Where the $\sim$ are nothing but the invertible jets $j_\omega^{-1}$ and $j'_\omega$
(see 1.2.5).
We denote $\textup{K}_\Sigma f:M\lra M'$ (or merely $\textup{K}f$
if there is no ambiguity about the valued monoid $\Sigma$)
the unique representative element of the above composite jet; 
and we call it the $\Sigma$-\textup{contact} of $f$.
\end{defi}

\begin{remk}
In other words, $f$ is $\Sigma$-contactable if
there exists a
$\Sigma$-Lhomogeneous
$h:M\lra M'$ 
such that $f\tang_\omega h|_V$ (where, here, $f$ is seen as
a map $V\lra M'$). In that case, $h=\textup{K}f$. 
If $f$ is itself $\Sigma$-Lhomogeneous,
it is $\Sigma$-contactable with $f=\textup{K}f$!
\end{remk}

\begin{prop}
Let $\sigma:\Sigma\lra\Sigma'$ be a morphism of valued monoids.
If $f:(M,V)\lra(M',V')$ is a $\Sigma'$-contactable centred map, it is also $\Sigma$-contactable, with
$\textup{K}_{\Sigma}f=\textup{K}_{\Sigma'} f$ ($M$ and $M'$ being   equipped with the $\Sigma$-contracting structure induced by $\Sigma'$).
\end{prop}

\proof
It arises from 5) in 2.1.17 and the $\Sigma$-uniqueness property (see 
2.1.19).
\cqfd

\begin{prop}
Let $E,E'$ be n.v.s., $U$ an open subset of $E$,\break
$a\in U$ and
$f:U\lra E'$ a map. If $f$ is differentiable at $a$, then
$f:(E_a,U)\lra(E'_{f(a)},E')$ is $\R$-contactable (for the standard 
$\R$-contracting structures)  with
$\textup{K}_\R f(x)=f(a)+\textup{d}f_a(x-a))$;
and, for every valued monoid $\Sigma$,
it is even $\Sigma$-contactable (for the canonical 
$\Sigma$-contracting structures: see 2.1.10) with
$\textup{K}_\Sigma f=\textup{K}_\R f$ as written above.
\end{prop}

\proof
Referring first to 2.1.8 and  2.1.17, we notice that the 
(standard) $\R$-contracting structures on $E_a$ and $E'_{f(a)}$
makes $\R$-Lhomogeneous
$E_a\lra E'_{f(a)}$
the affine map $h:x\mapsto f(a)+$d$f_a(x-a)$.
From the fact that $f\tang_a h|_U$,
we deduce that  $f$ is $\R$-contactable, with K$_{\R}f=h$;
this\break
$\R$-contact $h$ being $\R$-Lhomogeneous,
it is $\R_+$-Lhomogeneous (see 7) in 2.1.17), so that $f$ is 
$\R_+$-contactable, with 
K$_{\R_+}f=h$. We can now use
2. 2.16, to obtain that $f$ is also $\Sigma$-contactable 
($E_a$ and $E'_{f(a)}$ being now seen as 
canonical $\Sigma$-contracting spaces) and that K$_\Sigma f=$K$_{\R_+}f=h$.
\cqfd

\begin{remks}\par\hfill

1) The result of 2.2.17 is remarkable:
the continuous affine map which is tangent at $a$ to a map, 
differentiable at $a$, is its canonical\break
$\Sigma$-contact at $a$ for any valued monoid $\Sigma$!
 
2) We will see, in sections 2.4 and 2.5 (more precisely in 2.5.11), examples of functions $\R\lra\R$
which are $\R_+$-contactable (at 0) (which means that
$f:(\R_0,\R)\lra(\R_{f(0)},\R)$ is $\R_+$-contactable) but not 
differentiable at 0; and examples of functions $\N'_r$-contactable (at 0) but 
not $\R_+$-contactable (at 0).
\end{remks}

\begin{prop}
$\!\!$Let $(M,V),(M',V')$ be $\Sigma$-contracting domains
with $M'$ revertible, and $f:(M,V)\lra(M',V')$ a centred 
$\Sigma$-contactable map.
Then, for all $x\in M$, we have
$\textup{K}f(x)=\lim_{0\not= v(t)\rightarrow 0}\,
t\,\buildrel{-1}\over\star f(t\star x)$.
\end{prop}

\proof
We just have to set $h=$K$f$ in 2.1.20.
\cqfd

\vspace{3mm}
For the following propositions, we use
1.3.6, 1.3.7, 1.3.8 and 2.2.13.

\begin{prop}
Let $(M,V),(M',V'),(M'',V'')$ be $\Sigma$-contracting\break
domains and
$f:(M,V)\lra(M',V')$, $g:(M',V')\lra(M'',V'')$ two centred maps.
If $f$ and $g$ are $\Sigma$-contactable, so is $g.f$, and we have
$\textup{K}(g.f)=\textup{K}g.\textup{K}f$.
\end{prop}

\proof
Because T$(g.f)_\omega=$T$g_\omega.$T$f_\omega$.
\cqfd

\begin{prop}
Let $(M,V),(M_0,V_0),(M_1,V_1)$ be $\Sigma$-contracting\break
domains and
$f_0:(M,V)\lra(M_0,V_0)$,
$f_1:(M,V)\lra(M_1,V_1)$ two centred maps.
If $f_0$ and $f_1$ are $\Sigma$-contactable, so is the pair
$(f_0,f_1):\break
(M,V)\lra(M_0\times M_1,V_0\times V_1)$,
and we have
$\textup{K}(f_0,f_1)=(\textup{K}f_0,\textup{K}f_1)$.
\end{prop}

\proof
Because T$(f_0,f_1)_\omega=($T$f_{0\omega},$T$f_{1\omega})$
\cqfd

\begin{prop}
Let $(M_0,V_0),\!(M'_0,V'_0),\!
(M_1,V_1),\!(M'_1,V'_1)$ be
$\Sigma$-contrac-
ting domains, and
$f_0:(M_0,V_0)\lra(M'_0,V'_0)$,
$f_1:(M_1,V_1)\lra(M'_1,V'_1)$
two centred maps. 
If $f_0$ and $f_1$ are $\Sigma$-contactable, so is the product
$f_0\times f_1:(M_0\times M_1,V_0\times V_1)\lra
(M'_0\times M'_1,V'_0\times V'_1)$, and we have
$\textup{K}(f_0\times f_1)=\textup{K}f_0\times \textup{K}f_1$.
\end{prop}

\proof
Because T$(f_0\times f_1)_{(\omega,\omega)}=$T$f_{0\omega}\times$ T$f_{1\omega}$
\cqfd

%% file: jet_II2.tex
\section{Contactibility in the n.v.s.}
Since most of our examples of contactibility take place inside the n.v.s. world, we immerse ourselves, from now on, in this context.
It gives us the opportunity to lighten some previous notations 
and to make precise some others that
we will use for the two main great classes of examples studied all along the sections 2.4 and 2.5.

\vspace{1mm}
Let us briefly recall that, for each valued monoid $\Sigma$,
every n.v.s. $E$ can be considered as a
canonical $\Sigma$-contracting
(metric) space $E_a$ (where $a\in E$ is a chosen central point), the canonical external operation being $t\star x=a+v(t)(x-a)$ for all $t\in\Sigma$ and 
$x\in E$ (refer to 2.1.10).
Actually, referring to 2.1.14, we know that 
$E_a$ is a revertible (canonical)
$\Sigma$-contracting space (with
$t\buildrel{-1}\over\star x=a+(v(t))^{-1}(x-a)$ for $t\not=0$.

\vspace{1mm}
\texttt{Untill the end of this section 2.3, we only work with the canonical structures}.  

\begin{prop}
Let $E,E'$ be two n.v.s., and $a\in E$, $a'\in E'$.
We have an isometry:
$\Theta_{aa'}:\Sigma$-$\C\textup{ontr}(E_0,E'_0)\lra
\Sigma$-$\C\textup{ontr}(E_a,E'_{a'})$
defined by $\Theta_{aa'}(h)(x)=a'+h(x-a)$; by definition, $\Theta_{aa'}(h)$ is a translate of $h$.
\end{prop}

\proof
One easily verifies that $\Theta_{aa'}(h)$ is $\Sigma$-Lhomogeneous. The fact that
$\Theta_{aa'}$ is an isometry comes from the fact that (see 2.2.8), for $x\not=a$, we have: $\frac{d(\Theta_{aa'}(h_1)(a+x),\Theta_{aa'}(h_2)(a+x))}{d(a+x,a)}=
\frac{d(h_1(x),h_2(x))}{d(x,0)}$ (we use the fact that 
$x\in B'(0,1)\Llra x+a\in B'(a,1)$).
\cqfd

\vspace{2mm}
From now on, when $E$ and $E'$ are two n.v.s., the metric space\break
$\Sigma$-$\C$ontr$(E_0,E'_0)$ will be denoted
$\Sigma$-$\L$Hom$(E,E')$; its elements are thus the
$\Sigma$-Lhomogeneous maps $h:E_0\lra E'_0$, i.e the lipschitzian maps $h:E\lra E'$ which verify, for all $t\in\Sigma$ and $x\in E$:
$h(v(t)x)=v(t)h(x)$.

\begin{prop}
Let $E,E'$ be two n.v.s., then the distance on\break
$\Sigma$-$\L\textup{Hom}(E,E')$, defined around remark 2.2.5, derives from the norm\break
$\|h\|=\sup\{C(x)\,|\,\|x\|\leq 1\}$, where $C(x)=\frac{\|h(x)\|}{\|x\|}$ if $x\not=0$ and $C(0)=0$,
\end{prop}

\proof
First, $\Sigma$-$\L\textup{Hom}(E,E')$ is a vector space 
(it is a subvector space of $\L$L$((E,0),(E',0))$ (see 1.2.25)).
Now, we use the injective map 
$can:\Sigma$-$\L\textup{Hom}(E,E')\lra
\J$et$((E,0),(E',0))$ described just before 2.2.5; this $can$ is also linear
(being the following composite (see 1.2.25 for $q$):
$\Sigma$-$\L\textup{Hom}(E,E')
\buildrel j\over\lra\L$L$((E,0),(E',0))\buildrel q\over\lra
\J$et$((E,0),(E',0))$, where  the canonical surjection $q$ is, itself linear, just as the canonical injection $j$.
Thus, we can define a norm on $\Sigma$-$\L\textup{Hom}(E,E')$
setting (see 1.2.25) $\|h\|\!=\!\|can(h)\|\!=\!\|q(h)\|\!=\!d(q(h),q(0))\!=\!d(h,0)\!=\!
\sup\{C(x)\,|\,\|x\|\leq 1\}$, where $C(x)=\frac{\|h(x)\|}{\|x\|}$ if $x\not=0$ and $C(0)=0$
(using 2.2.8).
\cqfd

\begin{remks}\par\hfill

1) We notice that the previous norm on $\Sigma$-$\L\textup{Hom}(E,E')$ is an extension
of the operator norm on its subvector space $L(E,E')$.

2) Referring to 2.2.9 and 2.3.2, we notice that, for all\break
$h\in\Sigma$-$\L\textup{Hom}(E,E')$, we have
$\|h\|=\|\textup{T}h_0\|=d(\textup{T}h_0,O)\leq\rho(\textup{T}h_0)$, the last inequality being not surprising (it's a general result: see 1.2.10). 
\end{remks}

\begin{prop}
As for continuous linear maps, we have:

1) For $h\in\Sigma$-$\L\textup{Hom}(E,E')$ and $x\in E$:
$\ \|h(x)\|\leq\|h\|\,\|x\|$,

2) For $h\!\!\in\!\!\Sigma$-$\L\textup{Hom}(E,E')$ and
$k\!\!\in\!\!\Sigma$-$\L\textup{Hom}(E',E'')$:
$\ \|k.h\|\leq\|k\|\, \|h\|$.
\end{prop}

\proof
Let $t_0\in\Sigma$ verifying $0<v(t_0)<1$.

1) For $x\in E$, $x\not=0$, let $n\in\mathbb Z$ such that 
$\|x\|\leq v(t_0)^n$, thus such that $\|v(t_0)^{-n}x\|\leq 1$; 
so, we have $\frac{\|h(x)\|}{\|x\|}=
\frac{\|h(v(t_0)^{-n}x)\|}{\|v(t_0)^{-n}x\|}\leq \|h\|$.

2) Proceed as in the linear case.
\cqfd

\begin{remk}
We can also use the inequality (see 1.2.17): \break
$d(\textup{T}(k.h)_0,O)\leq d(\textup{T}k_0,O)
d(\textup{T}h_0,O)$, with $O=O_{00}$
(see 1.2.10 for the notation of 2) in 2.3.4.
\end{remk}

\vspace{3mm}
Here again, for convenience, we call \textit{normed domain} a pair $(E,U)$ where $E$ is a n.v.s. and $U$ an open subset of $E$. When $(E,U)$ and $(E',U')$ are normed
domains, a map $f:(E,U)\lra(E',U')$ is merely a map $f:U\lra U'$.
When, $E=U$ and $E'=U'$, we merely write $f:E\lra E'$.

\begin{defi}
Let $(E,U)$and $(E',U')$ two normed domains and\break
$f:(E,U)\lra(E',U')$ a map. We say that $f$ is $\Sigma$-\textup{contactable} at $a\in U$ if $f:(E_a,U)\lra(E'_{a'},U')$ is $\Sigma$-contactable, where $a'=f(a)$ (see 2.2.14 and 2.2.15). Then
we set $\textup{k}_\Sigma f_a=\Theta_{aa'}^{-1}(\textup{K}f)$
(or briefly $kf_a$).
We say that
$f$ is $\Sigma$-contactable on $U$ if it is $\Sigma$-contactable at every point of $U$.
\end{defi}

\begin{remks}\par\hfill

1) In other words, $f$ is $\Sigma$-contactable at $a$ iff
there exists\break
$h\in\Sigma$-$\L$\textup{Hom}$(E,E')$ such that the restriction, to $U$,
of its translate map $x\mapsto f(a)+h(x-a)$ (which is in fact 
$\textup{K}f$)
is tangent to $f$ at $a$.
Such a map $h$ is then unique and verifies $h=\textup{k}f_a$.

2) Actually $\textup{k}f_a(x)=
\textup{K}f(x+a)-a'$ for all $x\in E$, the $\Sigma$-contact of $f$ at $a$ being thus 
$\textup{K}f(x)=f(a)+\textup{k}f_a(x-a)$, a translate at $a$ of the 
$\Sigma$-Lhomogeneous map $\textup{k}f_a$.

3) Every $h\in\Sigma$-$\L$\textup{Hom}$(E,E')$ is $\Sigma$-contactable at 0 with
$\textup{k}h_0=h$.

4) Referring to 2.2.14 and 2.2.15, we know that every $\Sigma$-contactable
map at $a$ is tangentiable at $a$.

5) Let $\sigma:\Sigma\lra\Sigma'$ be a morphism of valued monoids.
If $f:(E,U)\lra(E',U')$ is $\Sigma'$-contactable at $a\in U$, 
thanks to 2.2.16 and 2.3.6, we know that $f$ is also
$\Sigma$-contactable at $a$ and that
$\textup{k}_\Sigma f_a=
\textup{k}_{\Sigma'}f_a$.
\end{remks}

\begin{prop}
$f\ $ differentiable at $a$ $\ \Lra\ $ $f\ $ $\Sigma$-contactable at $a$ for all  $\Sigma$, with $\textup{k}f_a=\textup{d}f_a$.
More precisely, if $f$ is $\Sigma$-contactable at $a$, then $f$ is 
differentiable at $a$ iff $\ \textup{k}f_a$ is linear. 
\end{prop}

\proof
Arises from 1) in 2.3.7, since $L(E,E')$ is a subspace of\break
$\Sigma$-$\L$Hom$(E,E')$.
The reader can refer to the remark 1) in 2.2.18!
\cqfd

\vspace{3mm}
\begin{exams}
{}
\end{exams}
The ``impatient'' reader could go straight to 2.4.7 and 2.5.8
to find such good examples.

\vspace{3mm}
\begin{prop}
Let $(E,U)$, $(E',U')$, $(E'',U'')$ be three normed domains,
$f:(E,U)\lra(E',U')$, $g:(E',U')\lra(E'',U'')$ two maps, $a\in U$ and 
$a'=f(a)$.
If $f$ is $\Sigma$-contactable at $a$ and $g$ $\Sigma$-contactable at
$a'$,\break 
then $g.f$ is $\Sigma$-contactable at $a$,
with $\textup{k}(g.f)_a=
\textup{k}g_{a'}\textup{k}f_a$. 
\end{prop}

\proof
Since the centred maps $f:(E_a,U)\lra(E'_{a'},U')$ and\break
$g:(E'_{a'},U')\lra(E''_{a''},U'')$ are $\Sigma$-contactable (where
$a''=g(a')$),\break
$g.f$ is $\Sigma$-contactable (refer to 2.2.20), with
K$(g.f)=$K$g.$K$f$.
Thus $g.f$ is $\Sigma$-contactable at $a$ and, for all $x\in E$, we have
$g(a')+$k$(g.f)_a(x)=$
K$(g.f)(a+x)\!=\!$K$g($K$f(x+a))\!=\!
$K$g(a'+$k$f_a(x))\!=\!g(a')+
$k$g_{a'}
(($k$f_a)(x))$;
which ends the proof.
\cqfd

\begin{prop}
Let $(E,U)$, $(E_1,U_1)$, $(E_2,U_2)$ be three normed domains,
$a\in U$, and
$f_1:(E,U)\lra(E_1,U_1)$, $f_2:(E,U)\lra(E_2,U_2)$ two maps.
If $f_1$ and $f_2$ are $\Sigma$-contactable at $a$, 
then the pair\break
$(f_1,f_2):(E,U)\lra (E_1\times E_2,U_1\times U_2)$ 
is $\Sigma$-contactable at $a$,
with $\textup{k}(f_1,f_2)_a=(\textup{k}f_{1a},
\textup{k}f_{2a})$. 
\end{prop}

\proof
Since, for each $i\in\{1,2\}$, $f_i:(E_a,U)\lra(E_{ia_i},U_i)$ is
$\Sigma$-contactable (where $a_i=f_i(a))$), we know that the pair
$(f_1,f_2):(E_a,U)\lra(E_{1a_1}\times E_{2a_2},U_1\times U_2)=
((E_1\times E_2)_{(a_1,a_2)},U_1\times U_2)$ is\break 
$\Sigma$-contactable and that 
K$(f_1,f_2)=($K$f_1,$K$f_2)$, so that, for all $x\in E$,\break
we have, using 2.2.21,
$\ (a_1,a_2)+$k$(f_1,f_2)_a(x)=$
K$(f_1,f_2)(a+x)=\break
($K$f_1(a+x),$K$f_2(a+x))=$
\noindent $(a_1+$k$f_{1a}(x),a_2+$k$f_{2a}(x))=
(a_1,a_2)+\break
($k$f_{1a},$k$f_{2a})(x)$;
which ends the proof.
\cqfd

\begin{cory}
Let $(E_1,U_1)$, $(E'_1,U'_1)$, $(E_2,U_2)$, $(E'_2,U'_2)$ be  normed domains,
$a_1\in U_1$, $a_2\in U_2$ and
$f_i:(E_i,U_i)\lra(E'_i,U'_i)$ maps for each $i\in\{1,2\}$.
If each $f_i$ is $\Sigma$-contactable at $a_i$, 
then the product\break
$f_1\times f_2:(E_1\times E_2,U_1\times U_2)\lra 
(E'_1\times E'_2,U'_1\times U'_2)$ 
is $\Sigma$-contactable at $(a_1,a_2)$,
with $\textup{k}(f_1\times f_2)_{(a_1,a_2)}=
\textup{k}f_{1a_1}\times \textup{k}f_{2a_2}$. 
\end{cory}

\proof
Arises from 2.2.22 (proceeding as in the previous proof),
or from 2.3.10 and  2.3.11.
\cqfd

\begin{prop}
When $f:(E,U)\lra(E',U')$ is $\Sigma$-contactable at $a\in U$,
we have $\textup{k}f_a(x)=\lim_{0\not= v(t)\rightarrow 0}
\frac{f(a+v(t)x)-f(a)}{v(t)}$ for all $x\in E$.
\end{prop}

\proof
Since $f:(E_a,U)\lra(E'_{f(a)},U')$ is $\Sigma$-contactable, with 
$E'_{f(a)}$ revertible (see the recalls at the beginning of this section 2.3), we use 2.2.19 to write, for all $x\in E$: 
K$f(x)=\lim_{0\not=v(t)\rightarrow 0}t\buildrel{-1}\over\star f(t\star x)=
\lim_{0\not=v(t)\rightarrow 0}(f(a)+\frac{f(a+v(t)(x-a))-f(a)}{v(t)})$;
we merely  need to use the equality
k$f_a(x)=$K$f(a+x)-f(a)$ to obtain the wished result.
\cqfd

\begin{prop}
Let $f:(E,U)\lra(E',U')$ be a $\Sigma$-contactable map at $a\in U$; 
then $\rho(\textup{T}f_a)=\sup\{\frac{\|h(x)-h(y)\|}{\|x-y\|}\,|\, x,y\in B'(0,1), 
x\not= y\}$ where $h=\textup{k}f_a$ which is 
$\rho(\textup{T}f_a)$-lipschitzian (see  2.2.9).
\end{prop}

\proof
Arises from 2.2.9, applied to the $\Sigma$-homogeneous map 
K$f:E_a\lra E'_{f(a)}$, using the fact that T$f_a=
$T$($K$f)_a$ 
(since K$f\in $T$f_a$), and that
$x,y\in B'(0,1)$ iff $x+a,y+a\in B'(a,1)$.
\cqfd

\vspace{3mm}

We are going now to generalize the results given in the section 1.5
which compared the differentials to the tangentials.
Here, we compare $tf$ to $kf$. 

Let $(E,U),(E',U')$ be two normed domains; $U$, $U'$ non empty
(the canonical injections are still denoted $j:U\lra E$ and $j':U'\lra E'$).
Consider then the following composite
(denoted $J$):
$$\xymatrix{
\Sigma\!\!-\!\!\L\textup{Hom}(E,E')\!\ar[r]^>>>>>{j} &\!\!\J\textup{et}
((E,0)\!,\!(E',0))
\!\!\ar[r]^<<<<{can} &\!\! \J\textup{et}\jf(E\!,\! E')
\!\ar[r]^{\widehat{\Gamma}^{-1}} &\!\! \J\textup{et}\jf
(U\!,\! U')
}$$
See just before 2.2.5 where this $j$ had been denoted $can$, and 1.4.15, 1.4.28 for these $can$ and $\widehat\Gamma$;
by composition, this $J$ is an isometric embedding.
By definition, we have $J(h)=\kappa(j')^{-1}.[$T$h_0,0,0].\kappa(j)$,
so that $J(0)=O$, the free zero jet (see 1.4.14).
We set $Im(J)=\break
J(\Sigma\!-\!\L\textup{Hom}(E,E'))$.

\begin{prop}
Let $f:(E,U)\lra(E',U')$ be a map and $a\in U$; we assume that
$f$ is $Tang_a$. Then $f$ is $\Sigma$-contactable at $a$ iff
$\textup{t}f_a\in Im(J)$; in this case 
$\textup{t}f_a=J(\textup{k}f_a)$ ... see 1.5.1 for \textup{t}$f_a$).
\end{prop}

\proof
We proceed as in the proof of 1.5.4.
Indeed, if $[\varphi,a,a']\in\J$et$\jf(U,U')$, then 
$[\varphi,a,a']\in Im(J)$
iff there exists
$h\in\Sigma$-$\L$Hom$(E,E')$ such that
$\Theta_{aa'}(h)|_U\in j_{a'}.\varphi$. In this case, $J(h)=[\varphi,a,a']$.
\cqfd

\begin{cory}
If $f:(E,U)\lra(E',U')$ is $\Sigma$-contactable (i.e. at every $x\in U$),
then $f$ is tangentiable and this diagram commutes (where 
$\textup{k}f(x)=\textup{k}f_x$):
$$\xymatrix{
{}&U\ar[dl]_{\textup{k}f}\ar[dr]^{\textup{t}f}&{}\\
\Sigma\!-\!\L\textup{Hom}(E,E')\ar[rr]_J && 
\J\textup{et}\jf(U,U')
}$$
\end{cory}

\begin{prop}
Let $(E,U),(E',U')$ be two normed domains,\break
$f:(E,U)\lra(E',U')$
a continuous map, $a,b\in U$ such that $[a,b]\subset U$, $F$ a finite subset of $]a,b[$
and $r$ a positive real number. We assume that, for all $x\in\,]a,b[-F$, the map $f$
is $\Sigma$-contactable at $x$  and satisfies 
$\|\textup{k}f_x\|\leq r$. Then
we have: $\|f(b)-f(a)\|\leq r\|b-a\|$.
\end{prop}

\proof
Arises from the fact that, for every $x\in\,]a,b[-F$, we have:
$d($t$f_x,O)=d(J($k$f_x),J(0))=
d($k$f_x,0)=\|$k$f_x\|\leq r$ ... it remains now to use 1.5.9.
\cqfd

\vspace{4mm}
We precise here, for \texttt{the standard $\R$-contracting spaces}, what we have yet recalled and precised (around 2.3.1), for the canonical\break
$\R$-contracting spaces: for instance,
every n.v.s. $E$ can be considered as a
standard $\R$-contracting space
$E_a$ ($a$ being its central point)
whose standard external operation is
$t\star x=a+t(x-a)$ for all $t\in\R$ and 
$x\in E$ (see 2.1.8 and 2.1.10). 
Actually, referring to 2.1.14, we know that 
$E_a$ is a revertible standard
$\R$-contracting space (with
$t\buildrel{-1}\over\star x=a+t^{-1}(x-a)$ for $t\not=0$).
The standard $\R$-Lhomogeneous maps $h:E_0\lra E'_0$, are the lipschitzian maps $h:E\lra E'$ which verify, for all $t\in\R$ and 
$x\in E$:
$h(tx)=th(x)$; see 2.1.16 and 2.1.17.

The reader will verify that the whole study of the ``canonical n.v.s. world''
made in this section 2.3 is still valid for the ``standard n.v.s. world'';
in particular, we can speak of standard $\R$-contactable maps
at $a$, with 
standard $\R$-contact K$_\R f$ (as used in 2.2.17) and its translate
k$_\R f_a$ which is standard $\R$-Lhomogeneous ; 
(the index $\R$ being there to distinguish them from their canonical 
analogues, which were simply denoted K$f$ and k$f_a$, like for every valued monoid).
In particular, we have a $\R$-standard analogue of 2.3.7.

\section{G-differentiability}
Here is the first of the two special cases tackled in this paper
(in the vector space context).
It is the case where $\Sigma=\R_+$. As we have said in the general introduction, the young French Ren\' e Gateaux
was the first mathematician, at the beginning of the XXth century\footnote{Ren\' e Gateaux was one of the first victims of the first world war, he was twenty five years old when he died on the third of october 1914.},
to be interested in maps (said ``differentiable in the sense of Gateaux'')
which are very close to our  
$\R_+$-contactable maps:
the essential difference between these two notions being the fact that
we replace the continuous maps (see Bouligand in [7]) by lipschitzian maps.
In homage to Gateaux,  we will choose the
word ``G-differentiable'' for 
``$\R_+$- contactable''\footnote{it is shorter than ``lipschitzian Gateaux-differentiable'' which would be more convenient.}.
Then, we will study, more specially, the continuously 
G-differentiable maps: we will state that, in finite dimension,  this continuity property
forces the linearity of the
G-differential at every point.
\vspace{1mm}

As in section 2.3, we remain in the n.v.s. framework.
Now,\break
$E$ and $E'$ being two such n.v.s., we merely write
$\L$Hom$(E,E')$ for the n.v.s. which we should denote
$\R_+$-$\L$Hom$(E,E')$, which is in fact\break 
$\R_+$-$\C$ontr$(E_0,E'_0)$: see around 2.3.2;
we thus recall that its elements are the
$\R_+$-Lhomogeneous maps
$h:E_0\lra E'_0$, i.e the maps $h:E\lra E'$
which are lipschitzian and verify
$h(tx)=th(x)$ for all $t\in\R_+$ and $x\in E$.
Obviously, standard $\R$-Lhomogeneous implies
$\R_+$-Lhomogeneous (for the standard notions, refer to the end of section 2.3, 
just above). 

\vspace{3mm}
\begin{exams}
{}
\end{exams}

1) The continuous linear maps from $E$ to $E'$ are in 
$\L$Hom$(E,E')$. 

2) Let $E$ be a n.v.s.; then every norm $N$ on $E$, which is equivalent to the given norm $\|\ \|$ on $E$, is in $\L$Hom$(E,\R)$, since it is lipschitzian $(E,\|\ \|)\lra(\R,|\ |)$. 
In particular, $\vartheta:\R\lra\R:x\mapsto|x|$. In higher dimensions, 
the usual norms on $\R^n$ : $N^2(x)=\|x\|_2=\sqrt{{\sum}_{i=1}^n|x_i|^2}$, $N^\infty(x)=\|x\|_\infty=$Max$_i|x_i|$ and $N^1(x)=\|x\|_1=\sum_{i=1}^n|x_i|$.
If $n=2$, the graphs of these norms are cones; a ``circular cone'' for 
$N^2$ and ``square cones'' for $N^\infty$ and $N^1$, since 
they are generated, when $z$ varies in $\R_+$,
by their horizontals sections at height $(0,0,z)$: the circle of equation $x^2+y^2=z^2$ (for $N^2$), and the squares of equations 
Max$(|x|,|y|)=z$ and $|x|+|y|=z$ (for $N^\infty$ and $N^1$).
For $N^\infty$, see Fig 2.

\includegraphics{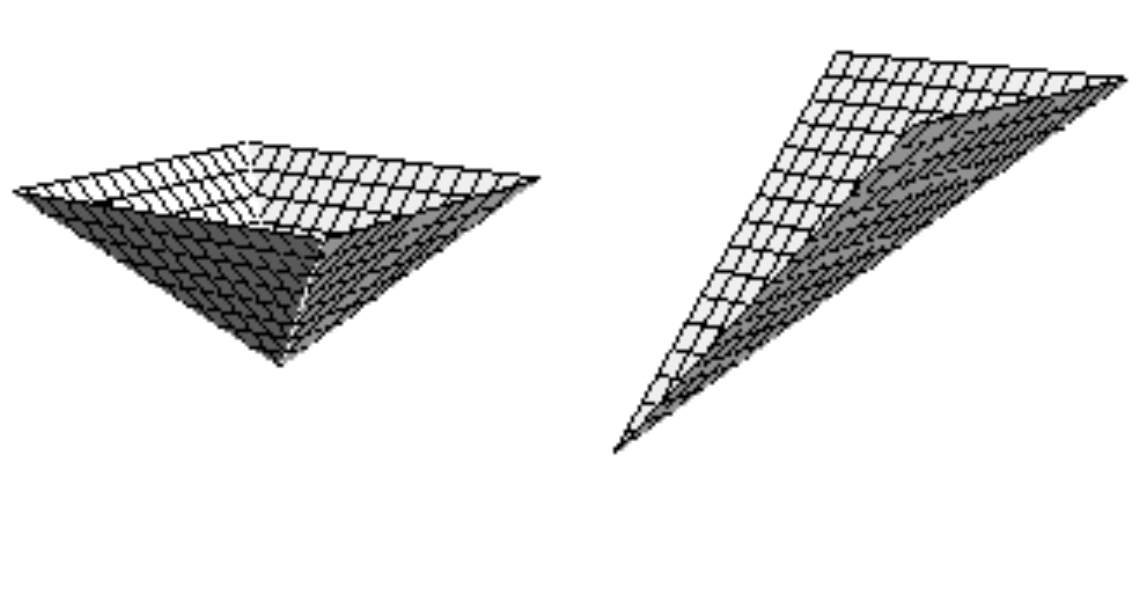}
\begin{picture}(0,0)
\put(-103,10){$N^\infty$}
\put(-33,10){Max}
\put(-72,5){Figure 2}
\end{picture}

\vfill\eject

3) The maps $\,$Max, Min: $\R^n\lra\R$ are in $\L$Hom$(\R^n,\R)$, since they are lipschitzian (indeed, for $x=(x_1,\cdots,x_n)$, $y=(y_1,\cdots,y_n)
\in\R^n$, let\break
$i,j\in\{1,\cdots,n\}$ such that $x_i=\textup{Max}(x)$ and
$y_j=\textup{Max}(y)$, then $\textup{Max}(x)-\textup{Max}(y)=
x_i-y_j\leq x_i-y_i\leq|x_i-y_i|\leq\|x-y\|_\infty$, the last term being the product distance on $\R^n$; similarly, we have also
$\textup{Max}(y)-\textup{Max}(x)\leq\|y-x\|_\infty=\|x-y\|_\infty$;
so that $|\textup{Max}(x)-\textup{Max}(y)|\leq\|x-y\|_\infty$).
If $n=2$, the graph of Max is generated, when $z$ varies in $\R$, by its horizontal section at 
height $(0,0,z)$ which is the right angle of respective sides  
$\{(z,y,z)\,|\, y\leq z\}$ and $\{(x,z,z)\,|\, x\leq z\}$, its vertex $(z,z,z)$ varying on the axis 
$\widehat\Delta$ of equation $z=y=x$.
Thus the graph of Max is the union of two halph-plans of $\R^3$, that we could call an ``angular cone'' (see Fig 2). 
The graph of Min is the ``angular cone'' obtained from the previous one by the reflection
with respect to their common axis $\widehat\Delta$.

\vspace{3mm}

\begin{prop}
Let $h\in\L\textup{Hom}(E,E')$, where $E\not=\{0\}$; we have\break
$\|h\|=\sup\{\|h(x)\|\, |\, \|x\|=1\}$.
\end{prop}

\proof
Let us set $k=\sup\{\|h(x)\|\,|\, \|x\|=1\}$ and take $x\in E$;
if $\|x\|=1$, then, using 2.3.2,  $\|h(x)\|=C(x)\leq \|h\|$, so that 
$k\leq\|h\|$. If $\|x\|\leq 1$ with $x\not=0$, then
$C(x)=\frac{\|h(x)\|}{\|x\|}=\|h(\frac{x}{\|x\|})\|\leq k$, so that $\|h\|\leq k$.
\cqfd

\begin{prop}
Let $E$ be a n.v.s.; then the following map is a linear isometry:
\quad$can:\L\textup{Hom}(\R,E)\lra E\times E:h\mapsto (-h(-1),h(1))$.
\end{prop}

\proof
The linearity of $can$ is obvious. As for the fact that $can$ is bijective:
if $(v,w)\in E\times E$, then $can^{-1}(v,w)=h$ where 
$h(t)=tw$ if $t\geq 0$
and $h(t)=tv$ if $t\leq 0$ (this $h$ is in $\L$Hom$(\R,E)$: $h$ is
$k$-lipschitzian with $k=\textup{Max}(\|v\|,\|w\|)$, since 
$\|h'(t)\|\leq\textup{Max}(\|v\|,\|w\|)$ for all $t\not=0$). Besides, $can$ is an isometry since, by 2.4.2, $\|h\|=\sup\{\|h(x)\|\ |\ |x|=1\}\break
=\textup{Max}(\|h(1)\|,\|h(-1)\|)=\|can(h)\|$.
\cqfd

\begin{prop}
Let $h\in\L\textup{Hom}(E,E')$, then, for all $\varepsilon>0$, 
we have 
$\rho(\textup{T}h_0)=\sup\{\frac{\|h(x)-h(y)\|}{\|x-y\|}\, |\, x\not= y;\ x,y\in C(\varepsilon)\}$
(see 1.2.8 for the lipschitzian ratio $\rho$), where 
$C(\varepsilon)=\{x\in E\, |\, 1-\varepsilon<\|x\|<
1+\varepsilon\}$. 
\end{prop}

\proof
Let us set $r(\varepsilon)=\sup\{\frac{\|h(x)-h(y)\|}{\|x-y\|}\, |\, x\not= y;\ x,y\in C(\varepsilon)\}$. It is clear that $r(\varepsilon)\leq\rho($T$h_0)$ since $h$ is $\rho($T$h_0)$-lipschitzian (by 2.3.14).
Conversely, take $a\in E$, $a\not=0$ and let us set 
$\varepsilon'=\|a\|\varepsilon$. Then, for every $x\in B(a,\varepsilon')$,
we have $\|\frac{x}{\|a\|}\|\leq
\frac{\|x-a\|+\|a\|}{\|a\|}=1+\frac{\|x-a\|}{\|a\|}<1+\frac{\varepsilon'}{\|a\|}=
1+\varepsilon$
and $\|\frac{x}{\|a\|}\|\geq \frac{\|a\|-\|x-a\|}{\|a\|}=1-\frac{\|x-a\|}{\|a\|}>
1-\frac{\varepsilon'}{\|a\|}=1-\varepsilon$; thus,
$\frac{x}{\|a\|}\in C(\varepsilon)$, $\frac{a}{\|a\|}\in C(\varepsilon)$. 
Then $r(\varepsilon)\geq
\frac{\|h(\frac{x}{\|a\|})-h(\frac{a}{\|a\|})\|}{\|\frac{x}{\|a\|}-\frac{a}{\|a\|}\|}=
\frac{\|h(x)-h(a)\|}{\|x-a\|}$ if $x\not=a$; so that $h$ is $r(\varepsilon)$-$LSL_a$.
In particular, for every $a\in B'(0,1)-\{0\}$, $h$ is $r(\varepsilon)$-
$LSL_a$, with $B'(0,1)$ convex. Thus, $h|_{B'(0,1)}$ is 
$r(\varepsilon)$-lipschitzian (by 1.1.23), so that
$\rho($T$h_0)\leq r(\varepsilon)$.
\cqfd

\begin{defi}
Let $(E,U)$,$(E',U')$ be two normed domains and\break
$f:(E,U)\lra(E',U')$ a map (refer just before 2.3.6). We say that $f$ is 
\textup{G-differentiable} at $a$ if
$f$ is $\R_+$-contactable at $a$ (see 2.3.6);
in this case, $\textup{K}_+f$ and $\textup{k}_+f_a$ will 
respectively merely denote $\textup{K}_{\R_+}f$ and $\textup{k}_{\R_+}f_a$.
\end{defi}

\begin{remks}
In each case, $E$ and $E'$ are n.v.s. and $f:U\lra E'$ where $U$ is an open subset of $E$, and $a\in U$.

1) If $f$ is differentiable at $a$, then $f$ is G-differentiable at $a$ with $\textup{k}_+f_a=\textup{d}f_a$ (see 2.3.8).

2) If $f$ is G-differentiable at $a$, then $f$ is differentiable at $a$ iff\break 
$\textup{k}_+f_a$ is linear 
(see 2.3.8).

3) If $f$ is G-differentiable at $a$, then $f$ is standard $\R$-contactable at $a$ iff \textup{k}$_+ f_a$ is standard $\R$-homogeneous with (referring to the end of section 2.3 for the standard notations) \textup{k}$_+ f_a=\textup{k}_\R f_a$
(see 2.3.7, since, referring to 2.1.3, the canonical injection $\R_+\hookrightarrow\R$ is a morphism of valued monoid).

4) Here, $E=\R$; then
$f$ is standard $\R$-contactable at $a$ iff $f$ is differentiable at $a$
with \textup{k}$_\R f_a=\textup{d}f_a$
(it comes from the fact that, for $h:\R\lra E$,  standard
$\R$-homogeneous $\ =\ $ linear!). 
When $E\not=\R$, see a counter-example in 2.5.12.

5) If $f$ is G-differentiable at $a$, then $f$ is 
canonically  $\Sigma$-contactable at $a$ for all valued monoid 
$\Sigma$,
with $\textup{k}_\Sigma f_a=\textup{k}_+f_a$
(see 2.3.7, since, referring to 2.1.5, every valuation is a morphism of valued monoid). However,\break 
G-differentiable
$\ \not\!\!\Lra$ standard $\ \R$-contactable at $a$: see 2.5.11.

\end{remks}

\vspace{3mm}
\begin{exams}
{}
\end{exams}
The following examples are all G-differentiable at every point:

1) The norm function $\vartheta:\R\lra\R:x\mapsto|x|$, which verifies
k$_+\vartheta_0=\vartheta$ and, for $a\not=0$, k$_+\vartheta_a=$d$\vartheta_a=sign(a)Id_\R$ where $sign(a)=\frac{a}{|a|}$.

2) The euclidian norm function $N^2:\R^n\lra\R:x\mapsto\|x\|_2$,
which verifies 
k$_+N^2_0=N^2$ and, for $a\not=0$, k$_+N^2_a=$d$N^2_a$.

3) The functions $\,$Max, Min:$\ \R^n\lra\R$, which verify:\par 
\noindent k$_+$Max$_a(x)=$Max$_{i\in\overline{\mathcal I}(a)}(x_i)$ and
k$_+$Min$_a(x)=$Min$_{i\in\underline{\mathcal I}(a)}(x_i)$, where
$\overline{\mathcal I}(a)=\break
\{i\in\{1,\dots,n\}\,|\, a_i=$Max$(a)\}$ and $\underline{\mathcal I}(a)=\{i\in\{1,\dots,n\}\,|\, a_i=$Min$(a)\}$;
which provides k$_+$Max$_0=$Max and k$_+$Min$_0=$Min.

4) The product norm function $N^\infty:\R^n\lra\R:x\mapsto\|x\|_\infty$,
which verifies 
k$_+N^\infty_0=N^\infty$ and, for $a\not=0$,
k$_+N^\infty_a(x)=$Max$_{i\in\overline{\mathcal I}(|a|)}(sign(a_i)x_i)$,
where $|a|=(|a_1|,\dots,|a_n|)$.

5) The norm function $N^1:\R^n\lra\R:x\mapsto\|x\|_1$, which veri-\break
fies 
k$_+N^1_a(x)=\sum_{i\in{\mathcal I}_0(a)}|x_i|+
\sum_{i\in{\mathcal I}_1(a)}sign(a_i)x_i$, where
${\mathcal I}_0(a)=\break
\{i\in\{1,\dots,n\}\,|\, a_i=0\}$ and
${\mathcal I}_1(a)=\{i\in\{1,\dots,n\}\,|\, a_i\not=0\}$;
which gives k$_+N^1_0=N^1$. 

\vspace{2mm}
\proof
First, all the above functions $h$ being $\R_+$-Lhomogeneous (see  2.4.1), they are all G-differentiable at 0 with k$_+h_0=h$ (see 2.3.7);
we thus have here easy functions which, however not differentiable at 0, are G-differentiable at 0, that is, instead of having at 0, a continuous linear tangent function, are their own tangent at 0 (i.e their $\R_+$-contact at 0), since they are 
$\R_+$-Lhomogeneous.
We use also 2.4.6.

1) Comes from the preceding lines. 
 
2) Still obvious for the same reason as the above 1).

3) We will show that, for all $a$, locally in a neighborhood $V$ of 0, we have,
for all $x\in V$, the equality
Max$(a+x)=$Max$(a)+$Max$_{i\in\overline{\mathcal I}(a)}x_i$
which will prove the foretold result for the map Max (since 
Max$_{i\in\overline{\mathcal I}(a)}$ is $\R_+$-Lhomogeneous).
Indeed:
 
\noindent - If $\overline{\mathcal I}(a)=\{1,\cdots,n\}$ (i.e, if $a$ is a constant 
$n$-uplet), then, for all $x\in\R^n$, we have 
Max$(a+x)=$Max$(a)+$Max$(x)=$Max$(a)+$
Max$_{i\in\overline{\mathcal I}(a)}x_i$. 

\noindent - If $\overline{\mathcal I}(a)\not=\{1,\cdots,n\}$.
Let $j\in\overline{\mathcal I}(a)$ such that
Max$_{i\in\overline{\mathcal I}(a)}(x_i)=x_j$\break
(by definition, this $j$ verifies
Max$(a)=a_j$);
for this $j$, we can write:\par
\noindent Max$(a)+$
Max$_{i\in\overline{\mathcal I}(a)}(x_i)=a_j+x_j\leq$Max$(a+x)$.
Conversely, 
we set
$r=\frac{1}{2}$Min$_{i\notin\overline{\mathcal I}(a)}($Max$(a)-a_i)$; 
then $r>0$.
Let $V$ be the open ball $B_\infty(0,r)$ (for the product norm $\|\ \|_\infty$); then, for 
$x\in V$, we have:

\quad - if $j\not\in
\overline{\mathcal I}(a)$,
 $\ x_j-$Max$_{i\in\overline{\mathcal I}(a)}(x_i)\leq
|x_j-$Max$_{i\in\overline{\mathcal I}(a)}(x_i)|
\leq|x_j|+
|$Max$_{i\in\overline{\mathcal I}(a)}(x_i)|\leq 2\|x\|_\infty<2r\leq$
Max$(a)-a_j$, and thus again\break
$a_j+x_j\leq$Max$(a)+$
Max$_{i\in\overline{\mathcal I}(a)}(x_i)$,

\quad - if $j\in\overline{\mathcal I}(a)$, 
$\ a_j+x_j=$Max$(a)+x_j\leq$
Max$(a)+$Max$_{i\in\overline{\mathcal I}(a)}(x_i)$,\break
 and thus 
Max$(a+x)\leq$Max$(a)+$
Max$_{i\in\overline{\mathcal I}(a)}(x_i)$.

Same for the function Min.

4) Since $N^\infty\!\!=$Max$.\vartheta^n$, 
(where $\vartheta^n=\vartheta\times\cdots\times\vartheta$, $n$ times),               
the function $N^\infty\!$ is G-differentiable, using 2.3.10 and 2.3.12.
For $a\not=0$,
we can write k$_+N^\infty_a(x)=$
k$_+$Max$_{|a|}($k$_+\vartheta_{a_1}(x_1),\cdots,$k$_+\vartheta_{a_n}(x_n))=$Max$_{i\in\overline{\mathcal I}(|a|)}($k$_+\vartheta_{a_i}(x_i))\break
=$
Max$_{i\in\overline{\mathcal I}(|a|)}(sign(a_i)x_i)$, since, for all 
$i\!\in\!\overline{\mathcal I}(|a|)$, we have $a_i\not=0$. This provides
k$_+N_a^\infty=$Max on $\{a\in\R^n| a_1=\cdots =a_n\!>\!0\}$.

5) Since $N^1=\sigma.\vartheta^n$, where $\sigma$ is the addition of $\R^n$ (which is linear and thus differentiable, with
k$_+\sigma_b=$d$\sigma_b=\sigma$ for all $b$), $N^1$
 is\break
G-differentiable with, for all $a\in\R^n$, 
k$_+N^1_a\,=\,$k$_+\sigma_{|a|}.($k$_+\vartheta_{a_1},\cdots,$k$_+\vartheta_{a_n})$
$=\sigma($k$_+\vartheta_{a_1},\cdots,$k$_+\vartheta_{a_n})$, so that, for every $x\in\R^n$, we have k$_+N^1_a(x)=\sum_{i=1}^n $k$_+\vartheta_{a_i}(x)=
\sum_{i\in{\mathcal I}_0(a)}|x_i|
+
\sum_{i\in{\mathcal I}_1(a)}sign(a_i)x_i$.
\cqfd

\vspace{2mm}
\begin{remks}\par\hfill

1) We are giving here, for each G-differentiable function $h$ studied in  2.4.7, 
the domain $D(h)$ on which $h$ is differentiable (using 2.4.6).
Anyhow, in each case, we give the $\R_+$-contact (we often speak in terms of graphs for $n=2$).

\qquad a) $D(\vartheta)=\R^*$; we can also use 2.4.9.

\qquad b) $D(N^2)=\R^n-\{0\}$.
Actually, at every point $a\in\R^n-\{0\}$, the graph of $N^2$ 
has a tangent plan; at 0, this graph is its own tangent, i.e $\R_+$-contact ($N^2$ being
$\R_+$-Lhomogeneous);
which can be easily ``seen'' for $n=2$.

\qquad c) First, the formula established in the proof of 3) of 2.4.7
can be locally written \textup{Max}$(x)=\textup{Max}(a)+$
\textup{k}$_+\textup{Max}_a(x-a)$, which gives 
\textup{K}$_+\textup{Max}_a$, the $\R_+$-contact of \textup{Max} at $a$
(as the translate of \textup{k}$_+\textup{Max}_a$ at $a$).
In particular, if $\overline{\mathcal I}(a)=\{1,\dots,n\}$, then
\textup{Max}$_{i\in\overline{\mathcal I}(a)}=\textup{Max}$, so that
\textup{k}$_+$\textup{Max}$_a=\textup{Max}=$
\textup{K}$_+\textup{Max}_a$, and this, at every point of the axis
$\Delta=\{a\in\R^n\,|
a_1=\cdots=a_n\}$. 
Now, let us notice that the map
\textup{Max} $:\R^n\lra\R$ is linear iff $n=1$.
Thus, we obtain
$D(\textup{Max})=\break
\{a\in\R^n\,|\, \exists ! i\leq n\ (a_i=\textup{Max}(a))\}$.
For instance, for $n=2$, D(\textup{Max})=$\Delta^c$ where
$\Delta=
\{a\in\R^2\,|\, a_1=a_2\}=
\{a\in\R^2\,|\, \overline{\mathcal I}(a)=\{1,2\}\}$;
so D(\textup{Max}) is the disjoint union of the two open half-plans $\Pi_1$ and 
$\Pi_2$ of $\R^2$, of respective inequations $a_1>a_2$ and $a_1<a_2$,
on which the restrictions \textup{Max}$|_{\Pi_i}$ are the $p_i|_{\Pi_i}$ (where the $p_i$ are the canonical projections),
so that, if $x\in\Pi_i$, \textup{k}$_+\textup{Max}_x=$
\textup{d}$\textup{Max}_x=p_i$.
We have yet studied the case of the points on $\Delta$ just above
(refer to 2.4.1 for this ``angular cone'').

\qquad d) Here $D(N^\infty)=\{a\in\R^n\,|\, \exists ! i\leq n\ (|a_i|=\|a\|_\infty)\}$.
For $n=2$, $D(N^\infty)=(\{a\in\R^2\,|\, |a_1|=|a_2|\})^c$, so that
$N^\infty$ is differentiable on the four open quarter-plans of $\R^2$ delimited by the two diagonals of $\R^2$
(we have recalled in 2.4.1 that the graph of $N^\infty$ is a square-cone). 
This cone admits a tangent plan on each of its sides; besides,
at 0, the\break
$\R_+$-contact is \textup{k}$_+N^\infty_0=N^\infty$ itself. As for the ridge lines, 
for example, the 
$\R_+$-contact on  
$N^\infty|_{\Delta_+}$, where $\Delta_+=\{a\in\R^2\,|\, a_1=a_2>0\}$,
is \textup{Max}!

\qquad e) For $N^1$, we have $D(N^1)=\{a\in\R^n\, |\, \forall i\leq n\ \, a_i\not=0\}$.
For $n=2$, The graph of $N^1$ is obtained, 
from the one of $N^\infty$, with the help of the rotation of axis 
``$x=y=0$'', and of angle $\pi/4$ (the two axis ``$y=0$''
and ``$x=0$'' in $\R^2$ playing here the part of the two diagonals for $N^\infty$).

2) Of course, there exist G-differentiable maps which are not\break
$\R_+$-homogeneous:
2.4.7 is crowded with such examples.
Indeed, a translate $g(x)=f(a+x)$ of a $\R_+$-homogeneous map
$f$ is not necessarily still $\R_+$-homogeneous, although,
in our examples, such a translate remains $G$-differentiable (by
composition). 
\end{remks}
\vspace{2mm}

\begin{prop}
Let $f:(E,U)\lra(E',U')$ be a map between normed domains, where $E=\R$ and 
$a\in U$. Then $f$ is G-differentiable iff $f$ admits left and right derivatives at $a$. In this case, referring to 2.4.3 for the linear isometry $can$, we have
\textup{k}$_+f_a=can^{-1}(f'_l(a),f'_r(a))$. 
\end{prop}

\proof
Let $h\in\L$Hom$(\R,E')$ and $(v,w)=can(h)$, i.e.  $v=-h(-1)$ and
$w=h(1)$.
Then, if (referring to 2.3.1 for the notations) we denote 
$k=\Theta_{af(a)}(h)$,
$V=U\cap\!\!\buildrel{-1}\over k\!\!(U')$
and $\underline k$ the restriction of $k$ to $V\lra U'$, we have the following equivalences:
$\ f|_V\tang_a\underline k \ \Llra$\par
\noindent$\lim_{a\not=x\rightarrow a}\frac{\|f(x)-f(a)-h(x-a)\|}{|x-a|}=0\! \Llra$
\noindent $\lim_{a\not=x\rightarrow a}\|\frac{f(x)-f(a)}{x-a}-\frac{h(x-a)}{x-a}\|=0
\Llra
(\lim_{a<x\rightarrow a}\|\frac{f(x)-f(a)}{x-a}-\frac{h(x-a)}{x-a}\|=0$ and
$\lim_{a>x\rightarrow a}\|\frac{f(x)-f(a)}{x-a}-\frac{h(x-a)}{x-a}\|=0)$ $\Llra
(\lim_{a<x\rightarrow a}\|\frac{f(x)-f(a)}{x-a}-w\|=0$ and
$\lim_{a>x\rightarrow a}\|\frac{f(x)-f(a)}{x-a}-v\|=0)$ $\Llra
(f'_r(a)=w$ and $f'_l(a)=v$).
Thus $h=$k$_+f_a$.
\cqfd

\vspace{5mm}
\centerline{\texttt{Continuous G-differentiable maps}}
\vspace{1mm}
Now, our aim is to prove that, in finite dimension, 
every \textit{continuous G-differentiable} map $f:(E,U)\lra(E',U')$ (i.e 
G-differentiable such that
k$_+ f:U\lra\L$Hom$(E,E')$ is continuous) is in fact of class $C^1$.
In order to prove it, we need some preliminary results. Let us begin by a particular presentation of the mean value theorem.

\begin{prop} \textup{(well-known)}

\noindent Let $f:[a,b]\lra\R$ be a continuous function which admits a right derivative at every point of $]a,b[$, and $k\in\R$. Then,

1) If, for all $t\in\,]a,b[$ $f'_r(t)\leq k$, then $f(b)-f(a)\leq k(b-a)$,

2) If, for all $t\in\,]a,b[$ $f'_r(t)\geq k$, then $f(b)-f(a)\geq k(b-a)$.
\end{prop}

\proof
2) comes from 1), considering the function
$[a,b]\lra\R:\break
t\mapsto -f(t)$.
\cqfd

\begin{cory}
Let $f:[a,b]\lra\R$ be a continuous function which admits a left derivative at every point of $]a,b[$, then:

1) If, for all $t\in\,]a,b[$ $f'_l(t)\leq k$, then $f(b)-f(a)\leq k(b-a)$,

2) If, for all $t\in\,]a,b[$ $f'_l(t)\geq k$, then $f(b)-f(a)\geq k(b-a)$.
\end{cory}

\proof
We just apply 2.4.10 to the function $[-b,-a]\lra\R:\break
t\mapsto f(-t)$.
\cqfd

\begin{prop}
Let $U$ be an open subset of $\R$ and $f:U\lra\R$ a continuous function admitting left and right derivatives at every point of $U$ and such that 
the functions $f'_l,f'_r:U\lra\R$ are continuous at $a\in U$.
Then $f'_l(a)=f'_r(a)$, so that $f$ is derivable at $a$.
\end{prop}

\proof
Let $\varepsilon>0$; let us prove that $f'_l(a)\leq f'_r(a)+\varepsilon$.
Since $f'_r$ is continuous at $a$, there exists $\eta>0$ such that
$]a-\eta,a+\eta[\,\subset U$ and $f'_r(x)<f'_r(a)+\varepsilon$ for all $x\in\,]a-\eta,a+\eta[$. Let $k=f'_r(a)+\varepsilon$. Then,\break
2.4.10 provides that, for all $x\in\R$ verifying $a-\eta<x<a$, 
we have\break
$f(a)-f(x)\leq k(a-x)$,
i.e $\frac{f(a)-f(x)}{a-x}\leq k$, which gives $f'_l(a)\leq k$ (doing $x\rightarrow a$); hence $f'_l(a)\leq f'_r(a)$, doing $\varepsilon\rightarrow 0$. Same for the reverse inequality applying 2.4.11.
\cqfd

\begin{prop}
Let $f:(E,U)\lra(\R,U')$ be a function which is\break
G-differentiable on $U$.
Besides, we assume that the map
$\textup{k}_+f:\break
U\lra\L\textup{Hom}(E,\R):x\mapsto \textup{k}_+f_x$ is continuous
at a point $a\in U$. Then, for all $v\in E$, the directionnal derivative at $a$ $\frac{\partial f}{\partial v}(a)=\lim_{0\not=t\rightarrow 0}
\frac{f(a+tv)-f(a)}{t}$ exists in $\R$.
\end{prop}

\proof
Let $v\in E$, $\alpha:\R\lra E:t\mapsto a+tv$ and 
$V=\,\buildrel{-1}\over\alpha(U)$.
The composite $(\R,V)\buildrel\alpha\over\lra(E,U)
\buildrel f\over\lra(\R,U')$, denoted $g$, is G-differentiable on $V$,
thus continuous on $V$, so that $g$ admits left and right derivatives at every point of $V$ (by 2.4.9);
actually, , we have $g'_r(t)=$k$_+g_t(1)$ and $g'_l(t)=-$k$_+g_t(-1)$ 
for all $t\in V$.
The map k$_+g:
V\lra\L$Hom$(\R,\R)$ is continuous at 0, since it is the following composite:
$$\xymatrix{
V\,\ar[r]^{\!\!\!\!\!\!\!\!(Id,\alpha)}&\,V\times U\,
\ar[r]^{\!\!\!\!\!\!\!\!\!\!\!\!\!\!\!\!\!\!\!\!\!\!\!\!\!\!
\textup{d}\alpha\times \textup{k}_+f}&\,L(\R,E)
\times\L\textup{Hom}(E,\R)\,
\ar[r]^{\,\,\,\,\,\,\,\ \ \ \ comp}&\,\L\textup{Hom}(\R,\R)
}$$
(Refer to 2.2.6 where we have seen that $comp$ is $LSL$, thus continuous). In consequence, $g'_l,g'_r:V\lra\R$ are continuous at 0 (we compose k$_+g$ with the isometry
$can:\L$Hom$(\R,\R)\lra\R^2$ given in 2.4.3).
It arises from this that, thanks to 2.4.12, the map $g$ is derivable at 0, which means that $\frac{\partial f}{\partial v}(a)$ exists in $\R$.
\cqfd

\begin{theo}
Let $(E,U)$, $(E',U')$ be normed domains, where $E$ and $E'$ have a finite dimension, and $f:(E,U)\lra(E',U')$ a map which is G-differentiable on $U$ and such that the map $\textup{k}_+f:\break
U\lra\L\textup{Hom}(E,E')$ is continuous everywhere. Then, $f$ is of class $C^1$. 
\end{theo}

\proof
1) We assume here that $E'=\R$. Then, for all $x\in U$ and\break
$v\in E$, we have k$_+f_x(v)=\lim_{0\not=t\rightarrow 0}
\frac{f(x+tv)-f(x)}{t}=\frac{\partial f}{\partial v}(x)$, the two equalities arising respectively from 2.3.13  and 2.4.13.
Thus, for all $v\in E$, the function $\frac{\partial f}{\partial v}:U\lra\R$
is continuous (since it is the following composite
$U\buildrel{\textup{k}_+f}\over\lra\L\textup{Hom}(E,\R)\buildrel{\widehat v}\over\lra
\R$, where $\widehat v:h\mapsto h(v)$ is continuous since,
by 2.3.4, we have
$|\widehat v(h)|=|h(v)|\leq\|h\|\,\|v\|$), which proves that\break
$f$ is of class $C^1$.

2) Now, if $E'$ is of dimension $n$, we know that $E'$
is isomorphic to $\R^n$ (in the category of the n.v.s.);  this allows to 
come to the same end: indeed, if $\gamma:E'\lra\R^n$ is such an isomorphism and $p_i:\R^n\lra\R$\break
is the $i^{th}$ canonical projection, then
the following composite $f_i$ is\break
G-differentiable:
$U\buildrel f\over\lra U'\hookrightarrow E'\buildrel\gamma\over\lra\R^n
\buildrel{p_i}\over\lra\R$, and the following diagram commutes
(since $\gamma$ and $p_i$ are continuous linear maps):

$$\xymatrix{
{}&U\ar[dl]_{\textup{k}_+f}\ar[dr]^{(\textup{k}_+f_i)_{i\leq n}}&{}\\
\L\textup{Hom}(E,E')\ar[r]^{\sim}&\L\textup{Hom}(E,\R^n)\ar[r]^{\sim}&
\L\textup{Hom}(E,\R)^n
}$$
where the horizontal maps are continuous (since $\R_+$-$\C$ontr
is a cartesian category enriched in $\M$et).
Thus, the map k$_+f_i$ is continuous, so that the map $f_i$ is of class $C^1$, by the above 1). Using now the fact that $f=\gamma^{-1}.(f_1,\dots,f_n)$,
with $(f_1,\dots,f_n):U\lra\R^n$, we finally obtain that $f$
is itself of class $C^1$.
\cqfd

%% file: jet_II3.tex
\section{Fractality and neo-fractality}
This second particular case (still in the n.v.s. context; with\break $\Sigma=\N'_r$: see examples 2.1.2 and 2.1.10) has an immediate interest, since it speaks of fractality. We give here a construction process of such fractal maps which allows, on the way, to notice that the vector space  
of the linked metric
jets: $(\R,0)\lra(\R,0)$ is of infinite dimension;
which shows the great complexity of this new type of jets.

\vspace{1mm}
As in sections 2.3 and 2.4, we remain in the n.v.s. context.
Let us fix a real number $0<r<1$. Now, $E,E'$ being two n.v.s., 
we recall that the
$\N'_r$-Lhomogeneous maps $h:E_0\lra E'_0$ are the maps $h:E\lra E'$ which are lipschitzian and which satisfy the following fractality property:
$h(rx)=rh(x)$ for all $x\in E$;
such maps will be called\break
``$r$-Lfractal''. Thus,
we merely write
$r$-$\L$Frac$(E,E')$ for the n.v.s. which we should denote 
$\N'_r$-$\L$Hom$(E,E')$, which is in fact
$\N'_r$-$\C$ontr$(E_0,E'_0)$: see around 2.3.2.

The valuation $v_r:\N'_r\lra\R_+$ being a morphism of valued monoids, it provides the following inclusion:
$\L$Hom$(E,E')\subset r$-$\L$Frac$(E,E')$ ... see the beginning of section 2.4 for the notation $\L$Hom$(E,E')$.

\begin{remk}\par\hfill

1) The definition of an $r$-Lfractal map is the good one, since it provides 
$h(r^nx)=r^nh(x)$ for all $n\in\N$ and $x\in E$; and also $h(0)=0$ by continuity at 0.
And thus, referring to 2.1.10,
$h(n\star x)=n\star h(x)$ and $h(\infty\star x)=\infty\star h(x)$.

2) When $h$ is $r$-Lfractal, we also have $h(r^{-1}x)=r^{-1}h(x)$, since $h(x)=h(rr^{-1}x)=rh(r^{-1}x)$
... see again 2.1.17.

3) Why fractal? Merely because we have the equivalence:\break
$(x,y)\in Graph(h)$ iff $(rx,ry)\in Graph(h)$, which means that 
$Graph(h)$ remains identical to itself when we zoom into 
0 with a ratio $r$ (and we can iterate this process as many times as we want!).
We can have an approximative idea of a
fractal function $h:\R\lra\R$, by considering 0 as a point at the infinity (i.e at the unreachable horizon point),
the graph of $h$, being then seen in perspective, infinitely decreasing towards this horizon point, and still remaining itself, but thinner and thinner.
\end{remk}

\vspace{3mm}
\begin{exams}
{}
\end{exams}

1) The $\R_+$-Lhomogeneous maps $E_0\lra E'_0$ are $r$-Lfractal.

2) Consider the function $f:\R\lra\R$ defined by $f(0)=0$ and 
$f(x)=x\sin\log |x|$ for all $x\not=0$. Then, $f$ is $r$-Lfractal for 
$r=e^{-2\pi}$.

3) More generally, for $p\in\{1,2,\infty\}$, the map
$f^p:\R^n\lra\R^n:\break
x\mapsto \lambda_p(x)x$, where 
$\lambda_p:\R^n\lra\R$ is the function defined by $\lambda_p(0)=0$ and
$\lambda_p(x)=\sin\log\|x\|_p$ for $x\not=0$.
Then $f^p$ is $r$-Lfractal for $r=e^{-2\pi}$.

4) $\K$ being the triadic  Kantor set, let $K_\infty=\cup_{n\in\N}3^n\K$
and\break
$g:\R\lra\R:x\mapsto d(x,K_\infty)$.
Then $g$ is $\frac{1}{3}$-Lfractal.

\vspace{2mm}
\proof
1) Comes from the inclusion noticed before 2.5.1.

2) $f$ is continuous on $\R$ (for it is $SL_0$); it is derivable on 
$\R^*=\R-\{0\}$ with $f'(x)=\sin\log|x|+\cos\log|x|$,
so that $|f'(x)|\leq 2$. By the mean value theorem, we know that $f$ is 
2-lipschitzian. It remains to verify that $f(e^{-2\pi}x)=e^{-2\pi}f(x)$, which is proved below for $p\in\{1,2,\infty\}$ (see Fig 3).

\includegraphics{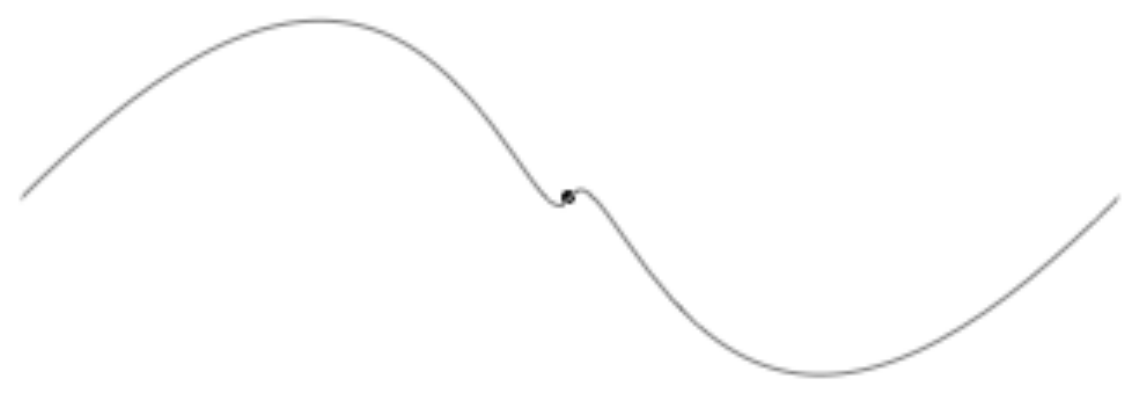}
\begin{picture}(0,0)
\put(-96,6){Figure 3}
\end{picture}

3) Let us set $r=e^{-2\pi}$.
For all $x\in\R^n-\{0\}$, we have $\lambda_p(rx)=
\sin\log\|rx\|_p=\sin\log r\|x\|_p=\sin(\log r+\log\|x\|_p)=
\sin(-2\pi+\log\|x\|_p)\break
=\sin\log\|x\|_p=\lambda_p(x)$, so that
$f^p(rx)=rf^p(x)$.

Thus, it remains to prove that $f^p$ is lipschitzian; we will prove it 
(with the help of the norm $\|\ \|_p$ on $\R^n$) using 
2.3.17 (with $\Sigma=\R_+$). First, $f^p$ is continuous (for it is $SL_0$);  
we also notice that, by composition, $\lambda_p$ is\break
G-differentiable
at every $x\not=0$ (since we have seen in the examples 2), 4) and 5) of 
2.4.7 that the norms $N^p=\|\ \|_p$ are G-differentiable on $\R^n$),
with $\textup{k}_+\lambda_{px}(y)=\frac{\cos\log\|x\|_p}{\|x\|_p}\, 
\textup{k}_+N^p_x(y)$
for all $y\in \R^n$, using 
the formula of\break
2.3.10 and 2.4.6.
Furthermore, if we notice that $f^p$ is the composite
$\R^n\buildrel{(\lambda_p,Id)}\over\lra\R\times\R^n
\buildrel b\over\lra\R^n$ where $b(t,x)=tx$, we deduce that $f^p$ is\break
G-differentiable at every $x\not=0$ (the map $b$ being bilinear,
it is differentiable at every $(t,x)\in\R\times\R^n$, and thus 
G-differentiable at such $(t,x)$) and, using also 2.3.11, we obtain, for $t=\lambda_p(x)$,
$\textup{k}_+f^p_x=\textup{d}b_{(t,x)}.(\textup{k}_+\lambda_{px},Id)$, 
so that
$\textup{k}_+f^p_x(y)=\textup{d}b_{(t,x)}((\textup{k}_+\lambda_{px}(y),y)=
ty+\textup{k}_+\lambda_{px}(y)x=
(\sin\log\|x\|_p)y\break
+
\frac{\cos\log\|x\|_p}{\|x\|_p}\, \textup{k}_+N^p_x(y)x$.
This provides the inequalities
$\|\textup{k}_+f^p_x(y)\|_p\leq\break
|\sin\log\|x\|_p|\,\|y\|_p+
\frac{|\cos\log\|x\|_p|}{\|x\|_p}\, |\textup{k}_+N^p_x(y)|\ \|x\|_p\leq
\|y\|_p+|\textup{k}_+N^p_x(y)|$.
We can now use 2.3.17 to obtain the fact that $f^p$ is 2-lipschitzian, once we have proved that, for $p\in\{1,2,\infty\}$, we really have
$|\textup{k}_+N^p_x(y)|\leq\|y\|_p$!\break
Indeed, we just have to use the calculus made in the examples 
of 2.4.7 with $x\not=0$:
if $p=2$, we have $\textup{k}_+N^2_x(y)=$d$N^2_x(y)=\frac{<x,y>}{\|x\|_2}$, so that
$|$k$_+N^2_x(y)|=|\frac{<x,y>}{\|x\|_2}|\leq\frac{\|x\|_2\|y\|_2}{\|x\|_2}
=\|y\|_2$;
if $p=1$, we have k$_+N^1_x(y)=\sum_{i\in{\mathcal I}_0(x)}|y_i|+
\sum_{i\in{\mathcal I}_1(x)}sign(x_i)y_i$, so that
$|$k$_+N^1_x(y)|\leq
\sum_{i\in{\mathcal I}_0(x)}|y_i|+
\sum_{i\in{\mathcal I}_1(x)}|y_i|=\|y\|_1$;
and for $p=\infty$, k$_+N^\infty_x(y)=$
Max$_{i\in\overline{\mathcal I}(|x|)}(sign(x_i)y_i)$, so that
$|$k$_+N^\infty_x(y)|\buildrel\star\over\leq$
Max$_{i\in\overline{\mathcal I}(|x|)}(|sign(x_i)y_i|)=$
Max$_{i\in\overline{\mathcal I}(|x|)}|y_i|
\leq\|y\|_\infty$ (the inequality $\buildrel\star\over\leq$  resulting from the inequality $|$Max(x)$|\leq$Max$|x|$, since Max is 1-lipschitzian).

\includegraphics{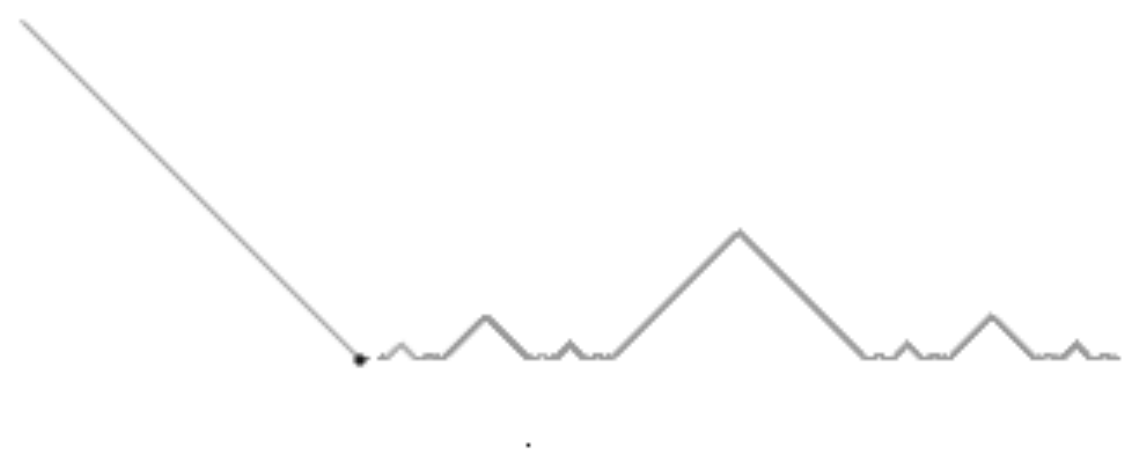}
\begin{picture}(0,0)
\put(-66,2){Figure 4}
\end{picture}

4) $g$ is 1-lipschitzian (as every distance to a non empty subset).
Besides, as $\frac{1}{3}K_\infty=K_\infty$, we have 
$g(\frac{1}{3}x)=d(\frac{1}{3}x,K_\infty)=d(\frac{1}{3}x,\frac{1}{3}K_\infty)
=\frac{1}{3}d(x,K_\infty)=\frac{1}{3}g(x)$.
For a more extensive study of this function, the reader should refer
farther to 2.5.8, where we speak of neo-fractal maps.
We could call $g$ the ``Giseh'' function, 
if, as Napoleon, we gaze at the Giseh pyramides diminishing at the horizon (see Fig 4)! 
\cqfd
\vspace{2mm}

\begin{remks}
We can easily calculate the lipschitzian ratio for the tangentials 
$\textup{T}f_0$ and $\textup{T}g_0$ of the functions $f$ and $g$ studied
in 2.5.2; we compare these lipschitzian ratios respectively to 
$d(\textup{T}f_0,O)$ and $d(\textup{T}g_0,O)$; we precisely obtain 
$1=d(\textup{T}f_0,O)<\rho(\textup{T}f_0)=\sqrt 2$
and $d(\textup{T}g_0,O)=\rho(\textup{T}g_0)=1$. 
Refer to 1.2.8 and 1.2.10... By the way, \textup{T}$f_0$ is the first example
of a jet which is not a good jet (see 1.2.11).

- Indeed, since $\sup_{\theta\in\R}|\sin\theta+\cos\theta|=\sqrt 2$, 
we have $\sup_{x\in\R^*}|f'(x)|=\sqrt 2$, so that $f$ is 
$\sqrt 2$-lipschitzian; thus $\rho(\textup{T}f_0)\leq\sqrt 2$.
Now, using 2.3.14, we know that
$|\frac{f(x)-f(y)}{x-y}|\leq\rho(\textup{T}f_0)$ for all $x,y\in B'(0,1)$ vérifying $x\not= y$. Fixing then $x\not=0$, and doing $x\not=y\rightarrow 0$, we obtain $|f'(x)|\leq\rho(\textup{T}f_0)$ for such $x$; setting
$x=e^{\frac{\pi}{4}+2k\pi}$, with $k\in\Z$ chosen such that 
$x<1$,
we obtain $f'(x)=\sqrt 2$, and thus $\sqrt 2\leq\rho(\textup{T}f_0)$.
So $\rho(\textup{T}f_0)=\sqrt 2$. On the 
other hand, using 2.3.2,
we can write $d(\textup{T}f_0,O)=\|f\|=
\sup\{\frac{|f(x)|}{|x|}\ |\ 0<|x|\leq 1\}=
\sup\{|\sin\log|x||\, |\
0<x\leq 1\}=
\sup\{\sin\theta\,|\, \theta\leq 0\}=1$.
So that, here, $d(\textup{T}f_0,O)<\rho(\textup{T}f_0)$. 

- We have seen that $g$ is 1-lipschitzian, so that
$\rho(\textup{T}g_0)\leq 1$. 
Besides, using 2.5.4 below, $d(\textup{T}g_0,O)=\|g\|=
\sup\{\frac{|g(x)|}{|x|}\ |\ \frac{1}{3}<|x|\leq 1\}\geq 1$
(since, if $x<0$, we have $g(x)=-x$, i.e $\frac{g(x)}{x}=-1$), 
and thus $1\leq d(\textup{T}g_0,O)\leq\rho(\textup{T}g_0)\leq 1$.
\end{remks}
\vspace{2mm}

\begin{prop}
For $h\in r$-$\L\textup{Frac}(E,E')$ where $E\not=\{0\}$, we have:
$\|h\|=\sup\{\frac{\|h(x)\|}{\|x\|}\, |\, r<\|x\|\leq 1\}$
\end{prop}

\proof
Let us put 
$R=\sup\{\frac{\|h(x)\|}{\|x\|}\, |\, r<\|x\|\leq 1\}$.
For every $x\in E$ verifying $r<\|x\|\leq 1$, we have 
$\frac{\|h(x)\|}{\|x\|}\leq\|h\|$, and thus $R\leq\|h\|$ (by 2.3.2).
Conversely, let us consider $x\in E$ verifying $0<\|x\|\leq 1$
and $n\in\N$ such that $r^{n+1}<\|x\|\leq r^n$, i.e
$r<\|r^{-n}x\|\leq 1$. The map $h$ being $r$-Lfractal,
we have $\frac{\|h(x)\|}{\|x\|}=\frac{\|h(r^{-n}x)\|}{\|r^{-n}x\|}\leq R$
(see 2.5.1). Thus $\|h\|\leq R$.
\cqfd

\begin{prop}
Let $h\in r$-$\L\textup{Frac}(E,E')$. Then, for every $\varepsilon>0$, 
we have $\rho(\textup{T}h_0)=\sup\{\frac{\|h(x)-h(y)\|}{\|x-y\|}\, |\, 
x\not= y\ x,y\in C(r,\varepsilon)\}$, with
$C(r,\varepsilon)=\{x\in E\,|\, r<\|x\|<1+\varepsilon\}$. 
\end{prop}

\proof
Let us put $R(\varepsilon)=\sup\{\frac{\|h(x)-h(y)\|}{\|x-y\|}\, |\, 
x\not= y\ x,y\in C(r,\varepsilon)\}$.
Clearly, , $R(\varepsilon)\leq\rho($T$h_0)$ since $h$ is 
$\rho($T$h_0)$-lipschitzian (by 2.3.14). 
Let us now show that $\rho($T$h_0)\leq R(\varepsilon)$.
Let $a\in E$ verifying $0<\|a\|\leq 1$
and $n\in\N$ such that $r^{n+1}<\|a\|\leq r^n$, i.e
$r<\|r^{-n}a\|\leq 1$.
Let us put $\varepsilon'=r^n$Min$\{\varepsilon,r^{-n}\|a\|-r\}$ and let
$x\in B(a,\varepsilon')$.
Then, $\|r^{-n}x\|=r^{-n}\|x\|\leq r^{-n}(\|x-a\|+\|a\|)<r^{-n}\varepsilon'+
r^{-n}\|a\|\leq \varepsilon+1$ and
$\|r^{-n}x\|=r^{-n}\|x\|\geq r^{-n}(\|a\|-\|x-a\|)>r^{-n}\|a\|-r^{-n}\varepsilon'
\geq r^{-n}\|a\|-(r^{-n}\|a\|-r)=r$.
So, $r^{-n}a$, $r^{-n}x\in C(r,\varepsilon)$; then, if $x\not= a$, we have
$R(\varepsilon)\geq\frac{\|h(r^{-n}x)-h(r^{-n}a)\|}{\|r^{-n}x-r^{-n}a\|}=
\frac{\|h(x)-h(a)\|}{\|x-a\|}$, so that $h$ is
$R(\varepsilon)$-$LSL_a$  (see 1.1.7).
This being true for every $a\in B'(0,1)-\{0\}$, we deduce (since $B'(0,1)$ is convex) that the restriction $h|_{B'(0,1)}$ is 
$R(\varepsilon)$-lipschitzian (by 1.1.22)
which finally implies that $\rho($T$h_0)\leq R(\varepsilon)$.
\cqfd

\begin{defi}
Let $(E,U)$, $(E',U')$ be two normed domains (see just before 2.3.6) and $f:(E,U)\lra(E',U')$ a map.
We say that $f$ is\break
$r$-\textup{neo-fractal} at $a\in U$ if $f$
is $\N'_r$-contactable at $a$
(see 2.3.6);
in this case, $\textup{k}_rf_a$ will merely denote 
$\textup{k}_{\N'_r}f_a$.
\end{defi}

\begin{remks}
In each case, $E$ and $E'$ are n.v.s. and $f:U\lra E'$ where $U$ is an open subset of $E$, and $a\in U$. In what follows, we essentially use 2.3.7.

1) When $f$ is G-differentiable at $a$, it is also $r$-\textup{neo-fractal}
for every $0<r<1$ (since $\L$\textup{Hom}$(E,E')
\subset r$-$\L$\textup{Frac}$(E,E')$),
and we have\break
$\textup{k}_rf_a=\textup{k}_+f_a$.

2) Every $r$-\textup{fractal} map $h:E\lra E'$ is $r$-\textup{neo-fractal}
at 0 and we have $\textup{k}_rh_0=h$.

3) Every $r$-\textup{neo-fractal} map at $a$ is tangentiable at $a$.

4) When $f$ is $r$-\textup{neo-fractal} at $a$, then $f$ is 
G-differentiable at $a$ iff $\textup{k}_rf_a\in\L\textup{Hom}(E,E')$;
in this case, $\textup{k}_+f_a=\textup{k}_rf_a$.
\end{remks}

\vspace{3mm}
\begin{exams}
{}
\end{exams}

We consider successively the examples 2),3),4) yet studied in 2.5.2:

1) The function $f$ is $e^{-2\pi}$-neo-fractal at 0, and
differentiable on $\R^*$.

2) For $p\in\{1,2,\infty\}$, the map $f^p$ is $e^{-2\pi}$-neo-fractal at 0, and\break
G-differentiable at every $x\not=0$.

3) The Giseh function $g$ is G-differentiable at every $x\notin K_\infty$.
Furthermore, if we denote $K^+_\infty$ and $K^-_\infty$ the subsets of $K_\infty$ defined, 
for $x\in K_\infty$, by:
\quad$x\in K^+_\infty\Llra \exists\varepsilon>0\ \,(]x-\varepsilon,x[\,
\cap K_\infty=\emptyset)$,\par

\quad\ \,$x\in K^-_\infty\Llra \exists\varepsilon>0\ \,(]x,x+\varepsilon[\,
\cap K_\infty=\emptyset)$,

\noindent then $g$ is $\frac{1}{3}$-neo-fractal at every point of
$K^+_\infty\cup K^-_\infty$, and we have:\break
for $a\in K^+_\infty$,
$\ $k$_{\frac{1}{3}}g_a=g$; and for $a\in K^-_\infty$,
$\ $k$_{\frac{1}{3}}g_a=g_-$, where $g_-(x)=g(-x)$.

\proof
1) Comes from remarks 2.5.7.

2) Still thanks to remarks 2.5.7, using the results of 2.5.2.

3) For the Giseh function, we have to work much harder!
Rememberring that, in the triadic expression of the element
of $\K$, we drop all the 1. 
We put $\widetilde\R=\{s\in\{0,1,2\}^{\Z}\, |\, \exists N\in\Z\ 
\forall n\in\Z\ \,(n\geq N\Lra s_n=0)\}$,  

\quad\quad\quad\ $\widehat\R=\{s\in\widetilde\R\, |\,  
\forall n\in\Z\ \exists m\in\Z\ \,(m\leq n\ \,\textup{and}\ \,s_m\not=2)\}$, 

\quad\quad\quad\ $\widetilde\K=\widetilde\R\cap\{0,2\}^\Z$,

\quad\quad\quad\ $\widehat{\widehat\R}=
\widehat\R-(\widehat\R\cap\widetilde\K)$,

\quad\quad\quad\ $Tri^+=\{s\in\widetilde\K\, |\, \exists N\in\Z\ 
\forall n\in\Z\ \,(n\leq N\Lra s_n=0)\}$,  

\quad\quad\quad\ $Tri^-=\{s\in\widetilde\K\, |\, \exists N\in\Z\ 
\forall n\in\Z\ \,(n\leq N\Lra s_n=2)\}$.

Then, we consider the function $\varphi:\widetilde\R\lra\R$ 
defined by\break
$\varphi(s)=\sum_{n\in\Z}s_n3^n$.
Actually, the restriction of $\varphi$ to $\widehat\R\lra\R_+$
is bijective, and we verify easily that
$\varphi(\widetilde\K)=K_\infty$, $\varphi(Tri^+)=K^+_\infty$
and $\varphi(Tri^-)=K^-_\infty$.
On the other hand, we consider the maps $p,q:\widehat\R\lra
\widetilde\R$ defined by (with, for $s\in\widehat{\widehat\R}$, 
$\ k(s)=\sup\{n\in\Z\,|\, s_n=1\}$):

- if $s\in\widehat\R\cap\widetilde\K$, we set $p(s)=s=q(s)$,

- if $s\not\in\widehat\R\cap\widetilde\K$, we set, for $n\in\Z$,

\quad\quad $p(s)_n=s_n=q(s)_n$, if $n>k(s)$,

\quad\quad $p(s)_n=0$, $q(s)_n=2$, , if $n=k(s)$,

\quad\quad $p(s)_n=2$, $q(s)_n=0$, if $n<k(s)$.

Then, we verify that $p(\widehat{\widehat\R})\subset Tri^-$, 
$q(\widehat{\widehat\R})\subset Tri^+$, and that, for all $s\in\widehat\R$, we have
$\varphi.p(s)\leq\varphi(s)$ and $\varphi.q(s)\geq \varphi(s)$; and
especially that $g(\varphi(s))=d(\varphi(s),K_\infty)=$
Min$(\varphi(s)-\varphi.p(s),\varphi.q(s)-\varphi(s))$.

When $x\in\,]\!-\infty,0[$, clearly $g(x)=-x$ (since $K_\infty\subset\R_+$ and $0\in K_\infty$), so that $g$ is G-differentiable at $x$.
When $x\in\R_+-K_\infty$, then (if, we set $x^+=\varphi(q(s))$ and $x^-=\varphi(p(s))$, where $s\in\widehat\R$ verifies $\varphi(s)=x$) we have
$]x^-,x^+[\,\subset\R_+-K_\infty$. Actually, for all $y\in\,]x^-,x^+[$, we have $y^-=x^-$, $y^+=x^+$  and $g(y)=$Min$(|y-x^-|,|y-x^+|)$,
so that the restriction $g|_{]x^-,x^+[}$ is G-differentiable.
Thus $g$ is G-differentiable at $x$.

On the other hand, for $i\in\Z$, we consider the maps
$\mu_i:\widetilde\R\lra\widetilde\R$ defined by $\mu_i(s)_n=s_{n+i}$,
and $m_i:\R\lra\R$ defined by $m_i(x)=3^{-i}x$; then, we have the 
following commutative diagrams:
$$\xymatrix{
\widetilde\R
\ar[d]_{\varphi}
\ \ar[r]^>>>>>{\mu_i}&\ 
\widetilde\R
\ar[d]^{\varphi}&&\widehat\R
\ar[d]_{p}
\ \ar[r]^>>>>>{\mu_i}&\ 
\widehat\R
\ar[d]^{p}&&\widehat\R
\ar[d]_{q}
\ \ar[r]^>>>>>{\mu_i}&\ 
\widehat\R
\ar[d]^{q}
\\
\R\ 
\ar[r]_>>>>>{m_i}&\ \R&&
\widetilde\K 
\ar[r]_>>>>>{\mu_i}&\ \widetilde\K&&\widetilde\K 
\ar[r]_>>>>>{\mu_i}&\ \widetilde\K
}$$

$$\xymatrix{
\widehat{\widehat\R}
\ar[d]_{k}
\ \ar[rr]^>>>>>>>>>>{\mu_i}&&\ 
\widehat{\widehat\R}
\ar[d]^{k}\\
\Z\ 
\ar[rr]_>>>>>>>>>>{(-)-i}&&\ \Z
}$$

This provides that $g$ is $\frac{1}{3}$-Lfractal (yet noticed in 
2.5.2). Now: 

 - If $s\in Tri^+$ and $N\in\Z$ are such that, for all $n\in\Z$,
we have $n\leq N\Lra s_n=0$, then, for $a=\varphi(s)$, we have
$g(x)=g(x-a)$ for all $x\in\,]a-\frac{3^{N+1}}{2},a+\frac{3^{N+1}}{2}[$;
so $g$ is $\frac{1}{3}$-neo-fractal and k$_{\frac{1}{3}}g_a=g$.

 - If $s\in Tri^-$ and $N\in\Z$ are such that, for all $n\in\Z$,
we have $n\leq N\Lra s_n=2$, then, for $a=\varphi(s)$, we have
$g(x)=g(a-x)$ for all $x\in\,]a-\frac{3^{N+1}}{2},a+\frac{3^{N+1}}{2}[$;
so $g$ is $\frac{1}{3}$-neo-fractal and k$_{\frac{1}{3}}g_a=g_-$
(actually, if $s\in Tri^+$ and $x\in\,]a-\frac{3^{N+1}}{2},a[$, then
$g(x)=a-x$; and if $s\in Tri^-$ and $x\in\,]a,a+\frac{3^{N+1}}{2}[$, then
$g(x)=x-a$).
\cqfd

\begin{remks}
In the previous examples, we notice that:\par

1) a) $f$ is $e^{-2\pi}$-neo-fractal at 0, but not G-differentiable at 0 
(we use 4) in 2.5.7, since $\textup{k}_r f_0$ is clearly not 
$\R_+$-homogeneous, by 2.4.3).

\quad$\,$b) Same remark for the $f^p$ where $p\in\{1,2,\infty\}$.

\quad$\,$c) $g$ is not G-differentiable at all $x\in K^+_\infty\cup K^-_\infty$
(since $g$ and $g_-$ are not $\R_+$-homogeneous), although it is
$\frac{1}{3}$-neo-fractal at these points.

2) a) Of course, there exist neo-fractal maps which are not 
Lfractal: we have just, as in 2.4.8, to translate our previous examples
at every point where they are neo-fractal.

\quad$\,$b) Concerning the function $f$ of 2.5.8, translate it is 
perhaps not the best idea since, although being no more
Lfractal, it remains differentiable at 0.
Having said this, we obtain a convincing example considering
the function $x^2+f(x)$ ... The reader, so guided, will find a lot of
other good examples of neo-fractal maps which are not Lfractal.
\end{remks}

\vspace{5mm}
\centerline{\texttt{Construction of fractal functions}}
\vspace{1mm}

Let $s$ and $T$ be strictly positive real numbers and
$f:\R\lra\R$ a $T$-periodic and $s$-lipschitzian function which admits a right derivative at every point (we know that, then, $|f'_r(x)|\leq s$ for all $x\in\R$). In particular, $f$ is bounded on $\R$ (since it is so on $[0,T]$).
Then, we associate to $f$ the function $\varphi:\R\lra\R$ defined by 
$\varphi(0)=0$ and, for $x\not=0$, $\varphi(x)=xf(\log|x|)$.
Then $\varphi$ admits a right derivative at every points of $]0,+\infty[$
(since $\log$ is strictly increasing and derivable) and 
$\varphi'_r(x)=\break
f(\log(x))+f'_r(\log(x))$;
so that we have $|\varphi'_r(x)|\leq R+s$ where\break
$R=\sup_{x\in\R}f(x)$. The function
$\varphi$ is thus also lipschitzian on $\R_+$ (it is continuous at 0
since it is $SL_0$); it is even lipschitzian on $\R$
(for, if $\rho=R+s$ and $x<0<y$, we have
$|\varphi(y)-\varphi(x)|\leq
|\varphi(y)-\varphi(0)|+\break
|\varphi(0)-\varphi(x)|\leq\rho|y|+\rho|x|=
\rho y-\rho x=\rho(y-x)=\rho|y-x|$).
Besides, if we set $r=e^{-T}$, we have $0<r<1$ and, for $x\not=0$,
$\varphi(rx)=rxf(\log|rx|)=rxf(\log r+\log|x|)=rxf(-T+\log|x|)=rxf(\log|x|)=
r\varphi(x)$, so that $\varphi$ is $r$-Lfractal.

Let us consider now the set ${\mathcal P}_T$ of the $T$-periodic functions
$\R\lra\R$ which are lipschitzian and which admit a right derivative at every point. Then, 
${\mathcal P}_T$ has a structure of vectorial subspace of $\R^\R$. The previous construction provides a map
$j:{\mathcal P}_T\lra r$-$\L\textup{Frac}(\R,\R)$, where $j(f)$ is the function $\varphi$ associated to $f$ as above.
This map is clearly linear; and injec-\break
tive for, 
$j(f)=0\Lra
\forall x\in\R^*
\ f(\log|x|)=0\Lra f=0$.
By composition, we have an injective linear map (referring to 2.2.4 for $J$ and to 1.4.15 for the isometry $can$) :

\qquad
${\mathcal P}_T\buildrel j\over\lra r$-$\F\textup{rac}(\R,\R)
\buildrel J\over\lra \J\textup{et}((\R,0),(\R,0))
\buildrel{can}\over\lra\J\textup{et}\jf(\R,\R)$

\begin{prop}
The space $\J\textup{et}\jf(\R,\R)$ is a vectorial space of\break
infinite dimension.
\end{prop}

\proof
We use the previous embedding, knowing that ${\mathcal P}_T$
is a vectorial space of infinite dimension (since the familly
$(f_n)_{n\in\N}$ is free,\break
with $f_n(x)=\sin(\frac{2\pi n}{T})x$).
\cqfd

\vspace{5mm}
\centerline{\texttt{Summary}}

The tangentiable maps have made brought out  new classes of maps.
Let us give here a recapitulative diagram of the various implications
proved all along this paper:
$$\xymatrix{
{}&\R\!-\!cont_a\ar@{=>}[dr]^{2''}&{}\\
Dif\! f_a\ar@{=>}[ru]^{2'}
\ar@{=>}[rr]^{2}&& G\!-\! di\! f\! f_a\ar@{=>}[r]^{3} & r\!-\! neo\! fr_a
\ar@{=>}[d]^{4}\\
C^1\ar@{=>}[rr]_{5}\ar@{=>}[u]^{1}&& LL_a\ar@{=>}[r]_{6}
 & Tang_a\ar@{=>}[r]_{7} & LSL_a
\ar@{=>}[r]_{8} & C^0_a
}$$
Where, here, $Di\! f\! f_a$, $G$-$di\!f\!f_a$ and $r$-$neo\! fr$ and
$\R\!-\! cont_a$ stand respectively for differentiable, G-differentiable , $r$-neo-fractal and standard $\R$-contactable (refer to the end of section 2.3) ... at $a$ for all of them
($C^1$ means ``of class $C^1\,$'', and
$C^0_a$ means continuous at $a$).

One can find the different proofs of these implications in:
2.2.17 (for 2'), 2.4.6 (for 2 and 2''), 2.5.7 (for 3 and 4), 1.1.12 (for 5),
1.3.2 (for 6 and 7) and 1.1.8 (for 8).

In the above diagram, every inverse implication is false; we give counter-examples below.

\vspace{2mm}
\begin{cexams} 
{}
\end{cexams}
In each case (except for 2'), we denote $f:\R\lra\R$ each given counter-example (here $a=0$ and $f(0)=0$; the ``number'' $i)$ corresponding to a counter-example to the $i^{ieth}$ above implication).

1) $f(x)=x^2\sin{\frac{1}{x}}$ (well-known),

2) $f(x)=|x|$ (see examples 2.4.7 and remark 2.4.8),

2') see 2.5.12,

2'') same as for 2) above,

3) $f(x)=x\sin(\log|x|)$ (see example 2.5.8 and remark 2.5.9),

4) $f(x)=x\sin(\log|\log|x||)$ if $x\not=0$ 
(will be studied later).

5) $f(x)=|x|$ (lipschitzian but not 
differentiable at 0).

6) $f(x)=x^2\sin\frac{1}{x^2}$ (see 1.3.9),

7) $f(x)=x\sin\frac{1}{x}$ (same as for 6)),

8) $f(x)=x^{\frac{1}{3}}$ (same as for 6)).

\begin{remks}\par\hfill

\includegraphics{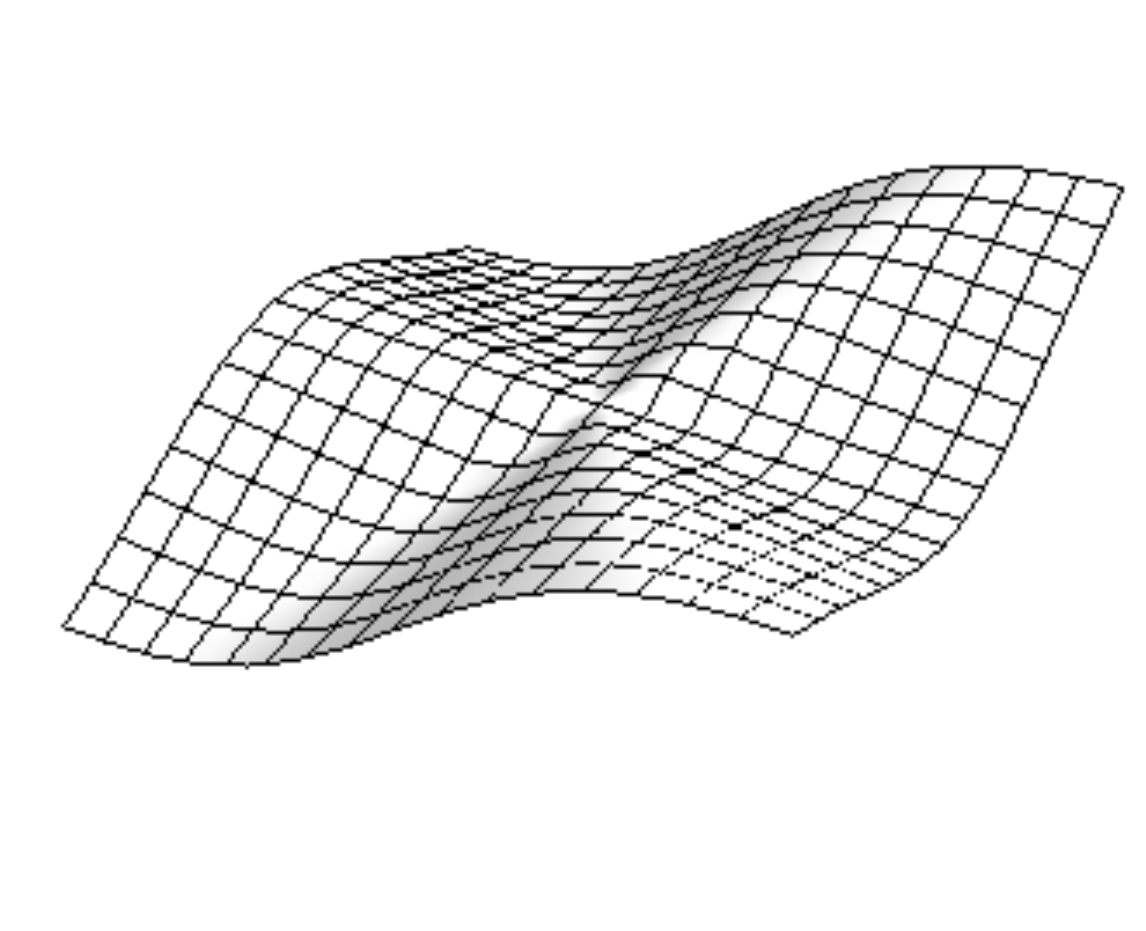}
\begin{picture}(0,0)
\put(-66,12){Figure 5}
\end{picture}

1) We prove here that $\R\!-\! cont_a\ \not\!\!\Lra Di\!f\!f_a$:
consider the function $f:\R^2\lra\R$ defined by $f(0,0)=0$ and $f(x,y)=
\frac{xy^2}{x^2+y^2}$ if $(x,y)\not=(0,0)$. This $f$ is
obviously 
differentiable on $\R^2-\{(0,0)\}$, with, for $(x,y)\not=(0,0)$,
$\frac{\partial f}{\partial x}(x,y)
=\frac{y^2(y^2-x^2)}{(x^2+y^2)^2}$ and
$\frac{\partial f}{\partial y}(x,y)=\frac{2x^3y}{(x^2+y^2)^2}$, so that
$|\frac{\partial f}{\partial x}(x,y)|\leq\frac{y^2(y^2+x^2)}{(x^2+y^2)^2}=\frac{y^2}{x^2+y^2}\leq 1$ and
$|\frac{\partial f}{\partial y}(x,y)|\leq 2$ (since $|x|,|y|\leq\sqrt{x^2+y^2)}$);
thus $f$ is lipschitzian on $\R^2$.
Furthermore, it is obviously $\R$-homogeneous, i.e, it verifies
$f(tx,ty)=tf(x,y)$ for all $t\in\R$, which implies that $f$ is
standard $\R$-contactable at 0, with $\textup{k}_\R f_0=f$. However, since $f$ is not linear, it cannot be differentiable at 0 (see 2.4.6 and Fig 5)!

2) We can complete the above diagram of implications,
by adding the following diagram (where $\R$-Lhom, $\R_+$\!-\! Lhom and 
$r$-Lfrac stand respectively for standard $\R$-$\!$Lhomogeneous, 
$\R_+$-Lhomogeneous and
$r$-Lfractal):

$$\xymatrix{
\R\!-\! Lhom\ar@{=>}[r]^{2'''}\ar@{=>}[d]_{9'}&
\R_+\!-\! Lhom\ar@{=>}[r]^{3'}\ar@{=>}[d]_9&r\!-\! Lfrac
\ar@{=>}[d]^{10}\\
\R\!-\! cont_0\ar@{=>}[r]_{2''}&
G\!-\! di\! f\! f_0\ar@{=>}[r]_3&r\!-\! neo\! fr_0
}$$

The proofs of these implications are in 2.5.7 (for 3), 2.5.2 (for 3'), just before 2.4.1
(for $2'''$), 2.4.6 (for 2''), 
and 2.3.7 (for 9 and 10); as for 9', refer the end of section 2.3.

Here again the inverse implications are false:
3' for the same reason as 3;
$2'''$ for the same reason as 2'';
for 9, refer to 2) in 2.4.8;
for 10, refer to 2) in 2.5.9; as for 9', we use the same arguments as for 
9 or 10.
\end{remks}


\section{Local extrema}
In this last section, we present nice generalisations of classical theorems about extrema of functions with codomain $\R$. In particular, we give a sufficient condition for having an extremum
which only needs hypotheses at order 1!

\begin{theo}
Let $\Sigma$ be a valued monoid, $b$ a fixed real number, $(M,U)$ a $\Sigma$-contracting
domain (where $U$ is a neigborhood oh $\omega$ in $M$)
and $f:(M,U)\lra(\R_b,\R)$ a centred function
(see just before 2.2.14). If $f$ is $\Sigma$-contactable and if $f$ admits a local minimum at $\omega$, then $\textup{K}_\Sigma f$ admits a global minimum at $\omega$.
\end{theo}

\proof
We recall that $\R_b$ is a canonical revertible $\Sigma$-contracting space whose external operation is given, for $t\in\Sigma$ and $y\in\R$, by $t\star y=v(t)(y-b)+b$ which verifies, for $t\in\Sigma^*$, 
$\ t\buildrel{-1}\over\star y=\frac{y-b}{v(t)}+b$ (see 2.1.10 and 2.1.14).
If $h=K_\Sigma f$; this $h$ verifying by definition $f\tang_\omega h|_U$, we can use 2.1.20 to obtain
$h(x)=
\lim_{0\not=v(t)\rightarrow 0}t\buildrel{-1}\over\star f(t\star x)=
\lim_{0\not=v(t)\rightarrow 0}\frac{f(t\star x)-b}{v(t)}+b$, for every $x\in M$.
Since $f$ is supposed to admit a local minimum at $\omega$, there exists a neighborhood $V$
of $\omega$ in $U$ such that $f(x)\geq f(\omega)=b$ for all $x\in V$.
Fixing then $x\in M$, and using the fact that $\lim_{v(t)\rightarrow 0}t\star x
=\omega$, there exists $\varepsilon>0$ such that, for all $t\in\Sigma$, 
we have the implication: $0<v(t)<\varepsilon\Lra t\star x\in V$.
So, when $0<v(t)<\varepsilon$, we have $f(t\star x)\geq b$ which implies
$t\buildrel{-1}\over\star f(t\star x)=\frac{f(t\star x)-b}{v(t)}+b\geq b$. Doing $v(t)\rightarrow 0$,
we obtain $h(x)\geq b=h(\omega)$, and this for all $x\in M$; we have so proved that $h=$K$_\Sigma f$ admits a global minimum at $\omega$.
\cqfd

\begin{cory}
Let $(E,U)$ a normed domain (see just before 2.3.6) and $f:(E,U)\lra(\R,\R)$
a $\Sigma$-contactable function at a point $a\in U$, and which admits a local minimum at $a$.
Then $\textup{k}_\Sigma f_a$ admits a global minimum at 0.
\end{cory}

\proof
It, comes from the fact that, for all $x\in E$, we have
k$_\Sigma f_a(x)=$\break
K$_\Sigma f(x+a)-f(a)\geq 0=$k$_\Sigma f_a(0)$, where K$_\Sigma f:E_a\lra\R_{f(a)}$. 
\cqfd

\begin{remk}
This gives back the well-known result of the differentiable case:
`` $f$ admits a local minimum at $a\Lra a$ is a critical point of $f$ (i.e $\textup{d}f_a=0$)'',
since there exists a unique linear function $E\lra\R$ which admits a global minimum at 0: the null function. 
\end{remk}

\begin{theo}
Let $(M,U)$ a $\R_+$-contracting domain, $b$ a fixed real number and $f:(M,U)\lra(\R_b,\R)$ a  
$\R_+$-contactable centred function. If $M$ is a Daniel space (i.e a metric space in 
which every closed and bounded subset is compact) and if 
$\textup{K}_\Sigma f:M\lra\R_b$ admits a 
strict global minimum at $\omega$, then $f$ admits a strict local minimum at $\omega$.
\end{theo}

\proof
The case where $M=\{\omega\}$ being immediate, we can suppose $M\not=\{\omega\}$.
Let us set $S=\{x\in M\, |\, d(x,\omega)=1\}$. Then $S$ is a non empty compact
(it is obviously closed and bounded, and if $x\in M-\{\omega\}$, we have 
$\frac{1}{d(x,\omega)}\star x\in S$). Since $h=$K$_\Sigma f$ is continuous (it is lipschitzian), $h$ reaches its inferior bound on $S$: there exists $x_0\in S$ such that $\inf_{x\in S}h(x)=h(x_0)$, so that
$h(x)\geq h(x_0)>b$ for all $x\in S$ (since $x_0\not=\omega$). Consider $\varepsilon=h(x_0)-b>0$. Since $f\tang_\omega h|_U$, there exists $\eta>0$ such that $B(\omega,\eta)\subset U$ and
verifying
the implication: $0<d(x,\omega)<\eta\Lra|f(x)-h(x)|<\varepsilon d(x,\omega)$ for all $x\in M$.
Let us fix $x\in B(\omega,\eta)-\{\omega\}$; it verifies $f(x)>h(x)-\varepsilon d(x,\omega)$.
If $y=\frac{1}{d(x,\omega)}\star x$, we have $y\in S$, so that  $h(y)\geq h(x_0)$ which implies\break
$h(y)-b-\varepsilon\geq h(x_0)-b-\varepsilon=0$. Hence 
(since $h:M\lra\R_b$ is\break
$\R_+$-homogeneous, where $\R_+$ is a quasi-group)
$\ h(x)-\varepsilon d(x,\omega)=\break
h(d(x,\omega)\star y)-\varepsilon d(x,\omega)=
d(x,\omega)(h(y)-b)+b-\varepsilon d(x,\omega)=\break
d(x,\omega)(h(y)-b-\varepsilon)+b\geq b$.
Thus, for all $x\in B(\omega,\eta)-\{\omega\}$, we have $f(x)>h(x)-\varepsilon d(x,\omega)\geq b=
f(\omega)$, which implies that $f$ admits a strict local minimum at $\omega$.
\cqfd

\begin{cory}
Let $(E,U)$ be a normed domain where $E$ is of finite dimension,
$a\in U$ and $f:(E,U)\lra(\R,\R)$ a map G-differentiable at $a$ such that $\textup{k}f_a>0$ (i.e verifying $\textup{k}f_a(x)>0$ for every $x\in E-\{a\}$).
Then, $f$ admits a strict local minimum at $a$.
\end{cory}

\proof
If $h=$K$_+f:E_a\lra\R_{f(a)}$, we have, for all $x\in E-\{a\}$,\break 
$h(x)=f(a)+$k$_+f_a(x-a)>f(a)$. Hence the wished conclusion
(since $E$ is a Daniel space).
\cqfd

\begin{remk}
This theorem has not its equivalent, at order 1, in differential calculus, since a linear function cannot have a strict minimum. It is rather inspired by theorems giving sufficient conditions, at order 2, 
for the existence of extrema
\end{remk}

\vspace{10mm}
\centerline{\small{BIBLIOGRAPHY}}
\vspace{3mm}

\noindent$[1]$ \small{F.Borceux}, \small\textit{Handbook of Categorical Algebra}, \small{3 vol., Cambridge University Press (1994)}.

\noindent$[2]$ \small{E.Burroni and J.Penon, 
Arr\^et sur la tangentialit\' e, conf\' erence au 
\small\textit{SIC(2006) \` a Calais}.

\noindent$[3]$ \small{E.Burroni and J.Penon, 
Un calcul diff\' erentiel m\' etrique, conf\' erence au}
\small\textit{CT(2008) \` a Calais}.

\noindent$[4]$ \small{E.Burroni and J.Penon, 
Repr\' esentation des jets m\' etriques, conf\' erence au}
\small\textit{SIC(2007) \` a Paris}.

\noindent$[5]$ \small{E.Burroni and J.Penon, 
Jets m\' etriques repr\' esentables, conf\' erence au}\break
\small\textit{CT(2008) \` a Calais}.

\noindent$[6]$ \small{C.Ehresmann,
Les prolongements d'une vari\' et\' e diff\' erentiable, I à V},
\small\textit{CRAS Paris} \small{233 (1951) 598-600, 777-779, 1081-1083, \small\textit{CRAS Paris} \small{234 (1952) 1028-1030, 1424-1425}.

\noindent$[7]$ \small{R.Gateaux (we give here three references, the first one historical, the third one more in accordance with the spirit of our article):
  
$\qquad$ - R.Gateaux, Sur la repr\' esentation des fonctionnelles continues},\break
$\qquad\ $\small\textit{Atti della reale accademia dei Lincei} \small{XXIII, 1,
310-345, (1914)}.

$\qquad$ - P.Levy, Leçons d'analyse fonctionnelle, 
\small\textit{GAUTHIER-VILLARS},\break
\small{chap. IV, p. 69, (1922)}. 

$\qquad$ - G.Bouligand, G\' eom\' etrie infinitésimale directe,
\small\textit{VUIBERT},\break
\small{p. 69, (1932)}. 

\vspace {5mm}
\centerline{\texttt{University address of the authors}}

\noindent Institut de Math\'ematiques de Jussieu;
Universit\'e Paris Diderot, Paris 7;\\
5, rue Thomas Mann;
75205 Paris cedex 13; France.\\

\vspace {1mm}
\centerline{\texttt{email addresses}}

\noindent {eburroni@math.jussieu.fr}\\
{penon@math.jussieu.fr}